\documentclass[11pt]{article}

\usepackage{amsmath,amsthm,amssymb,amsxtra,amsfonts,amsbsy,mathrsfs}
\usepackage[all,2cell]{xy} \UseAllTwocells
\SilentMatrices
\newtheorem{theorem}[subsection]{Theorem}
\newtheorem{corollary}[subsection]{Corollary}
\newtheorem{lemma}[subsection]{Lemma}
\newtheorem{proposition}[subsection]{Proposition}
\newtheorem{definition}[subsection]{Definition}
\newtheorem{sublemma}[subsubsection]{Lemma}
\newtheorem{subproposition}[subsubsection]{Proposition}
\theoremstyle{definition}
\newtheorem{notation-convention}[subsection]{Notations and
Conventions}
\newtheorem{notation}[subsubsection]{Notation}

\newtheorem{remark}[subsection]{Remark}
\newtheorem{example}[subsection]{Example}
\newtheorem{problem}[subsection]{Problem}
\newtheorem{anitem}[subsubsection]{}
\newtheorem{blank}[subsection]{}
\newtheorem{subremark}[subsubsection]{Remark}
\numberwithin{equation}{subsection}

\newcommand{\s}{\mathscr}
\newcommand{\bb}{\mathbb}
\newcommand{\fk}{\mathfrak}

\begin{document}

\title{$L$-SERIES OF ARTIN STACKS OVER FINITE FIELDS}
\author{Shenghao Sun}
\date{}
\maketitle

\begin{abstract}
We develop the notion of stratifiability in the context
of derived categories and the six operations for stacks in
\cite{LO1, LO2}. Then we reprove the Lefschetz trace formula
for stacks, and give the meromorphic continuation of
$L$-series (in particular, zeta functions) of $\mathbb
F_q$-stacks. We also give an upper bound for the weights of
the cohomology groups of stacks, and an ``independence of
$\ell"$ result for a certain class of quotient stacks.
\end{abstract}

\section{Introduction}

In topology there is the famous Lefschetz-Hopf trace formula,
which roughly says that if $f:X\to X$ is an endomorphism of a
compact connected oriented space $X$ with isolated fixed
points, then the number of fixed points of $f,$ counted with
multiplicity, is equal to the alternating sum of the traces of
$f^*$ on the singular cohomology groups $H^i(X,\mathbb Q).$

There is also a trace formula in algebraic geometry, for schemes
over finite fields, due to Grothendieck. It says that if $X_0$
is a scheme over $\bb F_q,$ separated and of finite type,
and $F_q$ is the $q$-geometric Frobenius map, then
$$
\#X_0(\bb F_q)=\sum_{i=0}^{2\dim X_0}(-1)^i\text{Tr}(F_q,
H^i_c(X,\overline{\bb Q}_{\ell})),
$$
where $H^i_c(X,\overline{\bb Q}_{\ell})$ is the $\ell$-adic
cohomology with compact support. In fact he proved the trace
formula for an arbitrary constructible sheaf. See \cite{Gro,
Ver, SGA4.5}.

Behrend conjectured the trace formula for smooth algebraic
stacks over $\mathbb F_q$ in his thesis and \cite{Beh1}, and
proved it in \cite{Beh2}. However, he used ordinary cohomology
and arithmetic Frobenius (rather than compact support
cohomology and geometric Frobenius) to prove the ``dual
statement", probably because at that time the theory of
dualizing complexes of algebraic stacks, as well as compact
support cohomology groups of stacks, were not developed.
Later Laszlo and Olsson developed the theory of the six
operations for algebraic stacks \cite{LO1, LO2}, which makes
it possible to reprove the trace formula, and remove the
smoothness assumption in Behrend's result. Also we will 
work with a fixed isomorphism of fields 
$\iota:\overline{\bb Q}_{\ell}\overset{\sim}{\to}
\bb C,$ namely we will work with \textit{$\iota$-mixed} 
complexes, rather than \textit{mixed} ones, and this 
is a more general setting (see (\ref{Laff})).

Once we have the trace formula, we get a factorization of the
zeta function into a possibly infinite product of $L$-factors,
and from this one can deduce the meromorphic continuation of
the zeta functions, generalizing a result of Behrend
(\cite{Beh1}, 3.2.4). Furthermore, to locate the zeros and
poles of the zeta functions, we give a result on the
weights of cohomology groups of stacks.

We briefly mention the technical issues. As pointed out in
\cite{Beh2}, a big difference between schemes and stacks is
the following. If $f:X_0\to Y_0$ is a morphism of
$\bb F_q$-schemes of finite type, and $K_0\in
D_c^b(X_0,\overline{\bb Q}_{\ell}),$ then $f_*K_0$ and
$f_!K_0$ are also bounded complexes. Since often we are
mainly interested in the simplest case when $K_0$ is a sheaf
concentrated in degree 0 (for instance, the constant sheaf
$\overline{\bb Q}_{\ell}),$ and $D^b_c$ is stable under
$f_*$ and $f_!,$ it is enough to consider $D^b_c$ only. But
for a morphism $f:\s X_0\to\s Y_0$ of $\bb
F_q$-algebraic stacks of finite type, $f_*$ and $f_!$ do not
necessarily preserve boundedness. For instance, the
cohomology ring $H^*(B\bb G_m,\overline{\bb
Q}_{\ell})$ is the polynomial ring $\overline{\bb
Q}_{\ell}[T]$ with $\deg(T)=2.$ So for stacks we have to
consider unbounded complexes, even if we are only interested
in the constant sheaf $\overline{\bb Q}_{\ell}.$ In order
to define the trace of the Frobenius on cohomology groups,
we need to consider the convergence of the complex series of
the traces. This leads to the notion of an $\iota$-convergent
complex of sheaves (see (\ref{D4.1})).

Another issue is the following. In the scheme case one
considers bounded complexes, and for any bounded complex
$K_0$ on a scheme $X_0,$ there exists a stratification of
$X_0$ that ``trivializes the complex $K_0"$ (i.e. the
restrictions of all cohomology sheaves $\s H^iK_0$ to
each stratum are lisse). But in the stack case we have to
consider unbounded complexes, and in general there might be
no stratification of the stack that trivializes every
cohomology sheaf. This leads to the notion of a stratifiable
complex of sheaves (see (\ref{D3.1})). We need the
stratifiability condition to control the dimensions of
cohomology groups (see (\ref{L3.11})). All bounded
complexes are stratifiable (\ref{L3.3}v).

Also we will have to impose the condition of 
$\iota$-mixedness, due to unboundedness. For bounded 
complexes on schemes, the trace formula can be proved 
without using this assumption, although the conjecture of 
Deligne (\cite{Del2}, 1.2.9) that all sheaves are 
$\iota$-mixed is proved by Laurent Lafforgue. See the 
remark (\ref{Laff}).

We briefly introduce the main results of this paper.

\begin{flushleft}
\textbf{Fixed point formula.}
\end{flushleft}

\begin{theorem}\label{T1.1}
Let $\s X_0$ be an Artin stack of finite type over
$\bb F_q,$ and let $[\s X_0(\bb F_q)]$ be the set of 
isomorphism classes of the groupoid of $\bb F_q$-points 
of $\s X_0.$ Then the series
$$
\sum_{n\in\bb Z}(-1)^n\emph{Tr}(F_q,H^n_c(\s
X,\overline{\bb Q}_{\ell})),
$$
regarded as a complex series via $\iota,$ is absolutely
convergent, and its limit is ``the number of $\mathbb
F_q$-points of $\mathscr X_0",$ namely
$$
\#\s X_0(\bb F_q):=\sum_{x\in[\s X_0(\bb
F_q)]}\frac{1}{\#\emph{Aut}_x(\bb F_q)}.
$$
\end{theorem}

Here $F_q$ denotes the $q$-geometric Frobenius. To
generalize, one wants to impose some condition (P) on
complexes $K_0\in D_c^-(\s X_0,\overline{\bb
Q}_{\ell})$ such that

(1) (P) is preserved by $f_!;$

(2) if a complex $K_0$ satisfies (P), then the ``naive local
terms" are well-defined, and

(3) trace formula holds in this case.

The condition (P) on $K_0$ turns out to be a combination of
three parts: $\iota$-convergence (which implies (2) for
$K_0$), $\iota$-mixedness and stratifiability (which, together
with the first part, implies (2) for $f_!K_0$). See
(\ref{T4.3}) for the general statement of the theorem. These
conditions are due to Behrend \cite{Beh2}.

\vskip.5truecm

\begin{flushleft}
\textbf{Meromorphic continuation.}
\end{flushleft}

The rationality in Weil conjecture was first proved by Dwork,
namely the zeta function $Z(X_0,t)$ of every variety
$X_0$ over $\bb F_q$ is a rational function in $t.$
Later, this was reproved using the fixed point formula
\cite{Gro, SGA5}. Following (\cite{Beh1}, 3.2.3), we define
the zeta functions of stacks as follows.

\begin{definition}\label{D1.2}
For an $\bb F_q$-algebraic stack $\s X_0$ of finite
type, define the \emph{zeta function}
$$
Z(\s X_0,t)=\exp\Big(\sum_{v\ge1}\frac{t^v}{v}\sum
_{x\in[\s X_0(\bb 
F_{q^v})]}\frac{1}{\#\emph{Aut}_x(\bb F_{q^v})}\Big),
$$
as a formal power series in the variable $t.$ 
\end{definition}

Notice that in general, the zeta function is not rational
(cf. $\S7$). Behrend proved that, if $\s X_0$ is a
smooth algebraic stack, and it is a quotient of an
algebraic space by a linear algebraic group, then its zeta
function $Z(\s X_0,t)$ is a meromorphic function in
the complex $t$-plane; if $\s X_0$ is a smooth
Deligne-Mumford stack, then $Z(\s X_0,t)$ is a
rational function (\cite{Beh1}, 3.2.4, 3.2.5). These
results can be generalized as follows.

\begin{theorem}\label{T1.3}
For every $\bb F_q$-algebraic stack $\s X_0$ of
finite type, its zeta function $Z(\s X_0,t)$ defines a
meromorphic function in the whole complex $t$-plane. If
$\s X_0$ is Deligne-Mumford, then $Z(\s X_0,t)$
is a rational function.
\end{theorem}

See (\ref{C7.2}) and (\ref{T8.1}) for the general statement.

\vskip.5truecm

\begin{flushleft}
\textbf{A theorem of weights.}
\end{flushleft}

One of the main results in \cite{Del2} is that, if $X_0$ is
an $\bb F_q$-scheme, separated and of finite type, and
$\s F_0$ is an $\iota$-mixed sheaf on $X_0$ of punctual
$\iota$-weights $\le w\in\bb R,$ then for every $n,$ the
punctual $\iota$-weights of $H^n_c(X,\s F)$ are $\le
w+n.$ The cohomology groups are zero unless $0\le n\le2\dim
X_0.$ We will see later (\ref{R7.1}) that the upper
bound $w+n$ for the punctual $\iota$-weights does not work in
general for algebraic stacks. We will give an upper bound
that applies to all algebraic stacks. Deligne's upper bound
of weights still applies to stacks for which all the
automorphism groups are affine.

\begin{theorem}\label{T1.4}
Let $\s X_0$ be an $\bb F_q$-algebraic stack of
finite type, and let $\s F_0$ be an $\iota$-mixed sheaf
on $\s X_0,$ with punctual $\iota$-weights $\le w,$ for
some $w\in\bb R.$ Then the $\iota$-weights of
$H^n_c(\s{X,F})$ are $\le\dim\s 
X_0+\frac{n}{2}+w,$ and they are congruent mod $\bb Z$
to weights that appear in $\s F_0.$ If
$n>2\dim\s X_0,$ then $H^n_c(\s X,-)=0$ on
sheaves. If for all points $\overline{x}\in\s X(\bb 
F)$ in the support of $\s F,$ the automorphism group 
schemes $\emph{Aut}_{\overline{x}}$ are affine, then the 
$\iota$-weights of $H^n_c(\s{X,F})$ are $\le n+w.$
\end{theorem}

\textbf{Organization.} In $\S2$ we review some preliminaries
on derived categories of $\ell$-adic sheaves on algebraic
stacks over $\bb F_q$ and $\iota$-mixed complexes, and
show that $\iota$-mixedness is stable under the six
operations.

In $\S3$ we develop the notion of stratifiable complexes in
the context of Laszlo and Olsson's $\ell$-adic derived
categories, and prove its stability under the six operations.

In $\S4$ we discuss convergent complexes, and show that they
are preserved by $f_!.$ In $\S5$ we prove the trace formula
for stacks. These two theorems are stated and proved in
\cite{Beh2} in terms of ordinary cohomology and arithmetic
Frobenius, and the proof we give here uses geometric
Frobenius.

In $\S6$ we discuss convergence of infinite products of
formal power series, which will be used in the proof of the
meromorphic continuation. In $\S7$ we give some examples of
zeta functions of stacks, and give the functional equation of
the zeta functions and independence of $\ell$ of Frobenius
eigenvalues for proper varieties with quotient singularities
(\ref{T7.3}).

In $\S8$ and $\S9,$ we prove the meromorphic continuation and
the weight theorem respectively. Finally in $\S10$ we discuss
``independence of $\ell"$ for stacks, and prove (\ref{P10.6})
that for the quotient stack $[X_0/G_0],$ where $X_0$ is a
proper smooth variety and $G_0$ is a linear algebraic group
acting on $X_0,$ the Frobenius eigenvalues on its cohomology
groups are independent of $\ell.$

\begin{notation-convention}\label{NC}

\begin{anitem}
We fix a prime power $q=p^a$ and an algebraic closure 
$\bb F$ of the finite field $\bb F_q$ with $q$ elements. Let $F$
or $F_q$ be the $q$-geometric Frobenius, namely the $q$-th
root automorphism on $\mathbb F.$ Let $\ell$ be a prime
number, $\ell\ne p,$ and fix an isomorphism of fields
$\overline{\bb Q}_{\ell}\overset{\iota}{\to}\bb C.$
For simplicity, let $|\alpha|$ denote the complex absolute
value $|\iota\alpha|,$ for $\alpha\in\overline{\bb
Q}_{\ell}.$
\end{anitem}

\begin{anitem}
In this paper, by an Artin stack (or an algebraic stack) over
a base scheme $S,$ we mean an $S$-algebraic stack in the sense
of M. Artin (\cite{LMB}, 4.1) \textit{of finite type}. When we
want the more general setting of Artin stacks \textit{locally
of finite type}, we will mention that explicitly.
\end{anitem}

\begin{anitem}
Objects over $\bb F_q$ will be denoted with an index $_0.$
For instance, $\s X_0$ will denote an $\bb 
F_q$-Artin stack; if $\s F_0$ is a lisse-\'etale 
sheaf (or more generally, a Weil sheaf (\ref{Weil-cplx})) 
on $\s X_0,$ then $\mathscr F$ denotes its inverse 
image on $\s X:=\s X_0\otimes_{\bb F_q}\bb F.$
\end{anitem}

\begin{anitem}
For a field $k,$ let $\text{Gal}(k)$ denote its absolute
Galois group $\text{Gal}(k^\text{sep}/k).$ By a variety over
$k$ we mean a separated reduced $k$-scheme of finite type. 
Let $W(\bb F_q)$ be the \textit{Weil group} 
$F_q^{\bb Z}$ of $\bb F_q.$ 
\end{anitem}

\begin{anitem}
For a profinite group $H,$ by $\overline{\bb 
Q}_{\ell}$-representations of $H$ we always mean
finite-dimensional continuous representations (\cite{Del2},
1.1.6), and denote by $\text{Rep}_{\overline{\bb
Q}_{\ell}}(H)$ the category of such representations.
\end{anitem}

\begin{anitem}
For a scheme $X,$ we denote by $|X|$ the set of its closed
points. For a category $\s C$ we write $[\s 
C]$ for the collection of isomorphism classes of objects
in $\s C.$ For example, if $v\ge1$ is an integer, then
$[\s X_0(\bb F_{q^v})]$ denotes the set of
isomorphism classes of $\bb F_{q^v}$-points of the stack
$\s X_0.$ This is a finite set.

For $x\in\s X_0(\bb F_{q^v})$ we will write
$k(x)$ for the field $\bb F_{q^v}.$ For an $\bb
F_q$-scheme $X_0$ (always of finite type) and $x\in|X_0|,$ we
denote by $k(x)$ the residue field of $x.$ In both cases, let
$d(x)$ be the degree of the field extension $[k(x):\bb
F_q],$ and $N(x)=q^{d(x)}=\#k(x).$ Also in both cases let
$x:\text{Spec }\bb F_{q^v}\to\s X_0$ (or $X_0$) be
the natural map ($v=d(x)$), although in the second case the
map is defined only up to an automorphism in
$\text{Gal}(k(x)/\bb F_q).$ Given a $K_0\in D_c(\s 
X_0,\overline{\bb Q}_{\ell})$ (cf. $\S2$), the pullback
$x^*K_0\in D_c(\text{Spec }k(x),\overline{\bb
Q}_{\ell})=D_c(\text{Rep}_{\overline{\bb
Q}_{\ell}}(\text{Gal}(k(x))))$ gives a complex
$K_{\overline{x}}$ of representations of $\text{Gal}(k(x)),$
and we let $F_x$ be the geometric Frobenius generator
$F_{q^{d(x)}}$ of this group, called ``the local Frobenius".
\end{anitem}

\begin{anitem}
Let $V$ be a finite dimensional $\overline{\bb
Q}_{\ell}$-vector space and $\varphi$ an endomorphism of $V.$
For a function $f:\overline{\bb Q}_{\ell}\to\bb C,$
we denote by $\sum_{V,\varphi}f(\alpha)$ the sum of values of
$f$ in $\alpha,$ with $\alpha$ ranging over all the
eigenvalues of $\varphi$ on $V$ with multiplicities. For
instance, $\sum_{V,\varphi}\alpha=\text{Tr}(\varphi,V).$

A $0\times0$-matrix has trace 0 and determinant 1. For $K\in
D^b_c(\overline{\bb Q}_{\ell})$ and an endomorphism
$\varphi$ of $K,$ we define (following \cite{SGA4.5})
$$
\text{Tr}(\varphi,K):=\sum_{n\in\bb Z}(-1)^n\text{Tr}
(H^n(\varphi),H^n(K))
$$
and
$$
\det(1-\varphi t,K):=\prod_{n\in\bb Z}\det(1-H^n(\varphi)
t,H^n(K))^{(-1)^n}.
$$
For unbounded complexes $K$ we use similar notations, if the
series (resp. the infinite product) converges (resp.
converges term by term; cf. (\ref{D6.1})).
\end{anitem}

\begin{anitem}
For a map $f:X\to Y$ and a sheaf $\mathscr F$ on $Y,$ we
sometimes write $H^n(X,\s F)$ for $H^n(X,f^*\s F).$
We will write $H^n(\s X)$ for $H^n(\s X,\overline{\bb 
Q}_{\ell}),$ and $h^n(\s{X,F})$ for $\dim H^n(\s{X,F}),$ 
and ditto for $H^n_c(\s X)$ and $h^n_c(\s{X,F}).$
\end{anitem}

\begin{anitem}
For an $\bb F_q$-algebraic stack $\s X_0$ and a 
Weil complex $K_0$ on $\s X_0,$ by $R\Gamma(\s 
X_0,K_0)$ (resp. $R\Gamma_c(\s X_0,K_0)$) we mean 
$Ra_*K_0$ (resp. $Ra_!K_0$), where $a:\s X_0\to 
\text{Spec }\bb F_q$ is the structural map.

The derived functors $Rf_*,Rf_!,Lf^*$ and $Rf^!$ are 
usually abbreviated as $f_*,f_!,f^*,f^!.$ But we reserve 
$\otimes,\s Hom$ and $Hom$ for the ordinary 
sheaf tensor product, sheaf Hom and Hom group respectively, 
and use $\otimes^L,R\s Hom$ and $RHom$ for their 
derived functors.
\end{anitem}
\end{notation-convention}

\begin{flushleft}
\textbf{Acknowledgment.}
\end{flushleft}

I would like to thank my advisor Martin Olsson for introducing
this topic to me, and giving so many suggestions during the
writing. Weizhe Zheng sent me many helpful comments, and 
Yves Laszlo helped me to correct some inaccuracies in the first 
version of the paper. Many people, especially Brian Conrad and 
Matthew Emerton, have helped me on mathoverflow during the 
writing. The revision of the paper was done during the stay in Ecole 
polytechnique CMLS and Universit\'e Paris-Sud, while I was supported 
by ANR grant G-FIB.

\section{Derived category of $\ell$-adic sheaves and
mixedness}

We briefly review the definition in \cite{LO1, LO2} for
derived category of $\ell$-adic sheaves on stacks. Then we
show that $\iota$-mixedness is stable under the six
operations. As a consequence of Lafforgue's result (\ref{Laff}), 
this is automatic, but we still want to give a much more elementary 
argument. The proof works for \textit{mixed} complexes as well 
(\ref{mixed-variant}). One can also generalize the 
structure theorem of $\iota$-mixed sheaves in \cite{Del2} 
to algebraic stacks (\ref{R2.7}). 

\begin{blank}\label{adic-setting}
Let $\Lambda$ be a complete discrete valuation ring with
maximal ideal $\fk m$ and residual characteristic
$\ell.$ Let $\Lambda_n=\Lambda/\fk m^{n+1},$ and
let $\Lambda_{\bullet}$ be the pro-ring $(\Lambda_n)_n.$ We
take the base scheme $S$ to be a scheme that satisfies the
following condition denoted (LO): it is noetherian affine
excellent finite-dimensional, in which $\ell$ is invertible,
and all $S$-schemes of finite type have finite
$\ell$-cohomological dimension. We denote by $\cal{X,Y}\cdots$
Artin stacks locally of finite type over $S.$

Consider the ringed topos $\s A=\s A(\mathcal
X):=\text{Mod}(\mathcal
X_{\text{lis-\'et}}^{\bb N},\Lambda_{\bullet})$ of
projective systems $(M_n)_n$ of $\text{Ab}(\mathcal
X_{\text{lis-\'et}})$ such that $M_n$ is a $\Lambda_n$-module
for each $n,$ and the transition maps are $\Lambda$-linear.
An object $M\in\s A$ is said to be \textit{AR-null}, if
there exists an integer $r>0$ such that for every integer
$n,$ the composed map $M_{n+r}\to M_n$ is the zero map. A
complex $K$ in $\s A$ is called \textit{AR-null}, if
all cohomology systems $\s H^i(K)$ are AR-null; it is
called \textit{almost AR-null}, if for every $U$ in
$\text{Lis-\'et}(\mathcal X)$ (assumed to be of finite type 
over $S$), the restriction of $\s 
H^i(K)$ to $\text{\'Et}(U)$ is AR-null. Let $\s{D(A)}$
be the ordinary derived category of $\s A.$ See 
(\cite{LMB}, 18.1.4) for the definition of constructible 
sheaves on $\mathcal X_{\text{lis-\'et}}.$ 
\end{blank}

\begin{definition}\label{D2.1}
An object $M=(M_n)_n\in\s A$ is \emph{adic} if all the
$M_n$'s are constructible, and for every $n,$ the natural map
$$
\Lambda_n\otimes_{\Lambda_{n+1}}M_{n+1}\to M_n
$$
is an isomorphism. It is called \emph{almost adic} if all the
$M_n$'s are constructible, and for every object $U$ in
$\emph{Lis-\'et}(\mathcal X),$ the restriction $M|_U$ is
AR-adic, i.e. there exists an adic $N_U\in\emph{Mod}(U_{\emph
{\'et}}^{\bb N},\Lambda_{\bullet})$ and a morphism $N_U\to
M|_U$ with AR-null kernel and cokernel.

A complex $K$ in $\s A$ is a $\lambda$-\emph{complex}
if $\s H^i(K)\in\s A$ are almost adic, for all
$i.$ Let $\s D_c(\s A)$ be the full triangulated
subcategory of $\s{D(A)}$ consisting of
$\lambda$-complexes, and let $D_c(\mathcal X,\Lambda)$ be
the quotient of $\s D_c(\s A)$ by the thick full
subcategory of those which are almost AR-null. This is called the
\emph{derived category of $\Lambda$-adic sheaves on}
$\mathcal X.$
\end{definition}

\begin{subremark}\label{R2.2}
(i) $D_c(\mathcal X,\Lambda)$ is a triangulated category with
a natural $t$-structure, and its heart is the quotient of the
category of almost adic systems in $\s A$ by the thick
full subcategory of almost AR-null systems. One can use this
$t$-structure to define the subcategories
$D_c^{\dagger}(\mathcal X,\Lambda)$ for $\dagger=\pm, b.$

If $\mathcal X/S$ is of finite type (in particular,
quasi-compact), it is clear that $K\in\s D_{\text{cart}}(\s A)$ is
AR-null if it is almost AR-null. Also if $M\in\s A$
is almost adic, the adic system $N_U$ and the map $N_U\to
M|_U$ in the definition above are unique up to unique
isomorphism, for each $U,$ so by (\cite{LMB}, 12.2.1) they
give an adic system $N$ of Cartesian sheaves on $\mathcal
X_{\text{lis-\'et}},$ and an AR-isomorphism $N\to M.$ This
shows that an almost adic system is AR-adic, and it is
clear (\cite{SGA5}, p.234) that the natural functor
$$
\Lambda\text{-Sh}(\mathcal X)\to\text{heart }D_c(\mathcal
X,\Lambda)
$$
is an equivalence of categories, where $\Lambda$-Sh$(\mathcal
X)$ denotes the category of $\Lambda$-adic (in particular,
constructible) systems.

(ii) $D_c(\mathcal X,\Lambda)$ is different from the ordinary
derived category of $\text{Mod}(\mathcal
X_{\text{lis-\'et}},\Lambda)$ with constructible cohomology;
the latter can be denoted by $\mathscr
D_c(\mathcal X,\Lambda).$ Here $\text{Mod}(\mathcal
X_{\text{lis-\'et}},\Lambda)$ denotes the abelian category of
$\Lambda_{\mathcal X}$-modules, not adic sheaves
$\Lambda\text{-Sh}(\mathcal X).$

(iii) When $S=\text{Spec }k$ for $k$ a finite field or an
algebraically closed field, and $\mathcal X=X$ is a separated
$S$-scheme, (\cite{LO2}, 3.1.6) gives a natural equivalence of
triangulated categories between $D^b_c(X,\Lambda)$ and
Deligne's definition $\mathscr D_c^b(X,\Lambda)$ in
(\cite{Del2}, 1.1.2).
\end{subremark}

\begin{blank}\label{normalization}
Let $\pi:\mathcal X^{\mathbb N}_{\text{lis-\'et}}\to\mathcal
X_{\text{lis-\'et}}$ be the morphism of topoi, where
$\pi^{-1}$ takes a sheaf $F$ to the constant projective
system $(F)_n,$ and $\pi_*$ takes a projective system to
the inverse limit. This morphism induces a morphism of
ringed topoi $(\pi^*,\pi_*):(\mathcal X^{\bb 
N}_{\text{lis-\'et}},\Lambda_{\bullet})\to(\mathcal
X_{\text{lis-\'et}},\Lambda).$ The functor $R\pi_*:\mathscr
D_c(\mathscr A)\to\mathscr D(\mathcal X,\Lambda)$ vanishes
on almost AR-null objects (\cite{LO2}, 2.2.2), hence factors
through $D_c(\mathcal X,\Lambda).$ In (\cite{LO2}, 3.0.8),
the normalization functor is defined to be
$$
K\mapsto\widehat{K}:=L\pi^*R\pi_*K:\ D_c(\mathcal X,\Lambda)
\to\s{D(A)}.
$$
This functor plays an important role in defining the six
operations \cite{LO2}. For instance:

$\bullet$ For $F\in D_c^-(\mathcal X,\Lambda)$ and $G\in
D_c^+(\mathcal X,\Lambda),\ R\s Hom(F,G)$ is defined
to be the image of $R\s Hom_{\Lambda_{\bullet}}(\widehat
{F},\widehat{G})$ in $D_c(\mathcal X,\Lambda).$

$\bullet$ For $F,G\in D_c^-(\mathcal X,\Lambda),$ the
derived tensor product $F\otimes^LG$ is defined to be the
image of
$\widehat{F}\otimes^L_{\Lambda_{\bullet}}\widehat{G}.$

$\bullet$ For a morphism $f:\mathcal X\to\mathcal Y$ and
$F\in D_c^+(\mathcal X,\Lambda),$ the derived direct image
$f_*F$ is defined to be the image of $f^{\bb 
N}_*\widehat{F}.$

Let $E_{\lambda}$ be a finite extension of $\bb
Q_{\ell}$ with ring of integers $\mathscr O_{\lambda}.$
Following \cite{LO2} we define $D_c(\mathcal X,E_{\lambda})$
to be the quotient of $D_c(\mathcal X,\s O_{\lambda})$
by the full subcategory consisting of complexes $K$ such
that, for every integer $i,$ there exists an integer
$n_i\ge1$ such that $\mathscr H^i(K)$ is annihilated by
$\lambda^{n_i}.$ Then we define
$$
D_c(\mathcal X,\overline{\bb Q}_{\ell})=\text{2-colim}
_{E_{\lambda}}D_c(\mathcal X,E_{\lambda}),
$$
where $E_{\lambda}$ ranges over all finite subextensions of
$\overline{\bb Q}_{\ell}/\bb Q_{\ell},$ and the
transition functors are
$$
E_{\lambda'}\otimes_{E_{\lambda}}-:D_c(\mathcal
X,E_{\lambda})\to D_c(\mathcal X,E_{\lambda'})
$$
for $E_{\lambda}\subset E_{\lambda'}.$
\end{blank}

\begin{blank}\label{Weil-cplx}
From now on in this section, $S=\text{Spec }\bb F_q.$ We
recall some notions of weights and mixedness from \cite{Del2},
generalized to $\mathbb F_q$-algebraic stacks. 

\begin{anitem}
\textbf{Frobenius endomorphism.} For an $\bb 
F_q$-scheme $X_0,$ let 
$F_{X_0}:X_0\to X_0$ be the morphism that is identity on 
the underlying topological space and $q$-th power on the 
structure sheaf $\s O_{X_0};$ this is an $\bb 
F_q$-morphism. Let $F_X:X\to X$ be the induced 
$\bb F$-morphism $F_{X_0}\times\text{id}_{\bb F}$ 
on $X=X_0\otimes\bb F.$

By functoriality of the maps $F_{X_0},$ we can extend it 
to stacks. For an $\bb F_q$-algebraic stack $\s 
X_0,$ define $F_{\s X_0}:\s X_0\to\s 
X_0$ to be such that for every $\bb F_q$-scheme 
$X_0,$ the map 
$$
\xymatrix@C=.7cm{
F_{\s X_0}(X_0):\s X_0(X_0) \ar[r] & \s X_0(X_0)}
$$
sends $x$ to $x\circ F_{X_0}.$ We also define 
$F_{\s X}:\s X\to\s X$ to be $F_{\s X_0}\times
\text{id}_{\bb F}.$ This morphism is a universal 
homeomorphism, hence $F_{\s X}^*$ and $F_{\s 
X*}$ are quasi-inverse to each other, and $F_{\s 
X}^*\simeq F_{\s X}^!,F_{\s X*}\simeq F_{\s X!}.$
\end{anitem}

\begin{anitem}
\textbf{Weil complexes.} A \textit{Weil complex $K_0$ on 
$\s X_0$} is a pair $(K,\varphi),$ where $K\in 
D_c(\s X,\overline{\bb Q}_{\ell})$ and 
$\varphi:F^*_{\s X}K\to K$ is an isomorphism. A morphism 
of Weil complexes on $\s X_0$ is a morphism of complexes on 
$\s X$ commuting with $\varphi.$ We 
also call $K_0$ a \textit{Weil sheaf} if $K$ is a 
sheaf. Let $W(\s X_0,\overline{\bb 
Q}_{\ell})$ be the category of Weil complexes on 
$\s X_0;$ it is a triangulated category with the 
standard $t$-structure, and its core is the category of 
Weil sheaves. There is a natural fully faithful triangulated 
functor 
$$
D_c(\s X_0,\overline{\bb Q}_{\ell})\to W(\s X_0,\overline{\bb 
Q}_{\ell}).
$$

The usual six operations are well-defined on Weil complexes. 

$\bullet$ Verdier duality. The Weil complex structure on 
$D_{\s X}K$ is given by the inverse of the isomorphism 
$$
\xymatrix@C=.6cm{
D_{\s X}K \ar[r]^-{D\varphi} & D_{\s X}
F_{\s X}^*K \ar[r]^-{\sim} & F_{\s X}^*D_{\s X}K.}
$$

$\bullet$ Tensor product. Let $K_0$ and $L_0$ be two Weil 
complexes such that $K\otimes^LL$ (which is $K\otimes L$ 
since they are of $\overline{\bb Q}_{\ell}$-coefficients) 
is constructible. This is the case when they are both 
bounded above. The Weil complex structure on $K\otimes L$ is 
given by 
$$
\xymatrix@C=.6cm{
F_{\s X}^*(K\otimes L) \ar[r]^-{\sim} & F_{\s 
X}^*K\otimes F_{\s X}^*L \ar[rr]^-{\varphi_K\otimes 
\varphi_L} && K\otimes L.}
$$

$\bullet$ Pullback. This is clear:
$$
\xymatrix@C=.6cm{
F_{\s X}^*f^*K \ar[r]^-{\sim} & f^*F_{\s Y}^*K 
\ar[r]^-{f^*\varphi} & f^*K.}
$$
Here $f:\s X_0\to\s Y_0$ is an $\bb 
F_q$-morphism and $(K,\varphi)$ is a Weil complex on 
$\s Y_0.$

$\bullet$ Pushforward. Let $f:\s X_0\to\s Y_0$ 
and $K_0\in W^+(\s X_0,\overline{\bb Q}_{\ell}).$ 
The Weil complex structure on $f_*K$ is given by 
$$
\xymatrix@C=.6cm{
F_{\s Y}^*f_*K \ar[r] & f_*F_{\s X}^*K 
\ar[r]^-{f_*\varphi} & f_*K,}
$$
where the first arrow is an isomorphism, because it is 
adjoint to 
$$
f_*K\to F_{\s Y*}f_*F_{\s X}^*K\simeq f_* 
F_{\s X*}F_{\s X}^*K
$$
obtained by applying $f_*$ to the adjunction morphism 
$K\to F_{\s X*}F_{\s X}^*K,$ which is an isomorphism. 

$\bullet$ The remaining cases $f^!,f_!$ and $R\s 
Hom$ follow from the previous cases.

In this article, when discussing stacks over $\bb F_q,$ 
by a ``sheaf" or ``complex of sheaves", we usually mean a 
``Weil sheaf" or ``Weil complex", whereas a ``lisse-\'etale 
sheaf or complex" will be an ordinary constructible 
$\overline{\bb Q}_{\ell}$-sheaf or complex on the 
lisse-\'etale site of $\s X_0.$

For $x\in\s X_0(\bb F_{q^v}),$ it is a fixed 
point of $F_{\s X}^v,$ hence there is a 
``local Frobenius automorphism" $F_x:K_{\overline{x}}\to 
K_{\overline{x}}$ for every Weil complex $K_0,$ defined 
to be 
$$
K_{\overline{x}}\simeq K_{F_{\mathscr X}(\overline{x})} 
=(F_{\s X}^*K)_{\overline{x}}\overset{\varphi}{\to}
K_{\overline{x}}.
$$
\end{anitem}

\begin{anitem}
\textbf{$\iota$-Weights and $\iota$-mixedness.} 
Recall that $\iota$ is a fixed isomorphism 
$\overline{\bb Q}_{\ell}\to\bb C.$ For
$\alpha\in\overline{\bb Q}_{\ell}^*,$ let
$w_q(\alpha):=2\log_q|\iota\alpha|,$ called the $\iota$-weight
of $\alpha$ relative to $q.$ For a real number $\beta,$ a
sheaf $\s F_0$ on $\s X_0$ is said to be
\textit{punctually $\iota$-pure of weight} $\beta,$ if for every
integer $v\ge1$ and every $x\in\s X_0(\bb F_{q^v}),$
and every eigenvalue $\alpha$ of $F_x$ acting on $\s 
F_{\overline{x}},$ we have $w_{N(x)}(\alpha)=\beta.$ We say
$\s F_0$ is \textit{$\iota$-mixed} if it has a finite filtration
with successive quotients punctually $\iota$-pure, and the
weights of these quotients are called the punctual
$\iota$-weights of $\s F_0.$ A complex $K_0\in 
W(\s X_0,\overline{\bb Q}_{\ell})$ is said to 
be $\iota$-mixed if all the cohomology sheaves $\s 
H^iK_0$ are $\iota$-mixed. Let $W_m(\s X_0,\overline{\bb 
Q}_{\ell})$ (resp. $D_m(\s X_0,\overline{\bb Q}_{\ell})$) 
be the full subcategory of $\iota$-mixed complexes in 
$W(\s X_0,\overline{\bb Q}_{\ell})$ (resp. $D_c(\s 
X_0,\overline{\bb Q}_{\ell})$).

One can also define \textit
{punctually pure sheaves, mixed sheaves} and \textit 
{mixed complexes} for algebraic stacks. 
\end{anitem}

\begin{anitem}
\textbf{Twists.} For $b\in\overline{\mathbb Q}_{\ell}^*,$ 
let $\overline{\mathbb Q}_{\ell}^{(b)}$ be the Weil sheaf 
on $\text{Spec }\mathbb F_q$ of rank one, where $F$ acts by 
multiplication by $b.$ This is an \'etale sheaf if and only 
if $b$ is an $\ell$-adic unit (\cite{Del2}, 1.2.7). For an 
algebraic stack $\mathscr X_0/\mathbb F_q,$ we also denote 
by $\overline{\mathbb Q}_{\ell}^{(b)}$ the inverse image on 
$\mathscr X_0$ of the above Weil sheaf via the structural 
map. If $\mathscr F_0$ is a sheaf on $\mathscr X_0,$ we 
denote by $\mathscr F_0^{(b)}$ the tensor product $\mathscr 
F_0\otimes\overline{\mathbb Q}_{\ell}^{(b)},$ and say that 
$\mathscr F_0^{(b)}$ is deduced from $\mathscr F_0$ by 
a generalized Tate twist. Note that the 
operation $\mathscr F_0\mapsto\mathscr F_0^{(b)}$ adds the 
weights by $w_q(b).$ For an integer $d,$ the usual Tate 
twist $\overline{\mathbb Q}_{\ell}(d)$ is 
$\overline{\mathbb Q}_{\ell}^{(q^{-d})}.$ We denote by
$\langle d\rangle$ the operation $(d)[2d]$ on complexes of
sheaves, where $[2d]$ means shifting $2d$ to the left. Note
that $\iota$-mixedness is stable under the operation $\langle
d\rangle.$
\end{anitem}
\end{blank}

\begin{lemma}\label{L2.3}
Let $\s X_0$ be an $\bb F_q$-algebraic stack.

(i) If $\s F_0$ is an $\iota$-mixed sheaf on $\s 
X_0,$ then so is every sub-quotient of $\s F_0.$

(ii) If $0\to\s F_0'\to\s F_0\to\s F_0''\to0$ is an exact 
sequence of sheaves on $\s X_0,$ and $\s F_0'$ and $\s F_0''$ 
are $\iota$-mixed, then so is $\s F_0.$

(iii) The full subcategory $W_m(\s X_0,\overline{\bb 
Q}_{\ell})$ (resp. $D_m(\s X_0,\overline{\bb Q}_{\ell})$) 
of $W(\s X_0,\overline{\bb Q}_{\ell})$ (resp. $D_c(\s
X_0,\overline{\bb Q}_{\ell})$) is a triangulated 
subcategory with induced standard $t$-structure. 

(iv) If $f$ is a morphism of $\bb F_q$-algebraic stacks,
then $f^*$ on complexes of sheaves preserves
$\iota$-mixedness.

(v) If $j:\s U_0\hookrightarrow\s X_0$ is an open
immersion and $i:\s Z_0\hookrightarrow\s X_0$ 
is its complement, then $K_0\in W(\s 
X_0,\overline{\bb Q}_{\ell})$ is $\iota$-mixed if and
only if $j^*K_0$ and $i^*K_0$ are $\iota$-mixed.
\end{lemma}

\begin{proof}
(i) If $\mathscr F_0$ is punctually $\iota$-pure of weight
$\beta,$ then so is every sub-quotient of it. Now suppose
$\mathscr F_0$ is $\iota$-mixed and $\mathscr F_0'$ is a
subsheaf of $\mathscr F_0.$ Let $W$ be a finite filtration
$$
0\subset\cdots\subset\mathscr F_0^{i-1}\subset\mathscr F_0^i
\subset\cdots\subset\mathscr F_0
$$
of $\mathscr F_0$ such that $\text{Gr}^W_i(\mathscr F_0):=
\mathscr F_0^i/\mathscr F_0^{i-1}$ is punctually $\iota$-pure
for every $i.$ Let $W'$ be the induced filtration
$W\cap\mathscr F_0'$ of $\mathscr F_0'.$ Then
$\text{Gr}^{W'}_i(\mathscr F_0')$ is the image of
$$
\mathscr F_0^i\cap\mathscr F_0'\subset\mathscr F_0^i
\twoheadrightarrow\text{Gr}^W_i(\mathscr F_0),
$$
so it is punctually $\iota$-pure. Let $\mathscr F_0''=\mathscr
F_0/\mathscr F_0'$ be a quotient of $\mathscr F_0,$ and let
$W''$ be the induced filtration of $\mathscr F_0'',$ namely
$(\mathscr F_0'')^i:=\mathscr F_0^i/(\mathscr
F_0^i\cap\mathscr F_0').$ Then $\text{Gr}^{W''}_i(\mathscr
F_0'')=\mathscr F_0^i/(\mathscr F_0^{i-1}+\mathscr
F_0^i\cap\mathscr F_0'),$ which is a quotient of $\mathscr
F_0^i/\mathscr F_0^{i-1}=\text{Gr}^W_i(\mathscr F_0),$ so it
is punctually $\iota$-pure. Hence every sub-quotient
of $\mathscr F_0$ is $\iota$-mixed.

(ii) Let $W'$ and $W''$ be finite filtrations
of $\mathscr F_0'$ and $\mathscr F_0''$ respectively, such
that $\text{Gr}^{W'}_i(\mathscr F_0')$ and
$\text{Gr}^{W''}_i(\mathscr F_0'')$ are punctually
$\iota$-pure for every $i.$ Then $W'$ can be regarded as a
finite filtration of $\mathscr F_0$ such that every member
of the filtration is contained in $\mathscr F_0',$ and
$W''$ can be regarded as a finite filtration of $\mathscr
F_0$ such that every member contains $\mathscr F_0'.$
Putting these two filtrations together,
we get the desired filtration for $\mathscr F_0.$

(iii) Being a triangulated subcategory means (\cite{SGA4.5},
p.271) that, if $K_0'\to K_0\to K_0''\to K_0'[1]$ is an exact
triangle in $W(\mathscr X_0,\overline{\mathbb 
Q}_{\ell}),$ and two of the three complexes are 
$\iota$-mixed, then so is the third. By the rotation axiom 
of a triangulated category, we can assume $K_0'$ and 
$K_0''$ are $\iota$-mixed. We have the exact sequence
$$
\xymatrix@C=.5cm{
\cdots \ar[r] & \mathscr H^nK_0' \ar[r] & \mathscr H^nK_0
\ar[r] & \mathscr H^nK_0'' \ar[r] & \cdots,}
$$
and by (i) and (ii) we see that $\mathscr H^nK_0$ is
$\iota$-mixed.

$W_m(\mathscr X_0,\overline{\mathbb Q}_{\ell})$ has the
induced $t$-structure because if $K_0$ is $\iota$-mixed, then
its truncations $\tau_{\le n}K_0$ and $\tau_{\ge n}K_0$ are
$\iota$-mixed.

(iv) On sheaves, $f^*$ preserves stalks, so it is exact and
preserves punctual $\iota$-purity on sheaves. Let $f:\mathscr
X_0\to\mathscr Y_0.$ Given an $\iota$-mixed sheaf $\mathscr
F_0$ on $\mathscr Y_0,$ let $W$ be a finite filtration of
$\mathscr F_0$ such that each $\text{Gr}^W_i(\mathscr F_0)$
is punctually $\iota$-pure. Then $f^*W$ gives a finite
filtration of $f^*\mathscr F_0$ and each
$\text{Gr}^{f^*W}_i(f^*\mathscr
F_0)=f^*\text{Gr}^W_i(\mathscr F_0)$ is punctually
$\iota$-pure. So $f^*\mathscr F_0$ is $\iota$-mixed.
For an $\iota$-mixed complex $K_0$ on $\mathscr Y_0,$ note
that $\mathscr H^n(f^*K_0)=f^*\mathscr H^n(K_0),$ hence
$f^*K_0$ is $\iota$-mixed.

(v) One direction follows from (iv). For the other direction, 
note that $j_!$ and $i_*$ are exact and preserve punctual 
$\iota$-purity on sheaves. If $\s F_0$ is an 
$\iota$-mixed sheaf on $\s U_0,$ with a finite
filtration $W$ such that each $\text{Gr}^W_i(\s F_0)$ 
is punctually $\iota$-pure, then for the induced 
filtration $j_!W$ of $j_!\s F_0,$ we see that
$\text{Gr}^{j_!W}_i(j_!\s F_0)=j_!\text{Gr}^W_i(\s F_0)$ 
is punctually $\iota$-pure, so 
$j_!\s F_0$ is $\iota$-mixed. For an $\iota$-mixed
complex $K_0$ on $\s U_0,$ use $\s H^n(j_!K_0)=j_!\s H^n(K_0).$ 
Similarly $i_*$ also preserves $\iota$-mixedness on complexes.
Finally the result follows from (iii) and the exact triangle
$$
\xymatrix@C=.5cm{
j_!j^*K_0 \ar[r] & K_0 \ar[r] & i_*i^*K_0 \ar[r] &.}
$$
\end{proof}

To show that $\iota$-mixedness is stable under the six
operations, we need to show that $\iota$-mixedness of
complexes on stacks can be checked locally on their
presentations. To descend a filtration on a presentation to
the stack, we generalize the structure theorem of
$\iota$-mixed sheaves to algebraic spaces. Recall the
following theorem of Deligne (\cite{Del2}, 3.4.1).

\begin{theorem}\label{T2.4}
Let $\mathscr F_0$ be an $\iota$-mixed sheaf on a scheme
$X_0$ over $\mathbb F_q.$

(i) $\mathscr F_0$ has a unique decomposition $\mathscr
F_0=\bigoplus_{b\in\mathbb R/\mathbb Z}\mathscr F_0(b),$
called the \emph{decomposition according to the weights mod
$\mathbb Z,$} such that the punctual $\iota$-weights of
$\mathscr F_0(b)$ are all in the coset $b.$ This
decomposition, in which almost all the $\mathscr F_0(b)$ are
zero, is functorial in $\mathscr F_0.$ Note that each
$\mathscr F_0(b)$ is deduced by twist from an
$\iota$-mixed sheaf with integer punctual weights.

(ii) If the punctual weights of $\mathscr F_0$ are integers
and $\mathscr F_0$ is lisse, $\mathscr F_0$ has a unique
finite increasing filtration $W$ by lisse subsheaves, called
the \emph{filtration by punctual weights,} such that
$\emph{Gr}_i^W(\mathscr F_0)$ is punctually $\iota$-pure of
weight $i.$ This filtration is functorial in $\mathscr F_0.$
More precisely, any morphism between $\iota$-mixed lisse
sheaves of integer punctual weights is strictly compatible
with their filtrations by punctual weights.

(iii) If $\mathscr F_0$ is lisse and punctually $\iota$-pure,
and $X_0$ is normal, then the sheaf $\mathscr F$ on $X$ is
semi-simple.
\end{theorem}

\begin{subremark}\label{R2.5}
(i) If $\mathscr C$ is an abelian category and $\mathscr D$ is
an abelian full subcategory of $\mathscr C,$ and $C$ is an
object in $\mathscr D,$ then every direct summand of $C$ in
$\mathscr C$ lies in $\mathscr D$ (or isomorphic to some
object in $\mathscr D$). This is because the kernel of the
composition
$$
\xymatrix@C=.7cm{
A\oplus B \ar @{->>}[r]^-{\emph{pr}_A} & A
\ar @{^{(}->}[r]^-{i_A} & A\oplus B}
$$
is $B.$ So direct summands of a lisse sheaf are lisse. If
$\mathscr F_0$ in (\ref{T2.4}i) is lisse, 
then each $\mathscr F_0(b)$ is lisse.

(ii) If the $\overline{\mathbb Q}_{\ell}$-sheaf $\mathscr F_0$
is defined over some finite subextension $E_{\lambda}$ of
$\overline{\mathbb Q}_{\ell}/\mathbb Q_{\ell},$ then its
decomposition in (\ref{T2.4}i) and filtration in
(\ref{T2.4}ii) are defined over $E_{\lambda}.$ This is because
the $E_{\lambda}$-action commutes with the Galois action.

(iii) In \cite{Del2} Deligne made the assumption that all
schemes are separated, at least in order to use Nagata
compactification to define $f_!.$ After the work of Laszlo and
Olsson \cite{LO1, LO2}, one can remove this assumption, and
many results in \cite{Del2}, for instance this one and
(3.3.1), remain valid. For (\cite{Del2}, 3.4.1) one can take a cover of a
not necessarily separated scheme $X_0$ by open affines (which
are separated), and use the functoriality to glue the
decomposition or filtration on intersections.
\end{subremark}

\begin{lemma}\label{L2.6}
Let $X_0$ be an $\mathbb F_q$-algebraic space, and $\mathscr
F_0$ an $\iota$-mixed sheaf on $X_0.$

(i) $\mathscr F_0$ has a unique decomposition $\mathscr
F_0=\bigoplus_{b\in\mathbb R/\mathbb Z}\mathscr F_0(b),$ the
\emph{decomposition according to the weights mod $\mathbb Z,$}
with the same property as in (\ref{T2.4}i). This
decomposition is functorial in $\mathscr F_0.$

(ii) If the punctual $\iota$-weights of $\mathscr F_0$ are
integers and $\mathscr F_0$ is lisse, $\mathscr F_0$ has a
unique finite increasing filtration $W$ by lisse subsheaves,
called the \emph{filtration by punctual weights,} with the
same property as in (\ref{T2.4}ii). This filtration is
functorial in $\mathscr F_0.$
\end{lemma}

\begin{proof}
Let $P:X_0'\to X_0$ be an \'etale presentation, and let
$\mathscr F_0'=P^*\mathscr F_0,$ which is also $\iota$-mixed
(\ref{L2.3}iv). Let $X_0''$ be the fiber product
$$
X_0''=\xymatrix@C=.7cm{
X_0'\times_{X_0}X_0' \ar[r]^-{p_1} \ar[d]_-{p_2} &
X_0' \ar[d]^-P \\
X_0' \ar[r]_-P & X_0.}
$$
Then $X_0''$ is an $\mathbb F_q$-scheme of finite type.

(i) Applying (\ref{T2.4}i) to $\mathscr F_0'$ we get a
decomposition $\mathscr F_0'=\bigoplus_{b\in\mathbb R/\mathbb
Z}\mathscr F_0'(b).$ For $j=1,2,$ applying $p_j^*$ we get a
decomposition
$$
p_j^*\mathscr F_0'=\bigoplus_{b\in\mathbb R/\mathbb
Z}p_j^*\mathscr F_0'(b).
$$
Since $p_j^*$ preserves weights, by the uniqueness in
(\ref{T2.4}i), this decomposition is the decomposition of
$p_j^*\mathscr F_0'$ according to the weights mod $\mathbb Z.$
By the functoriality in (\ref{T2.4}i), the canonical
isomorphism $\mu:p_1^*\mathscr F_0'\to p_2^*\mathscr F_0'$
takes the form $\bigoplus_{b\in\mathbb R/\mathbb Z}\mu_b,$
where $\mu_b:p_1^*\mathscr F_0'(b)\to p_2^*\mathscr F_0'(b)$
is an isomorphism satisfying cocycle condition as $\mu$ does.
Therefore the decomposition $\mathscr F_0'=\bigoplus_{b\in
\mathbb R/\mathbb Z}\mathscr F_0'(b)$ descends to a
decomposition $\mathscr F_0=\bigoplus_{b\in\mathbb R/\mathbb
Z}\mathscr F_0(b).$ We still need to show each direct summand
$\mathscr F_0(b)$ is $\iota$-mixed.

Fix a coset $b$ and consider the summand $\mathscr F_0(b).$
Twisting it appropriately, we can assume that its
inverse image $\mathscr F_0'(b)$ is $\iota$-mixed with
integer punctual $\iota$-weights. By (\ref{L2.3}v)and
noetherian induction, we can
shrink $X_0$ to a nonempty open subspace and assume $\mathscr
F_0(b)$ is lisse. Then $\mathscr F_0'(b)$ is also lisse, and
applying (\ref{T2.4}ii) we get a finite increasing filtration
$W'$ of $\mathscr F_0'(b)$ by lisse subsheaves $\mathscr
F_0'(b)^i,$ such that each $\text{Gr}^{W'}_i(\mathscr
F_0'(b))$ is punctually $\iota$-pure of weight $i.$ Pulling
back this filtration via $p_j,$ we get a finite increasing
filtration $p_j^*W'$ of $p_j^*\mathscr F_0'(b),$ and since
$\text{Gr}_i^{p_j^*W'}(p_j^*\mathscr
F_0'(b))=p_j^*\text{Gr}_i^{W'}(\mathscr F_0'(b))$ is
punctually $\iota$-pure of weight $i,$ it is the filtration by
punctual weights given by (\ref{T2.4}ii), hence functorial.
So the canonical isomorphism $\mu_b:p_1^*\mathscr F_0'(b)\to
p_2^*\mathscr F_0'(b)$ maps $p_1^*\mathscr F_0'(b)^i$
isomorphically onto $p_2^*\mathscr F_0'(b)^i,$ satisfying
cocycle condition. Therefore the filtration $W'$ of $\mathscr
F_0'(b)$ descends to a filtration $W$ of $\mathscr F_0(b),$
and $P^*\text{Gr}_i^W(\mathscr F_0(b))=\text{Gr}_i
^{W'}(\mathscr F_0'(b))$ is punctually $\iota$-pure of weight
$i.$ Note that $P$ is surjective, so every point $x\in
X_0(\mathbb F_{q^v})$ can be lifted to a point $x'\in
X_0'(\mathbb F_{q^{nv}})$ after some base extension $\mathbb
F_{q^{nv}}$ of $\mathbb F_{q^v}.$ This shows
$\text{Gr}_i^W(\mathscr F_0(b))$ is punctually $\iota$-pure
of weight $i,$ therefore $\mathscr F_0(b)$ is $\iota$-mixed.
This proves the existence of the decomposition in (i).

For uniqueness, let $\mathscr F_0=\bigoplus\widetilde{
\mathscr F}_0(b)$ be another decomposition with the desired
property. Then their restrictions to $X_0'$ are both equal to
the decomposition of $\mathscr F_0',$ which is unique
(\ref{T2.4}i), so they are both obtained by descending this
decomposition, and so they are isomorphic, i.e. for every
coset $b$ there exists an isomorphism making the diagram
commute:
$$
\xymatrix@C=.8cm{
\mathscr F_0(b) \ar[rr]^-{\sim} \ar@{^{(}->}[rd] &&
\widetilde{\mathscr F}_0(b) \ar@{^{(}->}[ld] \\
& \mathscr F_0. &}
$$

For functoriality, let $\mathscr G_0=\bigoplus\mathscr
G_0(b)$ be another $\iota$-mixed sheaf with decomposition on
$X_0,$ and let $\varphi:\mathscr F_0\to\mathscr G_0$ be a
morphism of sheaves. Pulling $\varphi$ back via $P$ we get
a morphism $\varphi':\mathscr F_0'\to\mathscr G_0'$ on $X_0',$
and the diagram
$$
\xymatrix@C=.9cm{
p_1^*\mathscr F_0' \ar[r]^-{\mu_{\mathscr F_0}}
\ar[d]_-{p_1^*\varphi'} & p_2^*\mathscr F_0'
\ar[d]^-{p_2^*\varphi'} \\
p_1^*\mathscr G_0' \ar[r]_-{\mu_{\mathscr G_0}} &
p_2^*\mathscr G_0'}
$$
commutes. By (\ref{T2.4}i) $\varphi'=\bigoplus\varphi'(b)$
for morphisms $\varphi'(b):\mathscr F_0'(b)\to\mathscr
G_0'(b),$ and the diagram
$$
\xymatrix@C=.9cm{
p_1^*\mathscr F_0'(b) \ar[r]^-{\text{can}}
\ar[d]_-{p_1^*\varphi'} & p_2^*\mathscr F_0'(b)
\ar[d]^-{p_2^*\varphi'} \\
p_1^*\mathscr G_0'(b) \ar[r]_-{\text{can}} &
p_2^*\mathscr G_0'(b)}
$$
commutes for each $b.$ Then the morphisms $\varphi'(b)$
descend to morphisms $\varphi(b):\mathscr F_0(b)\to\mathscr
G_0(b)$ such that $\varphi=\bigoplus\varphi(b).$

(ii) The proof is similar to part (i). Applying
(\ref{T2.4}ii) to $\mathscr F_0'$ on $X_0'$ we get a finite
increasing filtration $W'$ of $\mathscr F_0'$ by lisse
subsheaves $\mathscr F_0'^i$ with desired property. Pulling
back this filtration via $p_j:X_0''\to X_0'$ we get the
filtration by punctual weights of $p_j^*\mathscr F_0'.$ By
functoriality in (\ref{T2.4}ii), the canonical isomorphism
$\mu:p_1^*\mathscr F_0'\to p_2^*\mathscr F_0'$ maps
$p_1^*\mathscr F_0'^i$ isomorphically onto $p_2^*\mathscr
F_0'^i$ satisfying cocycle condition, therefore the
filtration $W'$ descends to a finite increasing filtration
$W$ of $\mathscr F_0$ by certain subsheaves $\mathscr F_0^i.$
By (\cite{Ols3}, 9.1) they are lisse subsheaves.

For uniqueness, if $\widetilde{W}$ is another filtration of
$\mathscr F_0$ by certain subsheaves $\widetilde{\mathscr
F}_0^i$ with desired property, then their restrictions to
$X_0'$ are both equal to the filtration $W'$ by punctual
weights, which is unique (\ref{T2.4}ii), so they are both
obtained by descending this filtration $W',$ and therefore
they are isomorphic.

For functoriality, let $\mathscr G_0$ be another lisse
$\iota$-mixed sheaf with integer punctual $\iota$-weights, and
let $V$ be its filtration by punctual weights, and let
$\varphi:\mathscr F_0\to\mathscr G_0$ be a morphism. Pulling
$\varphi$ back via $P$ we get a morphism $\varphi':\mathscr
F_0'\to\mathscr G_0'$ on $X_0',$ and the diagram
$$
\xymatrix@C=.9cm{
p_1^*\mathscr F_0' \ar[r]^-{\mu_{\mathscr F_0}}
\ar[d]_-{p_1^*\varphi'} & p_2^*\mathscr F_0'
\ar[d]^-{p_2^*\varphi'} \\
p_1^*\mathscr G_0' \ar[r]_-{\mu_{\mathscr G_0}} &
p_2^*\mathscr G_0'}
$$
commutes. By (\ref{T2.4}ii) we have $\varphi'(\mathscr
F_0'^i)\subset\mathscr G_0'^i,$ and the diagram
$$
\xymatrix@C=.9cm{
p_1^*\mathscr F_0'^i \ar[r]^-{\mu_{\mathscr F_0}}
\ar[d]_-{p_1^*\varphi'} & p_2^*\mathscr F_0'^i
\ar[d]^-{p_2^*\varphi'} \\
p_1^*\mathscr G_0'^i \ar[r]_-{\mu_{\mathscr G_0}} &
p_2^*\mathscr G_0'^i}
$$
commutes for each $i.$ Let $\varphi'^i:\mathscr
F_0'^i\to\mathscr G_0'^i$ be the restriction of $\varphi'.$
Then they descend to morphisms $\varphi^i:\mathscr
F_0^i\to\mathscr G_0^i,$ which are restrictions of $\varphi.$
\end{proof}

\begin{subremark}\label{R2.7}
One can prove a similar structure theorem of $\iota$-mixed
sheaves on algebraic stacks over $\bb F_q:$ the proof of
(\ref{L2.6}) carries over verbatim to the case of algebraic
stacks, except that for a presentation $X_0'\to\s 
X_0,$ the fiber product $X_0''=X_0'\times_{\s 
X_0}X_0'$ may not be a scheme, so we use the case for
algebraic spaces and replace every ``(\ref{T2.4})" in the
proof by ``(\ref{L2.6})". It turns out that (\ref{T2.4}iii)
also holds for algebraic stacks, as a consequence of the
proof of (\ref{T1.4}). As we will not use these results in
this paper, we do not give the proof in detail here. See 
(\cite{Decom}, 2.1).
\end{subremark}

\begin{proposition}\label{L2.8}
Let $\mathscr X_0$ be an $\mathbb F_q$-algebraic stack, and
let $P:X_0\to\mathscr X_0$ be a presentation (i.e. a smooth
surjection with $X_0$ a scheme). Then a complex $K_0\in
W(\mathscr X_0,\overline{\mathbb Q}_{\ell})$ is
$\iota$-mixed if and only if $P^*K_0$ (resp. $P^!K_0$) is
$\iota$-mixed.
\end{proposition}

\begin{proof}
We consider $P^*K_0$ first. The ``only if" part follows from
(\ref{L2.3}iv). For the ``if" part, since $P^*$ is exact on
sheaves and so $\mathscr H^i(P^*K_0)=P^*\mathscr H^i(K_0),$
we reduce to the case when $K_0=\mathscr F_0$ is a sheaf.
So we assume the sheaf $\mathscr F_0':=P^*\mathscr F_0$ on
$X_0$ is $\iota$-mixed, and want to show $\mathscr F_0$ is
also $\iota$-mixed. The proof is similar to the argument in
(\ref{L2.6}).

Let $X_0''$ be the fiber product
$$
X_0''=\xymatrix@C=.7cm{
X_0\times_{\mathscr X_0}X_0 \ar[r]^-{p_1}
\ar[d]_-{p_2} & X_0 \ar[d]^-P \\
X_0 \ar[r]_-P & \mathscr X_0.}
$$
Then $X_0''$ is an algebraic space of finite type. Applying
(\ref{T2.4}i) to $\mathscr F_0'$ we get a decomposition
$\mathscr F_0'=\bigoplus_{b\in\mathbb R/\mathbb Z}\mathscr
F_0'(b).$ For $j=1,2,$ applying $p_j^*$ we get a decomposition
$$
p_j^*\mathscr F_0'=\bigoplus_{b\in\mathbb R/\mathbb
Z}p_j^*\mathscr F_0'(b),
$$
which is the decomposition of $p_j^*\mathscr F_0'$ according
to the weights mod $\mathbb Z.$ By the functoriality in
(\ref{L2.6}i), the canonical isomorphism $\mu:p_1^*\mathscr
F_0'\to p_2^*\mathscr F_0'$ takes the form $\bigoplus_{b\in
\mathbb R/\mathbb Z}\mu_b,$ where $\mu_b:p_1^*\mathscr
F_0'(b)\to p_2^*\mathscr F_0'(b)$ is an isomorphism satisfying
cocycle condition as $\mu$ does. Therefore the decomposition
of $\mathscr F_0'$ descends to a decomposition $\mathscr
F_0=\bigoplus_{b\in\mathbb R/\mathbb Z}\mathscr F_0(b).$ The
$\iota$-weights of the local Frobenius eigenvalues of
$\mathscr F_0(b)$ at each point of $\mathscr X_0$ are in the
coset $b.$ Next we show that $\mathscr F_0(b)$'s are
$\iota$-mixed.

Replacing $\mathscr F_0$ by a direct summand $\mathscr
F_0(b)$ and then twisting it appropriately, we may assume its
inverse image $\mathscr F_0'$ is $\iota$-mixed with integer
punctual $\iota$-weights. By (\ref{L2.3}v) we can shrink
$\mathscr X_0$ to a nonempty open substack and assume
$\mathscr F_0$ is lisse. Then $\mathscr F_0'$ is also lisse,
and applying (\ref{T2.4}ii) we get a finite increasing
filtration $W'$ of $\mathscr F_0'$ by lisse subsheaves
$\mathscr F_0'^i,$ such that each $\text{Gr}^{W'}_i(\mathscr
F_0')$ is punctually $\iota$-pure of weight $i.$ Pulling back
this filtration via $p_j,$ we get a finite increasing
filtration $p_j^*W'$ of $p_j^*\mathscr F_0',$ and since
$\text{Gr}_i^{p_j^*W'}(p_j^*\mathscr
F_0')=p_j^*\text{Gr}_i^{W'}(\mathscr F_0')$ is punctually
$\iota$-pure of weight $i,$ it is the filtration by punctual
weights given by (\ref{L2.6}ii). By functoriality, the
canonical isomorphism $\mu:p_1^*\mathscr F_0'\to
p_2^*\mathscr F_0'$ maps $p_1^*\mathscr F_0'^i$
isomorphically onto $p_2^*\mathscr F_0'^i,$ satisfying
cocycle condition. Therefore the filtration $W'$ of $\mathscr
F_0'$ descends to a filtration $W$ of $\mathscr F_0,$ and
$P^*\text{Gr}_i^W(\mathscr F_0)=\text{Gr}_i^{W'}(\mathscr
F_0')$ is punctually $\iota$-pure of weight $i.$ Since $\pi$
is surjective, $\text{Gr}_i^W(\mathscr F_0)$ is also
punctually $\iota$-pure of weight $i,$ therefore $\mathscr
F_0$ is $\iota$-mixed.

Next we consider $P^!K_0.$ We know that $P$ is smooth of
relative dimension $d,$ for some function $d:\pi_0(X_0)\to
\mathbb N.$ Let $X_0^0$ be a connected component of $X_0.$
Since $\pi_0(X_0)$ is finite, $X_0^0$ is both open and
closed in $X_0,$ so $f:X_0^0\overset{j}{\to}X_0\overset{P}
{\to}\mathscr X_0$ is smooth of relative dimension
$d(X_0^0).$ Then $P^*K_0$ is $\iota$-mixed if and only if
$f^*K_0=j^*P^*K_0$ is $\iota$-mixed for the inclusion $j$ of
every connected component, if and only if $f^!K_0=f^*K_0
\langle d(X_0^0)\rangle$ is $\iota$-mixed, if and only if
$P^!K_0$ is $\iota$-mixed, since $f^!=j^!P^!=j^*P^!.$
\end{proof}

\begin{subremark}\label{Laff}
As a consequence of Lafforgue's theorem on the Langlands 
correspondence for function fields and a 
Ramanujan-Petersson type of result, one deduces that all 
complexes on any $\mathbb F_q$-algebraic stack is $\iota$-mixed, 
for any $\iota.$ To see this, by (\ref{L2.8}, \ref{L2.3}v,ii), 
we reduce to the case of an irreducible lisse sheaf on a 
smooth (in particular, normal) $\bb F_q$-scheme. By 
(\cite{Del2}, 1.3.6) we reduce to the case where the determinant 
of the lisse sheaf has finite order, and Lafforgue's result 
applies (\cite{Lau}, 1.3). 
In the following, when we want to emphasize the 
assumption of $\iota$-mixedness, we will still write 
$``W_m(\s X_0,\overline{\bb Q}_{\ell})",$ 
although it equals the full category $W(\s 
X_0,\overline{\bb Q}_{\ell}).$
\end{subremark}

Next we show the stability of $\iota$-mixedness, first for a
few operations on complexes on algebraic spaces, and then 
for all the six operations on stacks. Denote by 
$D_{\s X_0}$ or just $D$ the dualizing functor $R\s 
Hom(-,K_{\s X_0}),$ where $K_{\s X_0}$ is a
dualizing complex on $\s X_0$ (\cite{LO2}, $\S7$).

\begin{blank}\label{R2.9}
Recall (\cite{KW}, II 12.2) that, for $\bb F_q$-schemes 
and bounded complexes of sheaves on them, the operations 
$f_*,f_!,f^*,f^!,D$ and $-\otimes^L-$ all preserve 
$\iota$-mixedness. Since we are working with $\overline
{\bb Q}_{\ell}$-coefficients, $\otimes^L=\otimes.$ 
\end{blank}

\begin{lemma}\label{L2.10}
Let $f:X_0\to Y_0$ be a morphism of $\bb 
F_q$-algebraic spaces. Then the operations $-\otimes-,
D_{X_0}, f_*$ and $f_!$ all preserve $\iota$-mixedness,
namely, they induce functors
\begin{gather*}
-\otimes-:W^-_m(X_0,\overline{\bb Q}_{\ell})\times
W^-_m(X_0,\overline{\bb Q}_{\ell})\longrightarrow
W^-_m(X_0,\overline{\bb Q}_{\ell}), \\
D:W_m(X_0,\overline{\bb Q}_{\ell})\longrightarrow
W_m(X_0,\overline{\bb Q}_{\ell})^{\emph{op}}, \\
f_*:W^+_m(X_0,\overline{\bb Q}_{\ell})\longrightarrow
W^+_m(Y_0,\overline{\bb Q}_{\ell})\quad\emph{and}\quad
f_!:W^-_m(X_0,\overline{\bb Q}_{\ell})\longrightarrow
W^-_m(Y_0,\overline{\bb Q}_{\ell}).
\end{gather*}
\end{lemma}

\begin{proof}
We will reduce to the case of unbounded complexes on schemes,
and then prove the scheme case. Let
$P:X_0'\to X_0$ be an \'etale presentation.

\textbf{Reduction for $\otimes$.} For $K_0,L_0\in
W^-_m(X_0,\overline{\bb Q}_{\ell}),$ we have
$P^*(K_0\otimes L_0)=(P^*K_0)\otimes(P^*L_0),$ and the
reduction follows from (\ref{L2.8}).

\textbf{Reduction for $D$.} For $K_0\in W_m(X_0,\overline{
\bb Q}_{\ell}),$ we have $P^*DK_0=DP^!K_0,$ so the
reduction follows from (\ref{L2.8}).

\textbf{Reduction for $f_*$ and $f_!$.} By definition
(\cite{LO2}, 9.1) we have $f_*=Df_!D,$ so it suffices to
prove the case for $f_!.$ Let $K_0\in
W^-_m(X_0,\overline{\bb Q}_{\ell}),$ and let
$P':Y_0'\to Y_0$ and $X_0'\to X_0\times_{Y_0}Y_0'$ be
\'etale presentations:
$$
\xymatrix@C=1cm{
X_0' \ar @/^2pc/[rr]^-g \ar[r] \ar[rd]_-P & (X_0)_{Y_0'}
\ar[r]_-{f'} \ar[d]^-h & Y_0' \ar[d]^-{P'} \\
& X_0 \ar[r]^-f & Y_0.}
$$
By smooth base change (\cite{LO2}, 12.1) we have
$P'^*f_!K_0=f'_!h^*K_0.$ Replacing $f$ by $f'$ we can
assume $Y_0$ is a scheme. Let $j:U_0\to X_0$ be an open dense
subscheme (\cite{Knu}, II 6.7), with complement $i:Z_0\to
X_0.$ Applying $f_!$ to the exact triangle
$$
\xymatrix@C=.5cm{
j_!j^*K_0 \ar[r] & K_0 \ar[r] & i_*i^*K_0 \ar[r] &}
$$
we get
$$
\xymatrix@C=.5cm{
(fj)_!j^*K_0 \ar[r] & f_!K_0 \ar[r] & (fi)_!i^*K_0 \ar[r] &.}
$$
By (\ref{L2.3}iii) and noetherian induction, we can replace
$X_0$ by $U_0,$ and reduce to the case where $f$ is a
morphism between schemes.

This finishes the reduction to the case of unbounded
complexes on schemes, and now we prove this case.

For the Verdier dual $D_{X_0},$ since the dualizing complex
$K_{X_0}$ has finite quasi-injective dimension, for every
$K_0\in W_m(X_0,\overline{\bb Q}_{\ell})$ and every
integer $i,$ there exist integers $a$ and $b$ such that
$$
\s H^i(D_{X_0}K_0)\simeq\s H^i(D_{X_0}\tau_{[a,b]}K_0),
$$
and by (\ref{R2.9}), we see that $D_{X_0}K_0$ is
$\iota$-mixed.

Next we prove the case of $\otimes.$ For $K_0$ and 
$L_0\in W^-_m(X_0,\overline{\bb Q}_{\ell}),$ we have 
$$
\s H^r(K_0\otimes L_0)= 
\bigoplus_{i+j=r}\s H^i(K_0)\otimes\s H^j(L_0). 
$$
The result follows from (\ref{R2.9}).

Finally we prove the case of $f_*$ and $f_!.$ Let $K_0\in
W^+_m(X_0,\overline{\bb Q}_{\ell}).$ Then we have the
spectral sequence
$$
E_2^{ij}=R^if_*(\s H^jK_0)\Longrightarrow R^{i+j}f_*K_0,
$$
and the result follows from (\ref{R2.9}) and (\ref{L2.3}i,
ii). The case for $f_!=Df_*D$ also follows.
\end{proof}

Finally we prove the main result of this section. This
generalizes (\cite{Beh2}, 6.3.7).

\begin{theorem}\label{T2.11}
Let $f:\s X_0\to\s Y_0$ be a morphism of
$\bb F_q$-algebraic stacks. Then the operations
$f_*,f_!,f^*,f^!,D_{\s X_0},-\otimes-$ and
$R\s Hom(-,-)$ all preserve $\iota$-mixedness,
namely, they induce functors
\begin{gather*}
f_*:W_m^+(\s X_0,\overline{\bb Q}_{\ell})
\longrightarrow W_m^+(\s Y_0,\overline{\bb
Q}_{\ell}), \qquad
f_!:W_m^-(\s X_0,\overline{\bb Q}_{\ell})
\longrightarrow W_m^-(\s Y_0,\overline{\bb 
Q}_{\ell}), \\
f^*:W_m(\s Y_0,\overline{\bb Q}_{\ell})
\longrightarrow W_m(\s X_0,\overline{\bb 
Q}_{\ell}), \qquad
f^!:W_m(\s Y_0,\overline{\bb Q}_{\ell})
\longrightarrow W_m(\s X_0,\overline{\bb 
Q}_{\ell}), \\
D:W_m(\s X_0,\overline{\bb Q}_{\ell})
\longrightarrow W_m(\s X_0,\overline{\bb 
Q}_{\ell})^{\emph{op}}, \\
\otimes:W_m^-(\s X_0,\overline{\bb Q}_{\ell})
\times W_m^-(\s X_0,\overline{\bb Q}_{\ell})
\longrightarrow W_m^-(\s 
X_0,\overline{\bb Q}_{\ell})\quad\emph{and} \\
R\s Hom(-,-):W_m^-(\s X_0,\overline{\bb 
Q}_{\ell})^{\emph{op}}\times W_m^+(\s 
X_0,\overline{\bb Q}_{\ell}) \longrightarrow
W_m^+(\s X_0,\overline{\bb Q}_{\ell}).
\end{gather*}
\end{theorem}

\begin{proof}
Recall from (\cite{LO2}, 9.1) that $f_!:=Df_*D$ and
$f^!:=Df^*D.$ By (\cite{LO2}, 6.0.12, 7.3.1), for $K_0\in
W^-(\s X_0,\overline{\bb Q}_{\ell})$ and 
$L_0\in W^+(\s X_0,\overline{\bb Q}_{\ell}),$ 
we have 
\begin{equation*}
\begin{split}
D(K_0\otimes DL_0) &=R\s Hom(K_0\otimes DL_0,K_{\s 
X_0})=R\s Hom(K_0,R\s Hom(DL_0,K_{\s X_0})) \\
&=R\s Hom(K_0,DDL_0)=R\s Hom(K_0,L_0).
\end{split}
\end{equation*}
Therefore it suffices to prove the result for $f_*,f^*,D$ and
$-\otimes-.$ The case of $f^*$ is proved in (\ref{L2.3}iv).

For $D:$ let $P:X_0\to\s X_0$ be a presentation. Since
$P^*D=DP^!,$ the result follows from (\ref{L2.8}) and
(\ref{L2.10}).

For $\otimes:$ since
$P^*(K_0\otimes L_0)=P^*K_0\otimes P^*L_0,$ the result
follows from (\ref{L2.8}) and (\ref{L2.10}).

For $f_*$ and $f_!:$ we will start with $f_!,$ in order to
use smooth base change to reduce to the case when $\s 
Y_0$ is a scheme, and then turn to $f_*$ in order to use
cohomological descent.

Let $K_0\in W^-_m(\s X_0,\overline{\bb
Q}_{\ell}),$ and let $P:Y_0\to\s Y_0$ be a
presentation and the following diagram be 2-Cartesian:
$$
\xymatrix@C=.8cm{
(\s X_0)_{Y_0} \ar[r]^-{f'} \ar[d]_-{P'} & Y_0
\ar[d]^P \\
\s X_0 \ar[r]^-f & \s Y_0.}
$$
We have (\cite{LO2}, 12.1) that $P^*f_!K_0=f'_!P'^*K_0,$ so by
(\ref{L2.8}) we can assume $\mathscr Y_0=Y_0$ is a scheme.

Now we switch to $f_*,$ where $f:\s X_0\to Y_0,$ and
$K_0\in W_m^+(\s X_0,\overline{\bb Q}_{\ell}).$
Let $X_0\to\s X_0$ be a presentation. Then it
gives a strictly simplicial smooth hypercover 
$X_{0,\bullet}$ of $\s X_0:$
$$
X_{0,n}:=\underbrace{X_0\times_{\s 
X_0}\cdots\times_{\s X_0}X_0}_{n+1\text{\ factors}},
$$
where each $X_{0,n}$ is an $\bb F_q$-algebraic space of
finite type. Let $f_n:X_{0,n}\to Y_0$ be the restriction of
$f$ to $X_{0,n}.$ Then we have the spectral sequence
(\cite{LO2}, 10.0.9)
$$
E_1^{ij}=R^jf_{i*}(K_0|_{X_{0,i}})\Longrightarrow
R^{i+j}f_*K_0.
$$
Since $f_i$'s are morphisms of algebraic spaces, the result
follows from (\ref{L2.10}) and (\ref{L2.3}i, ii).
\end{proof}

\begin{remark}\label{mixed-variant}
In fact, we can take the dualizing complex $K_{\s X_0}$ 
to be \textit{mixed}, and results in this section hold 
(and can be proved verbatim) for \textit{mixed} complexes. 
In particular, mixedness is preserved by the six operations 
and the Verdier dualizing functor for stacks (if we take a 
mixed dualizing complex).
\end{remark}

\section{Stratifiable complexes}

In this section, we use the same notations and hypotheses in
(\ref{adic-setting}). For the purpose of this article, it
suffices to take $S$ to be $\text{Spec }k$ for an
algebraically closed field $k$ of characteristic not equal to
$\ell,$ but we want to work in the general setting (namely,
any scheme that satisfies (LO)) for future applications; for 
instance see \cite{Decom}. Let
$\mathcal{X,Y}\cdots$ be $S$-algebraic stacks of finite type.
By ``sheaves" we mean ``lisse-\'etale sheaves". 
``Jordan-H\"older" and ``locally constant constructible" are
abbreviated as ``JH" and ``lcc" respectively. A
\textit{stratification} $\s S$ of an $S$-algebraic stack
$\mathcal X$ is a finite set of disjoint
locally closed substacks that cover $\mathcal X.$ If
$\s F$ is a lcc $(\Lambda_n)_{\mathcal X}$-module, a
\textit{decomposition series} of $\s F$ is a
filtration by lcc $\Lambda_{\mathcal X}$-subsheaves, such
that the successive quotients are simple
$\Lambda_{\mathcal X}$-modules. Note that the filtration is
always finite, and the simple successive quotients, which
are $(\Lambda_0)_{\mathcal X}$-modules, are independent
(up to order) of the decomposition series chosen. They are
called the \textit{JH components of $\s F.$}

\begin{definition}\label{D3.1}
(i) A complex $K=(K_n)_n\in\s D_c(\s A)$ is said
to be \emph{stratifiable}, if there exists a pair $(\s 
S,\mathcal L),$ where $\s S$ is a stratification of
$\mathcal X,$ and $\mathcal L$ is a function that assigns to
every stratum $\mathcal U\in\s S$ a finite set
$\mathcal{L(U)}$ of isomorphism classes of simple (i.e.
irreducible) lcc $\Lambda_0$-modules on $\mathcal
U_{\emph{lis-\'et}},$ such that for each pair $(i,n)$ of
integers, the restriction of the sheaf $\s H^i(K_n)\in
\emph{Mod}_c(\mathcal X_{\emph{lis-\'et}},\Lambda_n)$ to
each stratum $\mathcal U\in\s S$ is lcc, with JH
components (as a $\Lambda_{\mathcal U}$-module) contained in
$\mathcal{L(U)}.$ We say that the pair $(\s S,\mathcal
L)$ \emph{trivializes} $K$ (or $K$ is $(\s 
S,\mathcal L)$\emph{-stratifiable}), and denote the full
subcategory of $(\s S,\mathcal L)$-stratifiable
complexes by $\s D_{\s S,\mathcal L}(\s 
A).$ The full subcategory of 
stratifiable complexes in $\s D_c(\s A)$ is
denoted by $\s D_c^{\emph{stra}}(\s A).$

(ii) Let $D_c^{\emph{stra}}(\mathcal X,\Lambda)$ be the
essential image of $\s D_c^{\emph{stra}}(\s A)$
in $D_c(\mathcal X,\Lambda),$ and we call the objects of
$D_c^{\emph{stra}}(\mathcal X,\Lambda)$ \emph{stratifiable
complexes of sheaves.}

(iii) Let $E_{\lambda}$ be a finite extension of $\bb 
Q_{\ell}$ with ring of integers $\s O_{\lambda}.$ Then
the definition above applies to $\Lambda=\s 
O_{\lambda}.$ Let $D_c^{\emph{stra}}(\mathcal X,E_{\lambda})$
be the essential image of $D_c^{\emph{stra}}(\mathcal
X,\s O_{\lambda})$ in $D_c(\mathcal X,E_{\lambda}).$
Finally we define
$$
D_c^{\emph{stra}}(\mathcal X,\overline{\bb Q}_{\ell})=
\emph{2-colim}_{E_{\lambda}}D_c^{\emph{stra}}(\mathcal X,
E_{\lambda}).
$$
\end{definition}

\begin{subremark}\label{R3.2}
(i) This notion is due to Beilinson, Bernstein and Deligne
\cite{BBD}, and Behrend \cite{Beh2} used it to define his
derived category for stacks. Many results in this section are
borrowed from \cite{Beh2}, but reformulated and reproved in
terms of the derived categories defined in \cite{LO2}.

(ii) Let $\s F$ be a $\Lambda_n$-sheaf trivialized by 
a pair $(\s S,\mathcal L),$ and let $\s G$ be a 
sub-quotient sheaf of $\s F.$ Then $\s G$ is not 
necessarily trivialized by $(\s S,\mathcal L).$ But if 
$\s G$ is lcc on each stratum in $\s S,$ then it 
is necessarily trivialized by $(\s S,\mathcal L).$
\end{subremark}

\begin{blank}\label{refinement}
We say that the pair $(\mathscr S',\mathcal L')$
\textit{refines the pair} $(\mathscr S,\mathcal L),$ if
$\mathscr S'$ refines
$\mathscr S,$ and for every $\mathcal V\in\mathscr S',\
\mathcal U\in\mathscr S$ and $L\in\mathcal{L(U)},$ such that
$\mathcal V\subset\mathcal U,$ the restriction $L|_{\mathcal
V}$ is trivialized by $\mathcal{L'(V)}.$ Given a pair
$(\mathscr S,\mathcal L)$ and a refined stratification
$\mathscr S'$ of $\mathscr S,$ there is a canonical way to
define $\mathcal L'$ such that $(\mathscr S',\mathcal L')$
refines $(\mathscr S,\mathcal L):$ for every $\mathcal
V\in\mathscr S',$ we take $\mathcal{L'(V)}$ to be the set of
isomorphism classes of JH components of the lcc sheaves
$L|_{\mathcal V}$ for $L\in\mathcal{L(U)},$ where $\mathcal U$
ranges over all strata in $\mathscr S$ that contains $\mathcal
V.$ It is clear that the set of all pairs $(\mathscr
S,\mathcal L)$ form a filtered direct system.

A pair $(\mathscr S,\mathcal L)$ is said to be
\textit{tensor closed} if for every $\mathcal U\in\mathscr
S$ and $L,M\in\mathcal{L(U)},$ the sheaf tensor product
$L\otimes_{\Lambda_0}M$ has JH components in $\mathcal{L(U)}.$

For a pair $(\mathscr S,\mathcal L),$ a \textit{tensor 
closed hull} of this pair is a tensor closed refinement.
\end{blank}

\begin{lemma}\label{L3.2.5}
Every pair $(\mathscr S,\mathcal L)$ can be refined to a
tensor closed pair $(\mathscr S',\mathcal L').$
\end{lemma}

\begin{proof}
First we show that, for a lcc sheaf of sets $\mathscr F$ on
$\mathcal X_{\text{lis-\'et}},$ there exists a finite \'etale
morphism $f:\mathcal Y\to\mathcal X$ of algebraic $S$-stacks
such that $f^{-1}\mathscr F$ is constant. Consider the total
space $[\mathscr F]$ of the sheaf $\mathscr F.$ Precisely,
this is the category fibered in groupoids over
$(\text{Aff}/S)$ with the underlying category described as
follows. Its objects are triples
$(U\in\text{obj(Aff}/S),u\in\text{obj }\mathcal
X(U),s\in(u^{-1}\mathscr F)(U)),$ and morphisms from
$(U,u,s)$ to $(V,v,t)$ are pairs $(f:U\to
V,\alpha:vf\Rightarrow u)$ such that $t$ is mapped to $s$
under the identification $\alpha:f^{-1}v^{-1}\mathscr F\cong
u^{-1}\mathscr F.$ The map $(U,u,s)\mapsto(U,u)$ gives a map
$g:[\mathscr F]\to\mathcal X,$ which is representable finite
\'etale (because it is so locally). The pullback sheaf
$g^{-1}\mathscr F$ on $[\mathscr F]$ has a global section,
so the total space breaks up into two parts, one part being
mapped isomorphically onto the base $[\mathscr F].$ By
induction on the degree of $g$ we are done.

Next we show that, for a fixed representable finite \'etale
morphism $\mathcal Y\to\mathcal X,$ there are only finitely
many isomorphism classes of simple lcc $\Lambda_0$-sheaves on
$\mathcal X$ that become constant when pulled back to
$\mathcal Y.$ We can assume that both $\mathcal X$ and 
$\mathcal Y$ are connected. By the following lemma 
(\ref{Gal}), we reduce to the case where $\mathcal 
Y\to\mathcal X$ is Galois with group $G,$ for some finite 
group $G.$
Then simple lcc $\Lambda_0$-sheaves on
$\mathcal X$ that become constant on $\mathcal Y$
correspond to simple left $\Lambda_0[G]$-modules, which are
cyclic and hence isomorphic to $\Lambda_0[G]/I$ for left
maximal ideals $I$ of $\Lambda_0[G].$ There are only finitely
many such ideals since $\Lambda_0[G]$ is a finite set.

Also note that, a lcc subsheaf of a constant constructible
sheaf on a connected stack is also constant. Let $L$ be a
lcc subsheaf on $\mathcal X$ of the constant sheaf
associated to a finite set $M.$ Consider their total spaces.
We have an inclusion of substacks
$i:[L]\hookrightarrow\coprod_{m\in M}\mathcal X_m,$ where
each part $\mathcal X_m$ is identified with $\mathcal X.$
Then $i^{-1}(\mathcal X_m)\to\mathcal X_m$ is finite \'etale,
and is the inclusion of a substack, hence is either an
equivalence or the inclusion of the empty substack, since
$\mathcal X$ is connected. It is clear that $L$ is also
constant, associated to the subset of those $m\in M$ for
which $i^{-1}(\mathcal X_m)\ne\emptyset.$

Finally we prove the lemma. Refining $\mathscr S$ if
necessary, we assume all strata are connected stacks. For
each stratum $\mathcal U\in\mathscr S,$ let $\mathcal
Y\to\mathcal U$ be a representable finite \'etale morphism,
such that all sheaves in $\mathcal{L(U)}$ become constant on
$\mathcal Y.$ Then define $\mathcal{L'(U)}$ to be the set of
isomorphism classes of simple lcc $\Lambda_0$-sheaves on
$\mathcal U_{\text{lis-\'et}}$ which become constant on
$\mathcal Y.$ For any $L$ and $M\in\mathcal{L'(U)},$ since
all lcc subsheaves of $L\otimes_{\Lambda_0}M$ are constant on
$\mathcal Y,$ we see that $L\otimes_{\Lambda_0}M$ has JH
components in $\mathcal{L'(U)}$ and hence $(\mathscr
S,\mathcal L')$ is a tensor closed refinement of $(\mathscr
S,\mathcal L).$
\end{proof}

\begin{sublemma}\label{Gal}
Let $\mathcal Y\to\mathcal X$ be a representable finite 
\'etale morphism between connected $S$-algebraic stacks. 
Then there exists a morphism 
$\mathcal Z\to\mathcal Y,$ such that $\mathcal Z$ is 
\emph{Galois} over $\mathcal X,$ i.e. it is a $G$-torsor 
for some finite group $G.$ 
\end{sublemma}

\begin{proof}
Assume $\mathcal X$ is non-empty, and take a geometric 
point $\overline{x}\to\mathcal X.$ Let $\mathscr C$ be the 
category $\text{F\'Et}(\mathcal X)$ of representable finite 
\'etale morphisms to $\mathcal X,$ and let 
$$
F:\mathscr C\to\text{FSet} 
$$
be the fiber functor to the category of finite sets, namely 
$F(\mathcal Y)=Hom_{\mathcal X}(\overline{x},\mathcal Y).$ 
Note that this Hom, which is a priori a category, is a 
finite set, since $\mathcal Y\to\mathcal X$ is 
representable and finite. Then one can verify that 
$(\mathscr C,F:\mathscr C\to\text{FSet})$ satisfies the 
axioms of Galois formalism in (\cite{SGA1}, Exp. V, 4), and 
use the consequence g) on p. 121 in \textit{loc. cit.} For the 
reader's convenience, we follow Olsson's suggestion and 
explain the proof briefly. Basically, we will 
verify certain axioms of (G1) -- (G6), and deduce the 
conclusion as in \textit{loc. cit.} 

First note that $\mathscr C,$ which is a priori a 
2-category, is a 1-category. This is because for any 
2-commutative diagram 
$$
\xymatrix@C=.6cm @R=.5cm{
\mathcal Y \ar[dr] \ar[rr]^-f && \mathcal Z \ar[dl] \\ 
& \mathcal X &}
$$
where $\mathcal Y,\mathcal Z\in\mathscr C,$ the morphism 
$f$ is also representable (and finite \'etale), so 
$Hom_{\mathcal X}(\mathcal Y,\mathcal Z)$ is discrete. By 
definition, the functor $F$ preserves fiber-products, and 
$F(\mathcal X)$ is a one-point set. 

Let $f:\mathcal Y\to\mathcal Z$ be a morphism in $\mathscr 
C,$ then it is finite \'etale. So if the degree of $f$ is 
1, then $f$ is an isomorphism. This implies that the 
functor $F$ is \textit{conservative}, i.e. $f$ is an 
isomorphism if $F(f)$ is. In particular, $f$ is a 
monomorphism if and only if $F(f)$ is. This is because $f$ 
is a monomorphism if and only if $p_1:\mathcal 
Y\times_{\mathcal Z}\mathcal Y\to\mathcal Y$ is an 
isomorphism, and $F$ preserves fiber-products. 

Since $f:\mathcal Y\to\mathcal Z$ is finite \'etale, its 
image stack $\mathcal Y'\subset\mathcal Z$ is both open and 
closed, hence $\mathcal Y'\to\mathcal Z$ is a monomorphism 
that is an isomorphism onto a direct summand of $\mathcal 
Z$ (i.e. $\mathcal Z=\mathcal Y'\coprod\mathcal Y''$ for 
some other open and closed substack $\mathcal Y''\subset 
\mathcal Z$). Also, since $\mathcal Y\to\mathcal Y'$ is epic 
and finite \'etale, it is \textit{strictly epic}, i.e. for 
every $\mathcal Z\in\mathscr C,$ the diagram 
$$
Hom(\mathcal Y',\mathcal Z)\to Hom(\mathcal 
Y,\mathcal Z)\rightrightarrows Hom(\mathcal 
Y\times_{\mathcal Y'}\mathcal Y,\mathcal Z)
$$
is an equalizer. 

Every object $\mathcal Y$ in $\mathscr C$ is artinian: for 
a chain of monomorphisms 
$$
\cdots\to\mathcal Y_n\to\cdots\to\mathcal Y_2\to\mathcal 
Y_1\to\mathcal Y,
$$
we get a chain of injections 
$$
\cdots\to F(\mathcal Y_n)\to\cdots\to F(\mathcal 
Y_1)\to F(\mathcal Y), 
$$
which is stable since $F(\mathcal Y)$ is a finite set, and 
so the first chain is also stable since $F$ is 
conservative. 

Since $F$ is left exact and every object in $\mathscr C$ is 
artinian, by (\cite{Gro2}, 3.1) the functor $F$ is 
\textit{strictly pro-representable}, i.e. there exists a 
projective system $P=\{P_i;i\in I\}$ of objects in $\mathscr 
C$ indexed by a filtered partially ordered set $I,$ with 
epic transition morphisms $\varphi_{ij}:P_j\to P_i\ (i\le 
j),$ such that there is a natural isomorphism of functors 
$$
F\overset{\sim}{\longrightarrow}Hom(P,-):=\text{colim}_IHom 
(P_i,-).
$$
Let $\psi_i:P\to P_i$ be the canonical projection in the 
category $\text{Pro}(\mathscr C)$ of pro-objects of 
$\mathscr C.$ We may assume that every epimorphism 
$P_j\to\mathcal Z$ in $\mathscr C$ is isomorphic to 
$P_j\overset{\varphi_{ij}}{\to}P_i$ for some $i\le j.$ This 
is because one can add 
$P_j\to\mathcal Z$ into the projective system $P$ without 
changing the functor it represents. Also one can show that 
the $P_i$'s are connected (cf. loc. cit.), and morphisms in 
$\mathscr C$ between connected stacks are strictly epic. 

Given $\mathcal Y\in\mathscr C,$ now we show that there 
exists an object $\mathcal Z\to\mathcal X$ that is Galois 
and factors through $\mathcal Y.$ Since $F(\mathcal Y)$ is a 
finite set, there exists an index $j\in I$ such that all 
maps $P\to\mathcal Y$ factors through $P\overset{\psi_j} 
{\to}P_j.$ This means that the canonical map 
$$
P\to\mathcal Y^J:=\underbrace{\mathcal Y\times_{\mathcal 
X}\cdots\times_{\mathcal X}\mathcal Y}_{\#J\text{ factors}}, 
\quad\text{where }J:=F(\mathcal 
Y)=Hom_{\text{Pro}(\mathscr C)}(P,\mathcal Y) 
$$
factors as 
$$
\xymatrix@C=.7cm{ 
P \ar[r]^-{\psi_j} & P_j \ar[r]^-A & \mathcal Y^J.}
$$
Let $P_j\to P_i\overset{B}{\to}\mathcal Y^J$ be the 
factorization of $A$ into a composition of an epimorphism 
and a monomorphism $B.$ We claim that $P_i$ is Galois over 
$\mathcal X.$ 

Since $F(P_i)$ is a finite set, there exists an index $k\in 
I$ such that all maps $P\to P_i$ factors through 
$P\overset{\psi_k}{\to}P_k.$ Fix any $v:P_k\to P_i.$ To show 
$P_i$ is Galois, it suffices to show that $\text{Aut}(P_i)$ 
acts on $F(P_i)=Hom(P_k,P_i)$ transitively, i.e. there 
exists a $\sigma\in\text{Aut}(P_i)$ making the triangle 
commute:
$$
\xymatrix@C=1cm @R=.7cm{
P_k \ar[dr]_-{\varphi_{ik}} \ar[r]^-v & P_i \ar[d]^-{\sigma} 
\\ 
& P_i.} 
$$
For every $u\in J=Hom(P_i,\mathcal Y),$ we have $u\circ 
v\in Hom(P_k,\mathcal Y),$ so there exists a $u'\in 
Hom(P_i,\mathcal Y)$ making the diagram commute:
$$
\xymatrix@C=1cm @R=.7cm{
P_k \ar[r]^-v \ar[d]_-{\varphi_{ik}} & P_i \ar[d]^-u \\ 
P_i \ar[r]_-{u'} & \mathcal Y.}
$$
Since $v$ is epic, the function $u\mapsto u':J\to J$ is 
injective, hence a bijection. Let $\alpha:\mathcal Y^J\to 
\mathcal Y^J$ be the isomorphism induced by the map 
$u\mapsto u'.$ Then the diagram 
$$
\xymatrix@C=1cm @R=.7cm{
P_k \ar[r]^-v \ar[dr]_-{\varphi_{ik}} & P_i \ar[r]^-B & 
\mathcal Y^J \ar[d]^-{\alpha} \\ 
& P_i \ar[r]_-B & \mathcal Y^J}
$$
commutes. By the uniqueness of the factorization of the map 
$P_k\to\mathcal Y^J$ into the composition of an epimorphism 
and a monomorphism, there exists a 
$\sigma\in\text{Aut}(P_i)$ such that $\sigma\circ 
v=\varphi_{ik}.$ This finishes the proof.
\end{proof}

We give some basic properties of stratifiable complexes.

\begin{lemma}\label{L3.3}
(i) $\s D_c^{\emph{stra}}(\s A)$ (resp.
$D_c^{\emph{stra}}(\mathcal X,\Lambda)$) is a triangulated
subcategory of $\s D_c(\s A)$ (resp. $D_c(\mathcal 
X,\Lambda)$) with the induced standard $t$-structure.

(ii) If $f:\mathcal X\to\mathcal Y$ is an $S$-morphism,
then $f^*:\s D_c(\s A(\mathcal Y))\to\s D_c(\s 
A(\mathcal X))$ (resp. $f^*:D_c(\mathcal Y,\Lambda)\to
D_c(\mathcal X,\Lambda)$) preserves stratifiability.

(iii) If $\s S$ is a stratification of $\mathcal X,$
then $K\in\s D_c(\s A(\mathcal X))$ is
stratifiable if and only if $K|_V$ is stratifiable for every
$V\in\s S.$

(iv) Let $P:X\to\mathcal X$ be a presentation, and let
$K=(K_n)_n\in\s D_c(\s A(\mathcal X)).$ Then $K$
is stratifiable if and only if $P^*K$ is stratifiable.

(v) $D_c^{\emph{stra}}(\mathcal X,\Lambda)$ contains
$D_c^b(\mathcal X,\Lambda),$ and the heart of
$D_c^{\emph{stra}}(\mathcal X,\Lambda)$ is the same as that
of $D_c(\mathcal X,\Lambda)$ (\ref{R2.2}i).

(vi) Let $K\in\s D_c(\s A)$ be a normalized
complex (\cite{LO2}, 3.0.8). Then $K$ is trivialized by a
pair $(\s S,\mathcal L)$ if and only if $K_0$ is
trivialized by this pair.

(vii) Let $K\in\s D_c^{\emph{stra}}(\s A).$ Then
its Tate twist $K(1)$ is also stratifiable.
\end{lemma}

\begin{proof}
(i) To show $\s D_c^{\text{stra}}(\s A)$ is a
triangulated subcategory, it suffices to show (\cite{SGA4.5},
p.271) that for every exact triangle $K'\to K\to K''\to K'[1]$
in $\s D_c(\s A),$ if $K'$ and $K''$ are
stratifiable, so also is $K.$

Using refinement we may assume that $K'$ and $K''$ are
trivialized by the same pair $(\s S,\mathcal L).$
Consider the cohomology sequence of this exact triangle at
level $n,$ restricted to a stratum $\mathcal U\in\s S.$
By (\cite{Ols3}, 9.1), to show that a sheaf is lcc on
$\mathcal U,$ one can pass to a presentation $U$ of the stack
$\mathcal U.$ Then by (\cite{Mil1}, 20.3) and 5-lemma, we see
that the $\s H^i(K_n)$'s are lcc on $\mathcal U,$ with
JH components contained in $\mathcal{L(U)}.$ Therefore
$\s D_c^{\text{stra}}(\s A)$ (and hence
$D_c^{\text{stra}}(\mathcal X,\Lambda)$) is a triangulated
subcategory.

The $t$-structure is inherited by $\s
D_c^{\text{stra}}(\s A)$ (and hence by
$D_c^{\text{stra}}(\mathcal X,\Lambda)$) because, if
$K\in\s D_c(\s A)$ is stratifiable, so also
are its truncations $\tau_{\le r}K$ and $\tau_{\ge r}K.$

(ii) $f^*$ is exact on the level of sheaves, and takes a lcc
sheaf to a lcc sheaf. If $(K_n)_n\in\s D_c(\s 
A(\mathcal Y))$ is trivialized by $(\s S,\mathcal L),$
then $(f^*K_n)_n$ is trivialized by $(f^*\s 
S,f^*\mathcal L),$ where $f^*\s 
S=\{f^{-1}(V)|V\in\s S\}$ and $(f^*\mathcal
L)(f^{-1}(V))$ is the set of isomorphism classes of JH
components of $f^*L,\ L\in\mathcal L(V).$ The case of
$D_c(-,\Lambda)$ follows easily.

(iii) The ``only if" part follows from (ii). The ``if" part
is clear: if $(\s S_V,\mathcal L_V)$ is a pair on $V$
that trivializes $(K_n|_V)_n,$ then the pair $(\s 
S_{\mathcal X},\mathcal L)$ on $\mathcal X,$ where
$\s S_{\mathcal X}=\cup\s S_V$ and $\mathcal
L=\{\mathcal L_V\}_{V\in\s S},$ trivializes $(K_n)_n.$

(iv) The ``only if" part follows
from (ii). For the ``if" part, assume $P^*K$ is trivialized
by a pair $(\s S_X,\mathcal L_X)$ on $X.$ Let
$U\in\s S_X$ be an open stratum, and let
$V\subset\mathcal X$ be the image of $U$ (\cite{LMB}, 3.7).
Recall that for every $T\in\text{Aff}/S,\ V(T)$ is the full
subcategory of $\mathcal X(T)$ consisting of objects $x$ that
are locally in the essential image of $U(T),$ i.e. such that
there exists an \'etale surjection $T'\to T$ in
$\text{Aff}/S$ and $u'\in U(T'),$ such that the image of $u'$
in $\mathcal X(T')$ and $x|_{T'}$ are isomorphic. Then $V$ is
an open substack of $\mathcal X$ (hence also an algebraic
stack) and $P|_U:U\to V$ is a presentation. Replacing
$P:X\to\mathcal X$ by $P|_U:U\to V$ and using noetherian
induction and (iii), we may assume $\s S_X=\{X\}.$

It follows from a theorem of Gabber \cite{Gab} that $P_*$ 
takes a bounded complex to a bounded complex. In fact, using 
base change by $P,$ we may assume that $P:Y\to X$ is a morphism 
from an $S$-algebraic space $Y$ to an $S$-scheme $X.$ Let 
$j:U\to Y$ be an open dense subscheme of $Y$ with complement 
$i:Z\to Y.$ For a bounded complex $L$ of $\Lambda_n$-sheaves 
on $Y,$ we have the exact triangle 
$$
\xymatrix@C=.5cm{
(Pi)_*i^!L \ar[r] & P_*L \ar[r] & (Pj)_*j^*L \ar[r] &.}
$$
Gabber's theorem implies that $(Pj)_*j^*L$ is bounded, 
since $Pj:U\to X$ is a morphism between schemes. Note that 
the dualizing functor preserves boundedness, so does 
$i^!=D_Zi^*D_Y,$ and therefore we may assume that 
$(Pi)_*i^!L$ is bounded by noetherian induction. It follows 
that $P_*L$ is bounded.

Now take a pair $(\mathscr S,\mathcal L)$ on $\mathcal X$ that
trivializes all $P_*L$'s, for $L\in\mathcal L_X;$ this is 
possible since each $P_*L$ is bounded and $\mathcal L_X$ 
is a finite set. We 
claim that $K$ is trivialized by $(\s S,\mathcal L).$

For each sheaf $\s F$ on $\mathcal X,$ the natural map
$\s F\to R^0P_*P^*\s F$ is injective. This follows 
from the sheaf axioms for the lisse-lisse topology, and 
the fact that the lisse-\'etale topos and the lisse-lisse 
topos are the same. Explicitly, to verify the injectivity 
on $X_U\to U,$ for any $u\in\mathcal X(U),$ since 
the question is \'etale local on $U,$ one can assume
$P:X_U\to U$ has a section $s:U\to X_U.$ Then the
composition $\s F_U\to R^0P_*P^*\s F_U\to
R^0P_*R^0s_*s^*P^*\s F_U=\s F_U$ of the two
adjunctions is the adjunction for $P\circ s=\text{id},$ so
the composite is an isomorphism, and the first map is
injective. 

We take $\s F$ to be the cohomology sheaves
$\s H^i(K_n).$ Since $P^*\s H^i(K_n)$ is an 
iterated extension of sheaves in $\mathcal L_X,$ we see that 
$P_*P^*\s H^i(K_n),$ and in particular $R^0P_*P^*
\s H^i(K_n),$ are trivialized by $(\s S,\mathcal L)$ by (i). 
Since $\s H^i(K_n)$ is lcc (\cite{Ols3}, 9.1), by 
(\ref{R3.2}ii) we see that $\s H^i(K_n)$ (hence $K$) is 
trivialized by $(\s S,\mathcal L).$

(v) It suffices to show, by (i) and (\ref{R2.2}i), that all
adic systems $M=(M_n)_n\in\s A$ are stratifiable. By
(iv) we may assume $\mathcal X=X$ is an $S$-scheme. Since $X$
is noetherian, there exists a stratification
(\cite{SGA5}, VI, 1.2.6) of $X$ such that $M$ is lisse on
each stratum. By (iii) we may assume $M$ is lisse on $X.$

Let $\mathcal L$ be the set of isomorphism classes of JH
components of the $\Lambda_0$-sheaf $M_0.$ We claim that
$\mathcal L$ trivializes $M_n$ for all $n.$ Suppose it
trivializes $M_{n-1}$ for some $n\ge1.$ Consider the
sub-$\Lambda_n$-modules $\lambda M_n\subset M_n[\lambda^n]
\subset M_n,$ where $M_n[\lambda^n]$ is the kernel of the
map $\lambda^n:M_n\to M_n.$ Since $M$ is adic, we have
exact sequences of $\Lambda_X$-modules
\begin{gather*}
\xymatrix@C=.5cm{
0 \ar[r] & \lambda M_n \ar[r] & M_n \ar[r] & M_0 \ar[r] &
0,} \\
\xymatrix@C=.5cm{
0 \ar[r] & M_n[\lambda^n] \ar[r] & M_n \ar[r] &
\lambda^nM_n \ar[r] & 0,}\quad\text{and} \\
\xymatrix@C=.5cm{
0 \ar[r] & \lambda^nM_n \ar[r] & M_n \ar[r] & M_{n-1}
\ar[r] & 0.}
\end{gather*}
The natural surjection $M_n/\lambda M_n\to M_n/M_n[\lambda
^n]$ implies that $\mathcal L$ trivializes $\lambda^nM_n,$
and therefore it also trivializes $M_n.$ By induction on $n$
we are done.

Since $D_c^b\subset D_c^{\text{stra}}\subset D_c,$ and
$D_c^b$ and $D_c$ have the same heart, it is clear that
$D_c^{\text{stra}}$ has the same heart as them.

(vi) Applying $-\otimes^L_{\Lambda_n}K_n$ to the
following exact sequence, viewed as an exact triangle in
$\s D(\mathcal X,\Lambda_n)$
$$
\xymatrix@C=.8cm{
0 \ar[r] & \Lambda_{n-1} \ar[r]^-{1\mapsto\lambda} &
\Lambda_n \ar[r] & \Lambda_0 \ar[r] & 0,}
$$
we get an exact triangle by (\cite{LO2}, 3.0.10)
$$
\xymatrix@C=.5cm{
K_{n-1} \ar[r] & K_n \ar[r] & K_0 \ar[r] &.}
$$
By induction on $n$ and (\ref{R3.4}) below, we see that $K$
is trivialized by $(\mathscr S,\mathcal L)$ if $K_0$ is.

(vii) Let $K=(K_n)_n.$ By definition $K(1)=(K_n(1))_n,$ where
$K_n(1)=K_n\otimes^L_{\Lambda_n}\Lambda_n(1).$ Note that the
sheaf $\Lambda_n(1)$ is a flat $\Lambda_n$-module: to show
that $-\otimes_{\Lambda_n}\Lambda_n(1)$ preserves injections,
one can pass to stalks at geometric points, over which we
have a trivialization $\Lambda_n\simeq\Lambda_n(1).$

Suppose $K$ is $(\s S,\mathcal L)$-stratifiable. Using
the isomorphism 
$$
\s H^i(K_n)\otimes_{\Lambda_n}\Lambda_n(1)=\s 
H^i(K_n\otimes^L_{\Lambda_n}\Lambda_n(1)),
$$
it suffices to show the existence of a pair $(\s 
S,\mathcal L')$ such that for each $\mathcal U\in\s S,$
the JH components of the lcc sheaves $L\otimes_{\Lambda_n}
\Lambda_n(1)$ lie in $\mathcal L'(\mathcal U),$ for all
$L\in\mathcal{L(U)}.$ Since $L$ is a $\Lambda_0$-module, we
have
$$
L\otimes_{\Lambda_n}\Lambda_n(1)=(L\otimes_{\Lambda_n}
\Lambda_0)\otimes_{\Lambda_n}\Lambda_n(1)=L\otimes_{\Lambda_n}
(\Lambda_0\otimes_{\Lambda_n}\Lambda_n(1))=L\otimes_{\Lambda_n}
\Lambda_0(1)=L\otimes_{\Lambda_0}\Lambda_0(1),
$$
and we can take $\mathcal L'(\mathcal U)$ to be a tensor 
closed hull of $\{\Lambda_0(1),L\in\mathcal{L(U)}\}.$
\end{proof}

\begin{subremark}\label{R3.4}
In fact the proof of (\ref{L3.3}i) shows that $\s 
D_{\s S,\mathcal L}(\s A)$ is a triangulated subcategory with 
induced standard $t$-structure, for each fixed pair $(\s 
S,\mathcal L).$ Let $D_{\s S,\mathcal L}(\mathcal
X,\Lambda)$ be the essential image of $\s 
D_{\s S,\mathcal L}(\s A)$ in $D_c(\mathcal
X,\Lambda),$ and this is also a triangulated subcategory 
with induced standard $t$-structure.

Also if $E^{ij}_r\Longrightarrow E^n$ is a
spectral sequence in the category $\s A(\mathcal X),$
and the $E_r^{ij}$'s are trivialized by $(\s S,\mathcal
L)$ for all $i,j,$ then all the $E^n$'s are trivialized by
$(\s S,\mathcal L).$
\end{subremark}

We denote by $D_c^{\dagger,\text{stra}}(\mathcal X,\Lambda),$
for $\dagger=\pm,b,$ the full subcategory of
$\dagger$-bounded stratifiable complexes, using the induced
$t$-structure.

The following is a key result for showing the stability of
stratifiability under the six operations later. Recall
that $M\mapsto\widehat{M}=L\pi^*R\pi_*M$ is the
normalization functor, where $\pi:\mathcal X^{\bb 
N}\to\mathcal X$ is the morphism of topoi in 
(\cite{LO2}, 2.1), mentioned in (\ref{normalization}).

\begin{proposition}\label{L3.4.5}
For a pair $(\s S,\mathcal L)$ on $\mathcal X,$ 
if $M\in\s D_{\s S,\mathcal L}(\s A),$ then
$\widehat{M}\in\s D_{\s S,\mathcal
L}(\s A),$ too. In particular, if $K\in
D_c(\mathcal X,\Lambda),$ then $K$ is stratifiable
if and only if its normalization $\widehat{K}\in\s 
D_c(\s A)$ is stratifiable.
\end{proposition}

\begin{proof}
First, we will reduce to the case where $M$ is essentially 
bounded (i.e. $\s H^iM$ is AR-null for $|i|\gg0$). Let 
$P:X\to\mathcal X$ be a presentation. The 
$\ell$-cohomological dimension of $X_{\text{\'et}}$ is 
finite, by the assumption (LO) on $S.$ Therefore, by 
(\cite{LO2}, 2.1.i), the normalization functor for $X$ 
has finite cohomological dimension, and the same is true 
for $\mathcal X$ since $P^*\widehat{M}=\widehat{P^*M},$ 
by (\cite{LO2}, 2.2.1, 3.0.11). 
This implies that for each integer $i,$ there exist 
integers $a$ and $b$ with $a\le b,$ such that 
$\s H^i(\widehat{M})=\s H^i(\widehat
{\tau_{[a,b]}M}).$ Since $\tau_{[a,b]}M$ is also 
trivialized by $(\s S,\mathcal L),$ we may assume 
$M\in\s D_{\s S,\mathcal L}^b(\s A(\mathcal X)).$ 

Since $\widehat{M}$ is normalized, by (\ref{L3.3}vi), it
suffices to show that $(\widehat{M})_0$ is trivialized by
$(\s S,\mathcal L).$ Using projection formula and the
flat resolution of $\Lambda_0$
$$
\xymatrix@C=.7cm{
0 \ar[r] & \Lambda \ar[r]^-{\lambda} & \Lambda
\ar[r]^-{\epsilon} & \Lambda_0 \ar[r] & 0,}
$$
we have (\cite{LO2}, p.176)
$$
(\widehat{M})_0=\Lambda_0\otimes^L_{\Lambda}R\pi_*M=
R\pi_*(\pi^*\Lambda_0\otimes^L_{\Lambda_{\bullet}}M),
$$
where $\pi^*\Lambda_0$ is the constant projective system
defined by $\Lambda_0.$ Let $C\in\s{D(A)}$ be
the complex of projective systems $\pi^*\Lambda_0\otimes^L
_{\Lambda_{\bullet}}M;$ it is a $\lambda$-complex, and
$C_n=\Lambda_0\otimes^L_{\Lambda_n}M_n\in\s D_c
(\mathcal X,\Lambda_0).$

Recall (\cite{SGA5}, V, 3.2.3) that, a projective system
$(K_n)_n$ ringed by $\Lambda_{\bullet}$ in an abelian
category is AR-adic if and only if

$\bullet$ it satisfies the condition (MLAR) (\cite{SGA5},
V, 2.1.1), hence (ML), and denote by $(N_n)_n$ the
projective system of the universal images of $(K_n)_n;$

$\bullet$ there exists an integer $k\ge0$ such that the
projective system
$(L_n)_n:=(N_{n+k}\otimes\Lambda_n)_n$ is adic.

Moreover, $(K_n)_n$ is AR-isomorphic to $(L_n)_n.$ Now
for each $i,$ the projective system $\mathscr H^i(C)$ is
AR-adic (\ref{R2.2}i). Let $N^i=(N^i_n)_n$ be the
projective system of the universal images of $\mathscr
H^i(C),$ and choose an integer $k\ge0$ such that the
system $L^i=(L^i_n)_n=(N^i_{n+k}\otimes\Lambda_n)_n$ is
adic. Since $N^i_{n+k}\subset\mathscr H^i(C_{n+k})$ is
annihilated by $\lambda,$ we have $L^i_n=N^i_{n+k},$ and
the transition morphism gives an isomorphism
$$
\xymatrix@C=.6cm{
L^i_n\simeq L^i_n\otimes_{\Lambda_n}\Lambda_{n-1}
\ar[r]^-{\sim} & L^i_{n-1}.}
$$
This means the projective system $L^i$ is the constant
system $\pi^*L^i_0.$ By (\cite{LO2}, 2.2.2) we have
$R\pi_*\mathscr H^i(C)\simeq R\pi_*L^i,$ which is just $L^i_0$
by (\cite{LO2}, 2.2.3).

The spectral sequence
$$
R^j\pi_*\mathscr H^i(C)\Longrightarrow\mathscr
H^{i+j}((\widehat{M})_0)
$$
degenerates to isomorphisms $L^i_0\simeq\mathscr
H^i((\widehat{M})_0),$ so we only need to show that $L^i_0$
is trivialized by $(\mathscr S,\mathcal L).$ Using the
periodic $\Lambda_n$-flat resolution of $\Lambda_0$
$$
\xymatrix@C=.5cm{
\cdots \ar[r] & \Lambda_n \ar[r]^-{\lambda} & \Lambda_n
\ar[r]^-{\lambda^n} & \Lambda_n \ar[r]^-{\lambda} &
\Lambda_n \ar[r]^-{\epsilon} & \Lambda_0 \ar[r] & 0,}
$$
we see that $\Lambda_0\otimes^L_{\Lambda_n}\mathscr
H^j(M_n)$ is represented by the complex
$$
\xymatrix@C=.5cm{
\cdots \ar[r] & \mathscr H^j(M_n) \ar[r]^-{\lambda^n} &
\mathscr H^j(M_n) \ar[r]^-{\lambda} & \mathscr H^j(M_n)
\ar[r] & 0,}
$$
so $\s H^i(\Lambda_0\otimes^L_{\Lambda_n}\s 
H^j(M_n))$ are trivialized by $(\s S,\mathcal L),$
for all $i,j.$ Since $M$ is essentially bounded, we have the 
spectral sequence
$$
\s H^i(\Lambda_0\otimes^L_{\Lambda_n}\s 
H^j(M_n))\Longrightarrow\s H^{i+j}(C_n),
$$
from which we deduce (by (\ref{R3.4})) that the 
$\s H^i(C_n)$'s are trivialized by $(\s S,
\mathcal L).$ The universal image $N^i_n$ is
the image of $\s H^i(C_{n+r})\to\s H^i(C_n)$
for some $r\gg0,$ therefore the $N^i_n$'s (and hence the
$L^i_n$'s) are trivialized by $(\s S,\mathcal L).$

For the second claim, let $K\in D_c(\mathcal X,\Lambda).$
Since $K$ is isomorphic to the image of $\widehat{K}$
under the localization $\s D_c(\s A)\to
D_c(\mathcal X,\Lambda)$ (\cite{LO2}, 3.0.14), we see
that $K$ is stratifiable if $\widehat{K}$ is. Conversely,
if $K$ is stratifiable, which means that it is isomorphic
to the image of some $M\in\s D_c^{\text{stra}}
(\s A),$ then $\widehat{K}=\widehat{M}$ is also
stratifiable.
\end{proof}

\begin{anitem}\label{D-stra}
For $K\in D_c(\mathcal X,\Lambda),$ we say that $K$ is
$(\mathscr S,\mathcal L)$\textit{-stratifiable} if
$\widehat{K}$ is, and (\ref{L3.4.5}) implies that $K\in
D_{\mathscr S,\mathcal L}(\mathcal X,\Lambda)$ (cf.
(\ref{R3.4})) if and only if $K$ is $(\mathscr
S,\mathcal L)$-stratifiable.
\end{anitem}

\begin{corollary}\label{C3.5}
(i) If $\mathscr S$ is a stratification of $\mathcal X,$
then $K\in D_c(\mathcal X,\Lambda)$ is stratifiable if and
only if $K|_V$ is stratifiable for every $V\in\mathscr S.$

(ii) Let $K\in D_c(\mathcal X,\Lambda).$ Then $K$ is
stratifiable if and only if its Tate twist $K(1)$ is.

(iii) Let $P:X\to\mathcal X$ be a presentation, and let $K\in
D_c(\mathcal X,\Lambda).$ Then $K$ is stratifiable if and
only if $P^*K$ (resp. $P^!K$) is stratifiable.
\end{corollary}

\begin{proof}
(i) The ``only if" part follows from (\ref{L3.3}ii). For the
``if" part, we first prove the following.

\begin{sublemma}\label{pull-normalized}
For an $S$-algebraic stack $\mathcal X$ locally of finite 
type, let $N\overset{u}{\to}M\to C\to N[1]$ be an exact 
triangle in $\mathscr D_c(\mathscr A),$ where $N$ is a 
normalized complex and $C$ is almost 
AR-null. Then the morphism $u$ 
is isomorphic to the natural map $\widehat{M}\to M.$
\end{sublemma}

\begin{proof}
Consider the following diagram 
$$
\xymatrix@C=1cm{
\widehat{N} \ar[r]^-{\widehat{u}} \ar[d]_-{\simeq} & 
\widehat{M} \ar[r] \ar[d] & \widehat{C} \ar[r] \ar[d] & 
\widehat{N}[1] \ar[d] \\ 
N \ar[r]^-u & M \ar[r] & C \ar[r] & N[1].}
$$
Since $C$ is almost AR-null, we have $\widehat{C}=0$ by 
(\cite{LO2}, 2.2.2), and so $\widehat{u}$ is an isomorphism.
\end{proof}

Now let $f:\mathcal V\to\mathcal X$ be a morphism of 
$S$-algebraic stacks, and let $M\in\mathscr D_c(\mathscr 
A(\mathcal X)).$ We claim that $f^*\widehat{M}\simeq 
\widehat{f^*M}.$ Applying $f^*$ to the exact triangle 
$$
\xymatrix@C=.7cm{
\widehat{M} \ar[r] & M \ar[r] & C \ar[r] &}
$$
we get 
$$
\xymatrix@C=.7cm{
f^*\widehat{M} \ar[r] & f^*M \ar[r] & f^*C \ar[r] &.}
$$
By (\cite{LO1}, 4.3.2), 
$\widehat{M}_n=\text{hocolim}_N\tau_{\le N} 
\widehat{M}_n,$ and $-\otimes^L_{\Lambda_n}\Lambda_{n-1}$ 
and $f^*$ preserve homotopy colimit because they 
preserve infinite direct sums. Now that $\tau_{\le N}
\widehat{M}_n$ and $\Lambda_{n-1}$ are bounded above 
complexes, we have $f^*(\tau_{\le N}\widehat{M}_n 
\otimes^L_{\Lambda_n}\Lambda_{n-1})\simeq 
f^*\tau_{\le N}\widehat{M}_n\otimes^L_{\Lambda_n} 
\Lambda_{n-1}$ (cf. the proof of (\cite{LO1}, 4.5.3)). 
Hence applying $f^*$ to the isomorphism 
$$
\xymatrix@C=.7cm{
\widehat{M}_n\otimes^L_{\Lambda_n}\Lambda_{n-1} 
\ar[r] & \widehat{M}_{n-1}}
$$
we get an isomorphism 
$$
\xymatrix@C=.7cm{
f^*\widehat{M}_n\otimes^L_{\Lambda_n}\Lambda_{n-1} 
\ar[r] & f^*\widehat{M}_{n-1},}
$$
and by (\cite{LO2}, 3.0.10), $f^*\widehat{M}$ is 
normalized. Also it is clear 
that $f^*C$ is AR-null. By (\ref{pull-normalized}) we 
have $f^*\widehat{M}\simeq\widehat{f^*M}.$ 

Therefore, the ``if" part follows from (\ref{L3.3}iii) 
and (\ref{L3.4.5}), since 
$\widehat{K}|_V\simeq\widehat{(K|_V)}.$

(ii) This follows from (\ref{L3.3}vii), since
$\widehat{K}(1)=\widehat{K(1)}.$

(iii) For $P^*K,$ the ``only if" part follows from
(\ref{L3.3}ii), and the ``if" part follows from
(\ref{L3.3}iv) and (\ref{L3.4.5}), since
$P^*\widehat{K}=\widehat{(P^*K)}$ (\cite{LO2}, 2.2.1,
3.0.11).

Since $P$ is smooth of relative dimension $d,$ for some
function $d:\pi_0(X)\to\mathbb N,$ we have $P^!K\simeq
P^*K(d)[2d],$ so by (ii), $P^*K$ is stratifiable if and only
if $P^!K$ is.
\end{proof}

Before proving the main result of this section, we prove some
special cases.

\begin{blank}\label{pushforward-unbounded}
Let $f:X\to Y$ be a morphism of $S$-schemes. Then the 
$\Lambda_n$-dualizing complexes $K_{X,n}$ and $K_{Y,n}$ 
of $X$ and $Y$ respectively have finite quasi-injective 
dimensions, and are bounded by some integer independent of 
$n.$ Together with the base change theorem for $f_!,$ we 
see that there exists an integer $N>0$ depending only on 
$X,Y$ and $f,$ such that for any integers $a,b$ and $n$ 
with $n\ge0$ and any $M\in\s D_c^{[a,b]}(X,\Lambda_n),$ we 
have $f_*M\in\s D_c^{[a,b+N]}(Y,\Lambda_n).$ This implies that 
for each $n,$ the functor (defined using $K$-injective 
resolutions; cf. (\cite{Spa}, 6.7)) 
$$
f_*:\s D(X,\Lambda_n)\to\s D(Y,\Lambda_n)
$$
restricts to 
$$
f_*:\s D_c(X,\Lambda_n)\to\s D_c(Y,\Lambda_n).
$$
Moreover, for $M\in\s{D(A}(X))$ with constructible 
$\s H^j(M_n)$'s (for all $j$ and $n$) and for each 
$i\in\bb Z,$ there exist integers $a<b$ such that 
$$
R^if_*M\simeq R^if_*\tau_{[a,b]}M.
$$
In particular, if $M$ is a $\lambda$-complex on $X,$ then 
$R^if_*M$ is AR-adic for each $i,$ and hence $f_*M=
(f_*M_n)_n$ is a $\lambda$-complex on $Y.$ 

This enables us to define 
$$
f_*:D_c(X,\Lambda)\to D_c(Y,\Lambda)
$$
to be $K\mapsto Qf_*\widehat{K},$ where $Q:\s D_c
(\s A(Y))\to D_c(Y,\Lambda)$ is the localization 
functor. It agrees with the definition in (\cite{LO2}, 8) 
when restricted to $D_c^+(X,\Lambda),$ and for each $i\in
\bb Z$ and $K\in D_c(X,\Lambda),$ there exist integers 
$a<b$ such that $R^if_*K\simeq R^if_*\tau_{[a,b]}K.$
\end{blank}

\begin{lemma}\label{L3.8}
(i) If $f:X\to Y$ is a morphism of $S$-schemes, and $K\in
D_c(X,\Lambda)$ is trivialized by $(\{X\},\mathcal L)$ for
some $\mathcal L,$ then $f_*K$ is stratifiable.

(ii) Let $\mathcal X$ be an $S$-algebraic stack that has a 
connected presentation (i.e. there exists a presentation 
$P:X\to\mathcal X$ with $X$ a connected $S$-scheme). 
Let $K_{\mathcal X}$ and $K'_{\mathcal X}$ be two
$\Lambda$-dualizing complexes on 
$\mathcal X,$ and let $D$ and $D'$ be the two associated 
dualizing functors, respectively. Let $K\in D_c(\mathcal 
X,\Lambda).$ If $DK$ is trivialized by a pair $(\mathscr 
S,\mathcal L),$ where all strata in $\mathscr S$ are 
connected, then $D'K$ is trivialized by $(\mathscr 
S,\mathcal L')$ for some other $\mathcal L'.$
In particular, for stacks with connected presentation, 
the property 
of the Verdier dual of $K$ being stratifiable is 
independent of the choice of the dualizing complex.

(iii) Let $\mathcal X$ be an $S$-algebraic stack that has 
a connected presentation, and assume that the constant 
sheaf $\Lambda$ on $\mathcal X$ is a dualizing complex. 
If $K\in D_c(\mathcal X,\Lambda)$ is 
trivialized by a pair $(\{\mathcal X\},\mathcal L),$ then 
$D_{\mathcal X}K$ is trivialized by $(\{\mathcal 
X\},\mathcal L')$ for some $\mathcal L'.$
\end{lemma}

\begin{proof}
(i) Since $f_*K$ is the image of $f_*\widehat{K},$ it 
suffices to show that $f_*\widehat{K}$ is stratifiable. 
Since $f_*L$ is bounded for each $L\in\mathcal L,$ there 
exists a pair $(\mathscr S_Y, \mathcal L_Y)$ on $Y$ that 
trivializes $f_*L,$ for all $L\in\mathcal L.$ We claim that 
this pair trivializes $R^if_*\widehat{K}_n,$ for each $i$ 
and $n.$

Since $R^if_*\widehat{K}_n=R^if_*\tau_{[a,b]}\widehat{K}_n$ 
for some $a<b,$ and $\tau_{[a,b]}\widehat{K}_n$ is trivialized 
by $(\{X\},\mathcal L),$ we may assume $\widehat{K}_n$ is 
bounded. The claim then follows from the spectral sequence
$$
R^pf_*\mathscr H^q((\widehat{K})_n)\Longrightarrow R^{p+q}
f_*((\widehat{K})_n)
$$
and (\ref{R3.4}). 

(ii) Recall that the dualizing complex $K_{\mathcal X}$
(resp. $K'_{\mathcal X}$) is defined to be the image of a
normalized complex $K_{\mathcal X,\bullet}$ (resp.
$K'_{\mathcal X,\bullet}$), where each $K_{\mathcal X,n}$ 
(resp. $K'_{\mathcal X,n}$) is a $\Lambda_n$-dualizing 
complex. See (\cite{LO2}, 7.2.3, 7.2.4). 

Let $P:X\to\mathcal X$ be a presentation where $X$ is a 
connected scheme. Then we have 
$$
P^*R\mathscr Hom(K_{\mathcal X,n},K'_{\mathcal X,n})=
R\mathscr Hom(P^*K_{\mathcal X,n},P^*K'_{\mathcal X,n})
=R\mathscr Hom(P^!K_{\mathcal X,n},P^!K'_{\mathcal X,n}).
$$
Since $P^!K_{\mathcal X,n}$ and $P^!K'_{\mathcal X,n}$ are 
$\Lambda_n$-dualizing complexes on $X,$ by (\cite{SGA5}, 
Exp. I, 2) we see that $P^*R\mathscr Hom(K_{\mathcal 
X,n},K'_{\mathcal X,n})$ (and hence $R\mathscr Hom
(K_{\mathcal X,n},K'_{\mathcal X,n})$) is cohomologically 
concentrated in one degree, therefore it is quasi-isomorphic 
to this non-trivial cohomology sheaf, appropriately shifted. 
So let $R\mathscr Hom(K_{\mathcal X,n},K'_{\mathcal X,n})
\simeq L_n[r_n]$ for some sheaf $L_n$ and integer $r_n.$ 
Since $P^*L_n$ is invertible and hence lcc (cf. (\cite
{SGA5}, p.19)), the sheaf $L_n$ is lcc (\cite{Ols3}, 9.1).

For every stratum $\mathcal U\in\mathscr S,$ let $\mathcal 
L_0(\mathcal U)$ be the union of $\mathcal{L(U)}$ and the 
set of isomorphism classes of the JH components of the lcc 
sheaf $L_0|_{\mathcal U}.$ 
Since all strata in $\mathscr S$ are connected,
there exists a tensor closed hull of $(\mathscr S,\mathcal
L_0)$ of the form $(\mathscr S,\mathcal L'),$ i.e. they
have the same stratification $\mathscr S.$

By (\cite{LO2}, 4.0.8), the
system $(L_n[r_n])_n=R\mathscr Hom((K_{\mathcal X,n})_n,
(K'_{\mathcal X,n})_n)$ is normalized, so by (\ref
{pull-normalized}), $\widehat{D'K_{\mathcal X}}=(L_n[r_n])
_n,$ and by (\ref{L3.3}vi), it is trivialized by 
$(\mathscr S,\mathcal L').$ Since $DK$ is trivialized 
by $(\mathscr S,\mathcal L'),$ so also is $D'K,$ 
because $\widehat{D'K}\simeq\widehat{DK}\otimes^L\widehat
{D'K_{\mathcal X}}.$ 

(iii) The assumption implies in particular 
that $\mathcal X$ is connected, 
so by (ii), the question is independent of 
the choice of the dualizing complex. 
By definition, $\widehat{K}$ is trivialized by 
$(\{\mathcal X\},\mathcal L),$ so are truncations of 
$\widehat{K}.$ The essential image of $R\mathscr Hom
(\widehat{K},\Lambda_{\bullet})$ in 
$D_c(\mathcal X,\Lambda)$ is $DK,$ so by (\ref{D-stra}) 
it suffices to show 
that $R\mathscr Hom(\widehat{K},\Lambda_{\bullet}) 
\in\mathscr D_{\{\mathcal X\},\mathcal L'}
(\mathscr A)$ for some $\mathcal L'.$ 

Since $\mathcal X$ is quasi-compact, 
each $\Lambda_n$-dualizing complex 
is of finite quasi-injective dimension, so 
for each integer $i,$ there exist 
integers $a$ and $b$ such that 
$$
\mathscr H^iR\mathscr Hom(\widehat{K}_n,
\Lambda_n)=\mathscr H^i
R\mathscr Hom(\tau_{[a,b]}\widehat{K}_n,\Lambda_n).
$$
Using truncation triangles, we may further replace 
$\tau_{[a,b]}\widehat{K}_n$ by the cohomology sheaves 
$\mathscr H^j\widehat{K}_n,$ and hence by their JH 
components. Therefore, it suffices to find an 
$\mathcal L'$ that trivializes $\mathscr 
H^iR\mathscr Hom(L,\Lambda_0),$ for 
all $i\in\mathbb Z$ and $L\in\mathcal L.$ 
Note that $R\mathscr Hom(L,\Lambda_0)=\mathscr Hom(L,
\Lambda_0)=L^{\vee}$ is a simple $\Lambda_0$-sheaf, 
so one can take $\mathcal 
L'=\{L^{\vee}|L\in\mathcal L\}.$ 
\end{proof}

\begin{subremark}
For any $S$-algebraic stack $\mathcal X,$ the Verdier 
dual of a complex $K\in D_c(\mathcal X,\Lambda)$ being 
stratifiable or not is independent of the choice of 
the dualizing complex. Let $K_{\mathcal X}$ and $K'
_{\mathcal X}$ be two dualizing complexes on $\mathcal 
X,$ defining dualizing functors $D$ and $D',$ 
respectively. Let 
$P:X\to\mathcal X$ be a presentation, and let 
$K_X=P^!K_{\mathcal X}$ and $K'_X=P^!K'_{\mathcal X},$ 
defining dualizing functors $D_X$ and $D'_X$ on $X,$ 
respectively. Suppose $DK$ is stratifiable. To show 
$D'K$ is also stratifiable, by (\ref{C3.5}iii) it 
suffices to show $P^!D'K=D'_XP^*K$ is stratifiable. Since 
$D_XP^*K=P^!DK$ is stratifiable by assumption, we may 
assume $\mathcal X=X$ is a scheme. Since $X$ is 
noetherian, it has finitely many connected components, 
each of which is both open and closed. Then the result 
follows from (\ref{C3.5}i) and (\ref{L3.8}ii).
\end{subremark}

Next we prove the main result of this section.

\begin{theorem}\label{T3.10}
Let $f:\mathcal X\to\mathcal Y$ be a morphism of
$S$-algebraic stacks. Then the operations
$f_*,f_!,f^*,f^!,D_{\mathcal X},-\otimes^L-$ and $R\mathscr
Hom(-,-)$ all preserve stratifiability, namely, they induce
functors
\begin{gather*}
f_*:D_c^{+,\emph{stra}}(\mathcal X,\Lambda)
\longrightarrow D_c^{+,\emph{stra}}(\mathcal
Y,\Lambda), \qquad
f_!:D_c^{-,\emph{stra}}(\mathcal X,\Lambda)
\longrightarrow D_c^{-,\emph{stra}}(\mathcal
Y,\Lambda), \\
f^*:D_c^{\emph{stra}}(\mathcal Y,\Lambda)
\longrightarrow D_c^{\emph{stra}}(\mathcal
X,\Lambda), \qquad
f^!:D_c^{\emph{stra}}(\mathcal Y,\Lambda)
\longrightarrow D_c^{\emph{stra}}(\mathcal
X,\Lambda), \\
D:D_c^{\emph{stra}}(\mathcal X,\Lambda)
\longrightarrow D_c^{\emph{stra}}(\mathcal
X,\Lambda)^{\emph{op}}, \\
\otimes^L:D_c^{-,\emph{stra}}(\mathcal
X,\Lambda)\times D_c^{-,\emph{stra}}(\mathcal X,
\Lambda)\longrightarrow D_c^{-,\emph{stra}}(\mathcal
X,\Lambda)\quad\emph{and} \\
R\mathscr Hom(-,-):D_c^{-,\emph{stra}}(\mathcal
X,\Lambda)^{\emph{op}}\times D_c^{+,\emph{stra}}
(\mathcal X,\Lambda)\longrightarrow D_c^{+,\emph{stra}}
(\mathcal X,\Lambda).
\end{gather*}
\end{theorem}

\begin{proof}
We may assume all stacks involved are reduced. 

We consider the Verdier dual functor $D$ first. Let
$P:X\to\mathcal X$ be a presentation. Since $P^*D=DP^!,$ by
(\ref{C3.5}iii) we can assume $\mathcal X=X$ is a scheme. Let
$K$ be a complex on $X$ trivialized by $(\s S,\mathcal
L).$ Refining if necessary, we may assume all strata in 
$\s S$ are connected and regular. Let $j:U\to X$ be
the immersion of an open stratum in $\s S$ 
with complement $i:Z\to X.$ Shrinking $U$ if necessary, 
we may assume there is a dimension function on $U$ 
(\cite{Rio}, D\'efinition 2.1), hence by a result of 
Gabber (loc. cit., Th\'eor\`eme 0.2), the constant sheaf 
$\Lambda$ on $U$ is a dualizing complex. 
Consider the exact triangle
$$
\xymatrix@C=.5cm{
i_*D_Z(K|_Z) \ar[r] & D_XK \ar[r] & j_*D_U(K|_U) \ar[r] &.}
$$
By (\ref{L3.8}iii) we see that $D_U(K|_U)$ is trivialized by 
$(\{U\},\mathcal L')$ for some $\mathcal L',$ so 
$j_*D_U(K|_U)$ is stratifiable (\ref{L3.8}i). By noetherian 
induction we may assume $D_Z(K|_Z)$ is stratifiable, and it 
is clear that $i_*$ preserves stratifiability. Therefore by 
(\ref{L3.3}i), $D_XK$ is stratifiable.

The case of $f^*$ (and hence $f^!$) is proved in
(\ref{L3.3}ii).

Next we prove the case of $\otimes^L.$ For $i=1,2,$ let
$K_i\in D_c^-(\mathcal X,\Lambda),$ trivialized by $(\s 
S_i,\mathcal L_i).$ Let $(\s S,\mathcal L)$ be a common
tensor closed refinement (by (\ref{L3.2.5})) of $(\s 
S_i,\mathcal L_i),\ i=1,2.$ The total tensor product
$K_1\otimes^LK_2$ is defined to be the image in $D_c(\mathcal
X,\Lambda)$ of $\widehat{K}_1\otimes^L_{\Lambda_{\bullet}}
\widehat{K}_2,$ which by (\cite{LO2}, 3.0.10) is normalized, 
so it suffices to show (by (\ref{L3.3}vi)) that
$$
\widehat{K}_{1,0}\otimes^L_{\Lambda_0}\widehat{K}_{2,0}=
\widehat{K}_{1,0}\otimes_{\Lambda_0}\widehat{K}_{2,0}
$$ 
is trivialized by $(\s S,\mathcal L).$ This follows from 
$$
\s H^r(\widehat{K}_{1,0}\otimes_{\Lambda_0}\widehat{K}
_{2,0})=\bigoplus_{i+j=r}\s H^i(\widehat{K}_{1,0})
\otimes_{\Lambda_0}\s H^j(\widehat{K}_{2,0})
$$
and the assumption that $(\s S,\mathcal L)$ is tensor closed.

The case of $R\s Hom(K_1,K_2)=D(K_1\otimes^LDK_2)$
follows.

Finally we prove the case of $f_*$ and $f_!.$ Let $f:\mathcal
X\to\mathcal Y$ be a morphism of $S$-algebraic stacks, 
and let $K\in D^-_{\s S,\mathcal L}(\mathcal X,\Lambda)$ 
for some pair $(\s S,\mathcal L).$ We want to show
$f_!K$ is stratifiable. Let $j:\mathcal U\to\mathcal X$ be
the immersion of an open stratum in $\s S,$ with
complement $i:\mathcal Z\to\mathcal X.$ From the exact
triangle
$$
\xymatrix@C=.5cm{
(fj)_!j^*K \ar[r] & f_!K \ar[r] & (fi)_!i^*K \ar[r] &}
$$
we see that it suffices to prove the claim for $fj$ and $fi.$
By noetherian induction we can replace $\mathcal X$ by
$\mathcal U.$ By (\ref{C3.5}iii) and smooth base change
(\cite{LO2}, 12.1), we can replace $\mathcal Y$ by a
presentation $Y,$ and by (\ref{C3.5}i) and (\cite{LO2}, 12.3)
we can shrink $Y$ to an open subscheme. After these reductions, 
we assume that $\mathcal Y=Y$ is a regular irreducible 
affine $S$-scheme that has a dimension function on it, that $K$ 
is trivialized by $(\{\mathcal X\},\mathcal L),$ and that the
relative inertia stack $\mathcal I_{f}:=\mathcal
X\times_{\Delta,\mathcal X\times_Y\mathcal X,\Delta}\mathcal
X$ is flat and has components over $\mathcal X$ (\cite{Beh2},
5.1.14). Therefore by (\cite{Beh2}, 5.1.13), $f$ factors as
$\mathcal X\overset{g}{\to}\mathcal Z\overset{h}{\to}Y,$
where $g$ is gerbe-like and $h$ is representable (cf.
(\cite{Beh2}, 5.1.3-5.1.6) for relevant notions). So we
reduce to two cases: $f$ is representable, or $f$ is
gerbe-like.

\textbf{Case when $f$ is representable}. By shrinking the
$S$-algebraic space $\mathcal X$ we can assume $\mathcal X=X$
is a regular connected scheme that has a dimension function, 
so that the constant sheaf $\Lambda$ on $X$ is 
a dualizing complex. By (\ref{L3.8}iii) we see that 
$DK$ is trivialized by some $(\{X\},\mathcal L'),$ and
by (\ref{L3.8}i), $f_*DK$ is stratifiable. Therefore
$f_!K=Df_*DK$ is also stratifiable.

\textbf{Case when $f$ is gerbe-like}. In this case $f$ is
smooth (\cite{Beh2}, 5.1.5), hence \'etale locally on $Y$ it
has a section. Replacing $Y$ by an \'etale cover, we may
assume that $f$ is a neutral gerbe, so $f:B(G/Y)\to Y$ is the
structural map, for some flat group space $G$ of finite type
over $Y$ (\cite{LMB}, 3.21). By (\cite{Beh2}, 5.1.1) and
(\ref{C3.5}i) we may assume $G$ is a $Y$-group scheme. Next
we reduce to the case when $G$ is smooth over $Y.$

By assumption $Y$ is integral. 
Let $k(Y)$ be the function field of $Y$ and $\overline{k(Y)}$
an algebraic closure. Then $G_{\overline{k(Y)},\text{red}}$
is smooth over $\overline{k(Y)},$ so there exists a finite
extension $L$ over $k(Y)$ such that $G_{L,\text{red}}$ is
smooth over $L.$ Let $Y'$ be the normalization of $Y$ in $L,$
which is a scheme of finite type over $S,$ and the natural
map $Y'\to Y$ is finite
surjective. It factors through $Y'\to Z\to Y,$ where $Z$ is
the normalization of $Y$ in the separable closure of $k(Y)$ in
$L=k(Y').$ So $Z\to Y$ is generically \'etale, and $Y'\to Z$
is purely inseparable, hence a universal homeomorphism, so
$Y'$ and $Z$ have equivalent \'etale sites. Replacing $Y'$ by
$Z$ and shrinking $Y$ we can assume $Y'\to Y$ is finite
\'etale. Replacing $Y$ by $Y'$ (by (\ref{C3.5}ii)) we assume
$G_{\text{red}}$ over $Y$ has smooth generic fiber, and by
shrinking $Y$ we assume $G_{\text{red}}$ is smooth over $Y.$

$G_{\text{red}}$ is a subgroup scheme of $G$ (\cite{SGA3},
VI$_{\text{A}},$ 0.2); let $h:G_{\text{red}}\hookrightarrow G$
be the closed immersion. Then $Bh:B(G_{\text{red}}/Y)\to
B(G/Y)$ is faithful and hence representable. It is also
radicial: consider the diagram where the square is
2-Cartesian
$$
\xymatrix@C=.8cm{
Y \ar@{^{(}->} [r]^-i & G/G_{\text{red}} \ar[r]^-g \ar[d] & Y
\ar[d]^-P \\
& B(G_{\text{red}}/Y) \ar[r]_-{Bh} & B(G/Y).}
$$
The map $i$ is a nilpotent closed embedding, so $g$ is
radicial. Since $P$ is faithfully flat, $Bh$ is also
radicial. This shows that
$$
(Bh)^*:D_c^-(B(G/Y),\Lambda)\to D_c^-(B(G_{\text{red}}/Y),
\Lambda)
$$
is an equivalence of categories. Replacing $G$ by $G_{\text
{red}}$ we assume $G$ is smooth over $Y,$ and hence $P:Y\to
B(G/Y)$ is a connected presentation.

Let $d$ be the relative dimension of $G$ over $Y.$ By 
assumption, the constant sheaf $\Lambda$ on $Y$ is a dualizing 
complex, and so $f^!\Lambda=\Lambda\langle-d\rangle$ (and 
hence the constant sheaf $\Lambda$ on $\mathcal X$) is a 
dualizing complex on $\mathcal X.$ By (\ref{L3.8}iii), 
we see that $DK$ is trivialized by a pair of the form 
$(\{\mathcal X\},\mathcal L').$ 
To show $f_!K$ is stratifiable is equivalent to
showing that $Df_!K=f_*DK$ is stratifiable. So replacing 
$K$ by $DK,$ it suffices to show that $f_*K$ is stratifiable, 
where $K\in D^+_{\{\mathcal X\},\mathcal L}(\mathcal X, 
\Lambda)$ for some $\mathcal L.$

Consider the strictly simplicial smooth hypercover associated 
to the presentation $P:Y\to B(G/Y),$ and let $f_i:G^i\to Y$ 
be the structural map. As in the proof 
of (\ref{L3.8}i), it suffices to show the existence of a pair
$(\mathscr S_Y,\mathcal L_Y)$ on $Y$ that trivializes
$R^nf_*L,$ for all $L\in\mathcal L$ and $n\in\mathbb Z.$
From the spectral sequence (\cite{LO2}, 10.0.9)
$$
E_1^{ij}=R^jf_{i*}f_i^*P^*L\Longrightarrow R^{i+j}f_*L,
$$
we see that it suffices for the pair $(\mathscr S_Y,\mathcal
L_Y)$ to trivialize all the $E_1^{ij}$-terms. Assume $i\ge1.$
If we regard the map $f_i:G^i\to Y$ as the product map
$$
\prod_if_1:\prod_iG\to\prod_iY,
$$
where the products are fiber products over $Y,$ then we can
write $f_i^*P^*L$ as
$$
f_1^*P^*L\boxtimes_{\Lambda_0}\Lambda_0\boxtimes
_{\Lambda_0}\cdots\boxtimes_{\Lambda_0}\Lambda_0.
$$
Note that the scheme $Y$ satisfies the condition (LO). 
By K\"unneth formula (\cite{LO2}, 11.0.14) we have 
$$
f_{i*}f_i^*P^*L=f_{1*}f_1^*P^*L\otimes_{\Lambda_0}
f_{1*}\Lambda_0\otimes_{\Lambda_0}\cdots\otimes_{\Lambda_0}
f_{1*}\Lambda_0.
$$
Since $f_{1*}f_1^*P^*L$ and $f_{1*}\Lambda_0$ are bounded
complexes (by a theorem of Gabber \cite{Gab}), there 
exists a tensor closed pair $(\s 
S_Y,\mathcal L_Y)$ that trivializes them, for all
$L\in\mathcal L.$ The proof is finished.
\end{proof}

Consequently, the theorem also holds for $\overline
{\bb Q}_{\ell}$-coefficients.

Finally we give a lemma which will be used in the next
section. This will play the same role as (\cite{Beh2},
6.3.16).

\begin{lemma}\label{L3.11}
Let $X$ be a connected variety over an algebraically closed
field $k$ of characteristic not equal to $\ell,$ and let
$\mathcal L$ be a finite set of isomorphism classes of simple
lcc $\Lambda_0$-sheaves on $X.$ Then there exists an integer
$d$ (depending only on $\mathcal L$) such that, for every
lisse $\Lambda$-adic sheaf $\s F$ on $X$ trivialized by
$\mathcal L,$ and for every integer $i,$ we have
$$
\dim_EH^i_c(X,\s F\otimes_{\Lambda}E)\le
d\cdot\emph{rank}_E(\s F\otimes_{\Lambda}E),
$$
where $E$ is the fraction field of $\Lambda.$
\end{lemma}

\begin{proof}
Since $\mathcal L$ is finite and $0\le i\le2\dim X,$ there
exists an integer $d>0$ such that $\dim_{\Lambda_0}H^i_c(X,
L)\le d\cdot\text{rank}_{\Lambda_0}L,$ for every $i$ and
every $L\in\mathcal L.$ For a short exact sequence of lcc
$\Lambda_0$-sheaves
$$
\xymatrix@C=.5cm{
0 \ar[r] & \s G' \ar[r] & \s G \ar[r] & 
\s G'' \ar[r] & 0}
$$
on $X,$ the cohomological sequence
$$
\xymatrix@C=.5cm{
\cdots \ar[r] & H^i_c(X,\s G') \ar[r] &
H^i_c(X,\s G) \ar[r] & H^i_c(X,\s G'') \ar[r] & \cdots}
$$
implies that $\dim_{\Lambda_0}H^i_c(X,\s G)\le\dim
_{\Lambda_0}H^i_c(X,\s G')+\dim_{\Lambda_0}H^i_c
(X,\s G'').$ So it is clear that if $\s G$ is
trivialized by $\mathcal L,$ then
$\dim_{\Lambda_0}H^i_c(X,\s G)\le d\cdot
\text{rank}_{\Lambda_0}\s G,$ for every $i.$

Since we only consider $\s F\otimes_{\Lambda}E,$ we may
assume $\s F=(\s F_n)_n$ is flat, of some
constant rank over $\Lambda$ (since $X$ is connected), and
this $\Lambda$-rank is equal to
$$
\text{rank}_{\Lambda_0}\s F_0=\text{rank}_E(\s 
F\otimes_{\Lambda}E).
$$
$H^i_c(X,\s F)$ is a finitely generated
$\Lambda$-module (\cite{SGA5}, VI, 2.2.2), so by Nakayama's
lemma, the minimal number of generators is at most
$\dim_{\Lambda_0}(\Lambda_0\otimes_{\Lambda}H^i_c(X,\s 
F)).$ Similar to ordinary cohomology groups (\cite{Mil1},
19.2), we have an injection
$$
\Lambda_0\otimes_{\Lambda}H^i_c(X,\s F)\hookrightarrow
H^i_c(X,\s F_0)
$$
of $\Lambda_0$-vector spaces. Therefore, $\dim_EH^i_c(X,
\s F\otimes_{\Lambda}E)$ is less than or equal to the
minimal number of generators of $H^i_c(X,\s F)$ over
$\Lambda,$ which is at most
$$
\dim_{\Lambda_0}(\Lambda_0\otimes_{\Lambda}H^i_c(X,\s 
F))\le\dim_{\Lambda_0}H^i_c(X,\s F_0)\le d\cdot\text
{rank}_{\Lambda_0}\s F_0=d\cdot\text{rank}_E(\s 
F\otimes_{\Lambda}E).
$$
\end{proof}

\section{Convergent complexes and finiteness}

We return to $\bb F_q$-algebraic stacks $\s 
X_0,\s Y_0,\cdots$ of finite type. A complex 
$K_0\in W(\s X_0,\overline{\bb Q}_{\ell})$ is 
said to be \textit{stratifiable} if $K$ on 
$\s X$ is stratifiable, and we denote by 
$W^{\text{stra}}(\s X_0,\overline{\bb 
Q}_{\ell})$ the full subcategory of such complexes. Note 
that if $K_0$ is a lisse-\'etale complex, and it is 
stratifiable on $\s X_0,$ then it is geometrically 
stratifiable (i.e. $K$ on $\s X$ is stratifiable). 
In turns out that in order for the trace formula to hold, 
it suffices to make this weaker assumption of geometric 
stratifiability. So we will only discuss stratifiable 
Weil complexes. Again, by a sheaf we mean a Weil 
$\overline{\bb Q}_{\ell}$-sheaf. 

\begin{definition}\label{D4.1}
(i) Let $K\in D_c(\overline{\bb Q}_{\ell})$ and
$\varphi:K\to K$ an endomorphism. The pair $(K,\varphi)$ is
said to be an \emph{$\iota$-convergent complex} (or just a
convergent complex, since we fixed $\iota$) if the complex
series in two directions
$$
\sum_{n\in\bb Z}\quad\sum_{H^n(K),H^n(\varphi)}|\alpha|^s
$$
is convergent, for every real number $s>0.$ In this case let
$\emph{Tr}(\varphi,K)$ be the absolutely convergent complex
series
$$
\sum_n(-1)^n\iota\emph{Tr}(H^n(\varphi),H^n(K))
$$
or its limit.

(ii) Let $K_0\in W^-(\s X_0,\overline{\bb 
Q}_{\ell}).$ We call $K_0$ an \emph{$\iota$-convergent 
complex of sheaves} (or just a convergent complex of 
sheaves), if for every integer $v\ge1$ and every point 
$x\in\s X_0(\bb F_{q^v}),$ the pair 
$(K_{\overline{x}},F_x)$ is a convergent complex. In 
particular, all bounded complexes are convergent. 

(iii) Let $K_0\in W^-(\s X_0,\overline{\bb 
Q}_{\ell})$ be a convergent complex of sheaves. Define
$$
c_v(\s X_0,K_0)=\sum_{x\in[\s X_0(\bb
F_{q^v})]}\frac{1}{\#\emph{Aut}_x(\bb
F_{q^v})}\emph{Tr}(F_x,K_{\overline{x}})\in\bb C,
$$
and define the \emph{$L$-series} of $K_0$ to be the formal
power series
$$
L(\s X_0,K_0,t)=\exp\Big(\sum_{v\ge1}c_v(\s 
X_0,K_0)\frac{t^v}{v}\Big)\in\bb C[[t]].
$$
\end{definition}

The zeta function $Z(\s X_0,t)$ in (\ref{D1.2}) is a
special case: $Z(\s X_0,t)=L(\s
X_0,\overline{\bb Q}_{\ell},t).$ It has rational
coefficients.

\begin{notation}\label{}
We sometimes write $c_v(K_0)$ for $c_v(\s X_0,K_0),$ if
it is clear that $K_0$ is on $\s X_0.$ We also write
$c_v(\s X_0)$ for $c_v(\s X_0,\overline{\bb Q}_{\ell}).$
\end{notation}

\begin{subremark}\label{R4.2}
(i) Behrend defined convergent complexes with respect to
arithmetic Frobenius elements (\cite{Beh2}, 6.2.3), and our
definition is for geometric Frobenius, and it is essentially
the same as Behrend's definition, except we work with 
$\iota$-mixed Weil complexes (which means \textit{all} 
Weil complexes, by (\ref{Laff})) for an arbitrary 
isomorphism $\iota:\overline{\bb Q}_{\ell}\to\bb C,$ 
while \cite{Beh2} works 
with pure or mixed lisse-\'etale sheaves with integer 
weights. Also our definition is a little different from that
in \cite{Ols1}; the condition there is weaker.

(ii) Recall that $\text{Aut}_x$ is defined to be the fiber
over $x$ of the inertia stack $\s I_0\to\s 
X_0.$ It is a group scheme of finite type (\cite{LMB}, 4.2)
over $k(x),$ so $\text{Aut}_x(k(x))$ is a finite group.

(iii) If we have the following commutative diagram
$$
\xymatrix@C=.5cm{
\text{Spec }\bb F_{q^{vd}} \ar[r] \ar[rd]_{x'} &
\text{Spec }\bb F_{q^v} \ar[d]^x \\ & \s X_0,}
$$
then $(K_{\overline{x}},F_x)$ is convergent if and only if
$(K_{\overline{x'}},F_{x'})$ is convergent, because
$F_{x'}=F_x^d$ and $s\mapsto sd:\bb R^{>0}\to\bb
R^{>0}$ is a bijection. In particular, for a lisse-\'etale 
complex of sheaves, the property of being a convergent 
complex is independent of $q$ and the structural morphism 
$\s X_0\to\text{Spec }\bb F_q.$ Also note that, 
for every integer $v\ge1,$ a complex $K_0$ on $\s 
X_0$ is convergent if and only if $K_0\otimes\bb 
F_{q^v}$ on $\s X_0\otimes\bb F_{q^v}$ is convergent.
\end{subremark}

We restate the main theorem in \cite{Beh2} using compactly
supported cohomology as follows. It generalizes
(\ref{T1.1}). We will prove it in this section and the
next.

\begin{theorem}\label{T4.3}
Let $f:\s X_0\to\s Y_0$ be a morphism of 
$\bb F_q$-algebraic stacks, and let $K_0\in 
W^{-,\emph{stra}}_m(\s X_0,\overline{\bb 
Q}_{\ell})$ be a convergent complex of sheaves. Then 

(i) \emph{(Finiteness)} $f_!K_0$ is a convergent complex of
sheaves on $\s Y_0,$ and 

(ii) \emph{(Trace formula)} $c_v(\s X_0,K_0)=c_v(\s 
Y_0,f_!K_0)$ for every integer $v\ge1.$
\end{theorem}

First we give a few lemmas.

\begin{lemma}\label{L4.5}
Let
$$
\xymatrix@C=.7cm @R=.5cm{
K' \ar[r] \ar[d]_{\varphi'} & K \ar[r] \ar[d]_{\varphi} &
K'' \ar[r] \ar[d]^{\varphi''} & K'[1] \ar[d]^{\varphi'[1]} \\
K' \ar[r] & K \ar[r] & K'' \ar[r] & K'[1]}
$$
be an endomorphism of an exact triangle $K'\to K\to K''\to
K'[1]$ in $D^-_c(\overline{\bb Q}_{\ell}).$ If any two of
the three pairs $(K',\varphi'), (K'',\varphi'')$ and
$(K,\varphi)$ are convergent, then so is the third, and
$$
\emph{Tr}(\varphi,K)=\emph{Tr}(\varphi',K')+\emph{Tr}
(\varphi'',K'').
$$
\end{lemma}

\begin{proof}
By the rotation axiom we can assume $(K',\varphi')$ and
$(K'',\varphi'')$ are convergent. We have the exact sequence
$$
\xymatrix@C=.5cm{
\cdots \ar[r] & H^n(K') \ar[r] & H^n(K) \ar[r] & H^n(K'') \ar[r]
& H^{n+1}(K') \ar[r] & \cdots.}
$$
Since $H^n(K)$ is an extension of a sub-object of $H^n(K'')$ 
by a quotient object of $H^n(K'),$ we have 
$$
\sum_{H^n(K),\varphi}|\alpha|^s\le\sum_{H^n(K'),\varphi'}
|\alpha|^s+\sum_{H^n(K''),\varphi''}|\alpha|^s
$$
for every real $s>0,$ so $(K,\varphi)$ is convergent.

Since the series $\sum_{n\in\bb Z}(-1)^n\sum_{H^n(K),
\varphi}\iota\alpha$ converges absolutely, we can change
the order of summation, and the second assertion follows 
if we split the long exact sequence above into short exact 
sequences. 
\end{proof}

\begin{corollary}\label{C4.6}
If $K_0'\to K_0\to K_0''\to K_0'[1]$ is an exact triangle in
$W^-(\s X_0,\overline{\bb Q}_{\ell}),$ and two 
of the three complexes $K_0',K_0''$ and $K_0$ are 
convergent complexes, then so is the third, and
$c_v(K_0)=c_v(K_0')+c_v(K_0'').$
\end{corollary}

\begin{proof}
For every $x\in\s X_0(\bb F_{q^v}),$ we have an
exact triangle
$$
\xymatrix@C=.5cm{
K'_{\overline{x}} \ar[r] & K_{\overline{x}} \ar[r] & K''_
{\overline{x}} \ar[r] &}
$$
in $D^-_c(\overline{\bb Q}_{\ell}),$ equivariant under the
action of $F_x.$ Then apply (\ref{L4.5}).
\end{proof}

\begin{lemma}\label{L4.9}
(\ref{T4.3}) holds for $f:\emph{Spec }\bb 
F_{q^d}\to\emph{Spec }\bb F_q.$
\end{lemma}

\begin{proof}
We have an equivalence of triangulated categories
$$
\xymatrix@C=.7cm{
W^-(\text{Spec }\bb F_q,\overline{\bb 
Q}_{\ell}) \ar[r]^-{\sim} & D^-_c(\text{Rep}_{ 
\overline{\bb Q}_{\ell}}(G)),}
$$
where $G$ is the Weil group $F^{\bb Z}$ of $\bb 
F_q.$ Let $H$ be the subgroup $F^{d\bb Z},$ the Weil 
group of $\bb F_{q^d}.$ Since $f:\text{Spec }\bb 
F_{q^d}\to\text{Spec }\bb F_q$ is finite, we have 
$f_!=f_*,$ and it is the induced-module functor
$$
\xymatrix@C=.5cm{
\text{Hom}_{\overline{\bb Q}_{\ell}[H]}\big(\overline{
\bb Q}_{\ell}[G],-\big):D_c^-\big(\text{Rep}_{\overline
{\bb Q}_{\ell}}(H)\big) \ar[r] &
D_c^-\big(\text{Rep}_{\overline{\bb 
Q}_{\ell}}(G)\big),}
$$
which is isomorphic to the coinduced-module functor
$\overline{\bb Q}_{\ell}[G]\otimes_{\overline{\bb
Q}_{\ell}[H]}-.$ In particular, $f_!$ is exact on the level
of sheaves.

Let $A$ be a $\overline{\bb Q}_{\ell}$-representation
of $H,$ and $B=\overline{\bb Q}_{\ell}[G]\otimes_
{\overline{\bb Q}_{\ell}[H]}A.$ Let $x_1,\cdots,x_m$ be
an ordered basis for $A$ with respect to which $F^d$ is an
upper triangular matrix
$$
\begin{bmatrix}\alpha_1 & * & * \\ & \ddots & * \\ &&
\alpha_m \end{bmatrix}
$$
with eigenvalues $\alpha_1,\cdots,\alpha_m.$ Then $B$ has 
an
ordered basis
\begin{gather*}
1\otimes x_1,\ F\otimes x_1,\ \cdots,\ F^{d-1}\otimes x_1,
\\
1\otimes x_2,\ F\otimes x_2,\ \cdots,\ F^{d-1}\otimes x_2,
\\
\cdots\ \cdots \\
1\otimes x_m,\ F\otimes x_m,\ \cdots,\ F^{d-1}\otimes x_m,
\end{gather*}
with respect to which $F$ is the matrix
$
\begin{bmatrix}
M_1 & * & * \\
& \ddots & * \\
&& M_m\end{bmatrix},
$
where $M_i=\begin{bmatrix}0 & \cdots & 0 & \alpha_i \\
1 & \ddots & 0 & 0 \\
& \ddots & \ddots & \vdots \\
&& 1 & 0\end{bmatrix}.$ The characteristic polynomial of 
$F$
on $B$ is $\prod_{i=1}^m(t^d-\alpha_i).$

Let $K_0$ be a complex of sheaves on $\text{Spec }\bb
F_{q^d}.$ The eigenvalues of $F$ on $\s
H^n(f_!K)=f_!\s H^n(K)$ are all the $d$-th roots of 
the eigenvalues of $F^d$ on $\s H^n(K),$ so for every 
$s>0$ we have
$$
\sum_n\quad\sum_{\s H^n(f_!K),F}|\alpha|^s=d\sum_n\quad
\sum_{\s H^n(K),F^d}|\alpha|^{s/d}.
$$
This shows that $f_!K_0$ is a convergent complex on
$\text{Spec }\bb F_q$ if and only if $K_0$ is a
convergent complex on $\text{Spec }\bb F_{q^d}.$

Next we prove
$$
c_v(\text{Spec }\bb F_{q^d},K_0)=c_v(\text{Spec 
}\bb F_q,f_!K_0)
$$
for every 
$v\ge1.$ Since $H^n(f_!K)=f_!H^n(K),$ and both sides are
absolutely convergent series so that one can change the 
order 
of summation without changing the limit, it suffices to 
prove 
it when $K=A$ is a single representation concentrated in
degree 0. Let us review this classical calculation. Use the
notation as above. For each $i,$ fix a $d$-th root
$\alpha_i^{1/d}$ of $\alpha_i,$ and let $\zeta_d$ be a
primitive $d$-th root of unity. Then the eigenvalues of $F$
on $B$ are $\zeta_d^k\alpha_i^{1/d},$ for $i=1,\cdots,m$ 
and $k=0,\cdots,d-1.$

If $d\nmid v,$ then $\text{Hom}_{\bb F_q}(\bb
F_{q^d},\bb F_{q^v})=\emptyset,$ so $c_v(\text{Spec
}\bb F_{q^d},A)=0.$ On the other hand,
$$
c_v(\text{Spec }\bb 
F_q,f_!A)=\text{Tr}(F^v,B)=\sum_{i,k}\zeta_d
^{vk}\alpha_i^{v/d}=\sum_i\alpha_i^{v/d}\sum_{k=0}^{d-1}
\zeta_d^{vk}=0.
$$
If $d|v,$ then $\text{Hom}_{\bb F_q}(\bb F_{q^d},
\bb F_{q^v})=\text{Hom}_{\bb F_q}(\bb F_{q^d},
\bb F_{q^d})=\bb Z/d\bb Z.$ So
$$
c_v(\bb 
F_{q^d},A)=d\text{Tr}(F^v,A)=d\sum_i\alpha_i^{v/d}.
$$
On the other hand,
$$
c_v(\bb F_q,B)=\text{Tr}(F^v,B)=\sum_{i,k}\zeta_d^{vk}
\alpha_i^{v/d}=\sum_{i,k}\alpha_i^{v/d}=d\sum_i\alpha_i^{v/d}.
$$
\end{proof}

Next, we consider $BG_0,$ for a finite group scheme $G_0$ 
over $\bb F_q.$

\begin{lemma}\label{L4.8}
Let $G_0$ be a finite $\bb F_q$-group scheme, and let 
$\s F_0$ be a sheaf on $BG_0.$ Then 
$H^r_c(BG,\s F)=0$ for all $r\ne0,$ and 
$H^0_c(BG,\s F)\simeq H^0(BG,\s F)$ has 
dimension at most $\emph{rank}(\s F).$ Moreover, 
the set of $\iota$-weights of $H^0_c(BG,\s F)$ is a subset 
of the $\iota$-weights of $\s F_0.$
\end{lemma}

\begin{proof}
By (\cite{Ols1}, 7.12-7.14) we can replace $G_0$ by its
maximal reduced closed subscheme, and assume $G_0$ is
reduced, hence \'etale. Then $G_0$ is the same as a finite
group $G(\mathbb F)$ with a continuous action of
$\text{Gal}(\mathbb F_q)$ (\cite{Mil2}, 7.8). We will also
write $G$ for the group $G(\bb F),$ if there is no
confusion. Since $\text{Spec }\bb F\to BG$ is 
surjective, we see that there is no non-trivial 
stratification on $BG.$ In particular, all sheaves on $BG$ 
are lisse, and they are just $\overline{\bb 
Q}_{\ell}$-representations of $G.$

$BG$ is quasi-finite and proper over $\bb F,$ with 
finite diagonal, so by (\cite{Ols1}, 5.8), 
$H^r_c(BG,\s F)=0$ for all $r\ne0.$ From 
(\cite{Ols1}, 5.1), we see that if $\s F$ is a sheaf 
on $BG$ corresponding to the representation $V$ of $G,$ 
then $H^0_c(BG,\s F)=V_G$ and $H^0(BG,\s 
F)=V^G,$ and there is a natural isomorphism
$$
v\mapsto\sum_{g\in G}gv: V_G\longrightarrow V^G.
$$
Therefore 
$$
h^0_c(BG,\s F)=\dim V_G\le\dim V=\text{rank}(\s 
F),
$$
and the weights of $V_G$ form a subset of the weights 
of $V$ (counted with multiplicities).
\end{proof}

\begin{blank}\label{alg-gp}
(i) If $k$ is a field, by a \textit{$k$-algebraic group} 
$G$ we mean a smooth $k$-group scheme of finite type. If 
$G$ is connected, then it is geometrically connected 
(\cite{SGA3}, VI$_{\text{A}},$ 2.1.1). 

(ii) For a connected $k$-algebraic group $G,$ let 
$a:BG\to\text{Spec }k$ be the structural map. Then 
$$
a^*:\Lambda\text{-Sh(Spec }k)\longrightarrow 
\Lambda\text{-Sh}(BG)
$$
is an equivalence of categories. This is because 

$\bullet BG$ has no non-trivial stratifications (it is 
covered by $\text{Spec }k$ which has no non-trivial 
stratifications), and therefore 

$\bullet$ any constructible $\Lambda$-adic sheaf on $BG$ is 
lisse, given by an adic system $(M_n)_n$ of sheaves on 
$\text{Spec }k$ with $G$-actions, and these actions are 
trivial since $G$ is connected. See 
(\cite{Beh2}, 5.2.9).

(iii) Let $G_0$ be a connected $\bb F_q$-algebraic 
group. By a theorem of Serge Lang (\cite{Lan}, Th. 2), 
every $G_0$-torsor over $\text{Spec }\bb F_q$ is 
trivial, with automorphism group $G_0,$ therefore 
$$
c_v(BG_0)=\frac{1}{c_v(G_0)}=\frac{1}{\#G_0(\bb 
F_{q^v})}.  
$$
\end{blank}

Recall the following theorem of Borel as in (\cite{Beh2}, 
6.1.6).

\begin{theorem}\label{Borel}
Let $k$ be a field and $G$ a connected $k$-algebraic group. 
Consider the Leray spectral sequence given by the 
projection $f:\emph{Spec }k\to BG,$
$$
E_2^{rs}=H^r(BG_{\overline{k}})\otimes
H^s(G_{\overline{k}})\Longrightarrow\overline{\bb 
Q}_{\ell}.
$$
Let $N^s=E_{s+1}^{0,s}\subset H^s(G_{\overline{k}})$ be the 
transgressive subspaces, for $s\ge1,$ and let $N$ be the 
graded $\overline{\bb Q}_{\ell}$-vector space 
$\bigoplus_{s\ge1}N^s.$ We have 

(a). $N^s=0$ if $s$ is even,

(b). the canonical map
$\xymatrix@C=.5cm{\bigwedge N \ar[r] & 
H^*(G_{\overline{k}})}$ is an isomorphism of graded 
$\overline{\bb Q}_{\ell}$-algebras. 

(c). The spectral sequence above induces an epimorphism of
graded $\overline{\bb Q}_{\ell}$-vector spaces
$
\xymatrix@C=.5cm{
H^*(BG_{\overline{k}}) \ar@{->>}[r] & N[-1].}
$
Any section induces an isomorphism
$$
\xymatrix@C=.5cm{
\emph{Sym}^*(N[-1]) \ar[r]^-{\sim} & 
H^*(BG_{\overline{k}}).}
$$
\end{theorem}

\begin{subremark}
(i) The $E_2^{rs}$-term in (\ref{Borel}) should be 
$H^r(BG_{\overline{k}},R^sf_*\overline{\bb 
Q}_{\ell}),$ and $R^sf_*\overline{\bb Q}_{\ell}$ is a 
constructible sheaf on $BG.$ By (\ref{alg-gp}ii), the sheaf 
$R^sf_*\overline{\bb Q}_{\ell}$ is isomorphic to
$a^*f^*R^sf_*\overline{\bb 
Q}_{\ell}=a^*H^s(G_{\overline{k}}),$ where 
$a:BG\to\text{Spec }k$ is the structural map and 
$H^s(G_{\overline{k}})$ is the $\text{Gal}(k)$-module 
regarded as a sheaf on $\text{Spec }k.$ Therefore by 
projection formula, $E_2^{rs}=H^r(BG_{\overline{k}}) 
\otimes H^s(G_{\overline{k}}).$ 

(ii) Since the spectral sequence converges to 
$\overline{\bb Q}_{\ell}$ sitting in degree 0, all 
$E_{\infty}^{rs}$ are zero, except $E_{\infty}^{00}.$ For 
each $s\ge1,$ consider the differential map 
$d_{s+1}^{0,s}:E_{s+1}^{0,s}\to E_{s+1}^{s+1,0}$ on the 
$(s+1)$st page. This map must be injective (resp. 
surjective) because it is the last possibly non-zero map 
from $E_*^{0,s}$ (resp. into $E_*^{s+1,0}$). Therefore, 
it is an isomorphism. Note that $N^s=E_{s+1}^{0,s}$ is a 
subspace of $E_2^{0,s}=H^s(G_{\overline{k}}),$ and 
$E_{s+1}^{s+1,0}$ is a quotient of 
$E_2^{s+1,0}=H^{s+1}(BG_{\overline{k}}).$ Using the 
isomorphism $d_{s+1}^{0,s}$ we get the surjection 
$H^{s+1}(BG_{\overline{k}})\to N^s.$ 
\end{subremark}

\begin{anitem}\label{4.6.1}
Let $G_0$ be a connected $\mathbb F_q$-algebraic group of
dimension $d.$ We apply (\ref{Borel}) to investigate the
compact support cohomology groups $H^*_c(BG).$

We have graded Galois-invariant subspaces
$N=\bigoplus_{r\ge1}N^r\subset\bigoplus_{r\ge0}H^r(G)$
concentrated in odd degrees, such that the induced map
$$
\xymatrix@C=.5cm{
\bigwedge N \ar[r] & H^*(G)}
$$
is an isomorphism, and $H^*(BG)\cong\text{Sym}^*N[-1].$ Let
$n_r=\dim N^r,$ and let $v_{r1},\cdots,v_{rn_r}$ be a basis
for $N^r$ with respect to which the Frobenius acting on $N^r$
is upper triangular
$$
\begin{bmatrix}\alpha_{r1} & * & * \\
& \ddots & * \\
& & \alpha_{rn_r}\end{bmatrix}
$$
with eigenvalues $\alpha_{r1},\cdots,\alpha_{rn_r}.$ By
(\cite{Del2}, 3.3.5), the weights of $H^r(G)$ are $\ge r,$ so
$|\alpha_{ri}|\ge q^{r/2}>1.$ We have
$$
H^*(BG)=\text{Sym}^*\overline{\mathbb Q}_{\ell}\langle
v_{ij}|\text{for all }i,j\rangle=\overline{\mathbb
Q}_{\ell}[v_{ij}],
$$
with $\deg(v_{ij})=i+1.$ Note that all $i+1$ are even. In
particular, $H^{2r-1}(BG)=0$ and
\begin{equation*}
\begin{split}
H^{2r}(BG) &=\{\text{homogeneous polynomials of degree }2r
\text{ in }v_{ij}\} \\
&=\overline{\bb Q}_{\ell}\langle
\prod_{i,j}v_{ij}^{m_{ij}};\sum_{i,j}m_{ij}(i+1)=2r\rangle.
\end{split}
\end{equation*}
With respect to an appropriate order of the basis, the
matrix representing $F$ acting on $H^{2r}(BG)$ is upper
triangular, with eigenvalues
$$
\prod_{i,j}\alpha_{ij}^{m_{ij}},\quad\text{ for }\sum
_{i,j}m_{ij}(i+1)=2r.
$$
By Poincar\'e duality, the eigenvalues of $F$ acting on
$H^{-2r-2d}_c(BG)$ are
$q^{-d}\prod_{i,j}\alpha_{ij}^{-m_{ij}},$
for tuples of non-negative integers $(m_{ij})_{i,j}$ such
that $\sum_{i,j}m_{ij}(i+1)=2r.$ Therefore the reciprocal
characteristic polynomial of $F$ on $H^{-2r-2d}_c(BG)$ is
$$
P_{-2r-2d}(BG_0,t)=\prod_{\substack{m_{ij}\ge0 \\
\sum_{i,j}m_{ij}(i+1)=2r}}\Big(1-q^{-d}\prod_{i,j}\alpha
_{ij}^{-m_{ij}}\cdot t\Big).
$$
\end{anitem}

In the following two lemmas we prove two key cases of 
(\ref{T4.3}i).

\begin{lemma}\label{L4.10}
Let $G_0$ be an $\bb F_q$-group scheme of finite type. 
Then (\ref{T4.3}i) holds for the structural map 
$f:BG_0\to\emph{Spec }\bb F_q$ and any convergent 
complex $K_0\in W^-(BG_0,\overline{\bb Q}_{\ell}).$
\end{lemma}

\begin{proof}
By (\cite{Ols1}, 7.12-7.14) we may assume $G_0$ is reduced 
(hence smooth), so that the natural projection $\text
{Spec }\bb F_q\to BG_0$ is a presentation. Note that then 
the assumptions of $\iota$-mixedness and stratifiability 
on $K_0$ are verified automatically, by (\ref{L2.8}, 
\ref{C3.5}iii), 
even though we will not use them explicitly in the proof. 

Let $G_0^0$ be the identity component of
$G_0$ and consider the exact sequence of algebraic groups
$$
\xymatrix@C=.5cm{
1 \ar[r] & G_0^0 \ar[r] & G_0 \ar[r] & \pi_0(G_0) \ar[r] & 1.}
$$
The fibers of the induced map $BG_0\to B\pi_0(G_0)$ are
isomorphic to $BG_0^0,$ so we reduce to prove two cases:
$G_0$ is finite \'etale (or even a finite constant group
scheme, by (\ref{R4.2}iii)), or $G_0$ is connected and
smooth.

\textbf{Case of $G_0$ finite constant.} Let $G_0/\bb F_q$
be the finite constant group scheme associated with a finite
group $G,$ and let $K_0\in W^-(BG_0,\overline{\bb 
Q}_{\ell}).$ Again we denote by $G$ both the group scheme
$G_0\otimes\bb F$ and the finite group $G_0(\bb F),$
if no confusion arises. Let $y$ be the unique point in
$\text{Spec }\bb F_q,$
$$
\xymatrix@C=.5cm{
BG \ar[r] \ar[d]_{f_{\overline{y}}} & BG_0 \ar[d]^f \\
\text{Spec }\bb F \ar[r]^-{\overline{y}} & \text{Spec
}\bb F_q.}
$$
Then $D^-_c(BG,\overline{\bb Q}_{\ell})$ is
equivalent to $D^-_c(\text{Rep}_{\overline{\bb 
Q}_{\ell}}(G)),$ and the functor
$$
(f_{\overline{y}})_!:D^-_c(BG/\bb F,\overline{\bb Q}
_{\ell})\longrightarrow D^-_c(\text{Spec }\bb F,
\overline{\bb Q}_{\ell})
$$
is identified with the coinvariance functor
$$
(\ )_G:D^-_c(\text{Rep}_{\overline{\bb Q}_{\ell}}(G))
\longrightarrow D^-_c(\overline{\bb Q}_{\ell}),
$$
which is exact on the level of modules, since the category
$\text{Rep}_{\overline{\bb Q}_{\ell}}(G)$ is semisimple.
So $(f_!K_0)_{\overline{y}}
=(f_{\overline{y}})_!K=K_G$ and $\s H^n(K_G)=\s 
H^n(K)_G.$ Therefore
$$
\sum_{\s H^n((f_{\overline{y}})_!K),F}|\alpha|^s\le
\sum_{\s H^n(K),F}|\alpha|^s
$$
for every $n\in\bb Z$ and $s>0.$ Therefore, if $K_0$ is a
convergent complex, so is $f_!K_0.$

\textbf{Case of $G_0$ smooth and connected.} In this case
$$
f^*:\overline{\bb Q}_{\ell}\text{-Sh}(\text{Spec
}\bb F_q)\longrightarrow\overline{\bb 
Q}_{\ell}\text{-Sh}(BG_0)
$$
is an equivalence of categories (\ref{alg-gp}ii). Let
$d=\dim G_0,$ and let $\s F_0$ be a sheaf on $BG_0,$
corresponding to a representation $V$ of the Weil group 
$W(\bb F_q),$ with 
$\beta_1,\cdots,\beta_m$ as eigenvalues of $F.$ By the
projection formula (\cite{LO2}, 9.1.i) we have
$H^n_c(BG,\s F)\simeq H^n_c(BG)\otimes V,$ and by
(\ref{4.6.1}) the eigenvalues of $F$ on
$H_c^{-2r-2d}(BG)\otimes V$
are (using the notation in (\ref{4.6.1}))
$$
q^{-d}\beta_k\prod_{i,j}\alpha_{ij}^{-m_{ij}},
$$
for $k=1,\cdots,m$ and tuples $(m_{ij})$ such that
$\sum_{i,j}m_{ij}(i+1)=2r.$ For every $s>0,$
$$
\sum_{n\in\bb Z}\ \sum_{H^n_c(BG)\otimes V,F}|\alpha|^s=
\sum_{m_{ij},k}q^{-ds}|\beta_k|^s\prod_{i,j}|\alpha_{ij}
^{-m_{ij}}|^s=\big(\sum_{k=1}^m
|\beta_k|^s\big)\big(\sum_{m_{ij}}
q^{-ds}\prod_{i,j}|\alpha_{ij}|^{-m_{ij}s}\big),
$$
which converges to
$$
q^{-ds}\Big(\sum_{k=1}^m
|\beta_k|^s\Big)\prod_{i,j}\frac{1}{1-|\alpha_{ij}|^{-s}},
$$
since $|\alpha_{ij}|^{-s}<1$ and the product above is taken
over finitely many indices.

Let $K_0$ be a convergent complex on $BG_0,$ and for each
$k\in\bb Z,$ let $V_k$ be a $W(\bb 
F_q)$-module corresponding to $\s H^kK_0.$ For every
$x\in BG_0(\bb F_q)$ (for instance the trivial
$G_0$-torsor), the pair $(\s 
H^k(K)_{\overline{x}},F_x)$ is isomorphic to $(V_k,F).$
Consider the $W(\bb F_q)$-equivariant spectral 
sequence 
$$
H^r_c(BG,\s H^k(K))\Longrightarrow H^{r+k}_c(BG,K).
$$
We have
\begin{equation*}
\begin{split}
\sum_{n\in\bb Z}\quad\sum_{H^n_c(BG,K),F}|\alpha|^s
&\le\sum_{n\in\bb Z}\ \sum_{r+k=n}\quad
\sum_{H^r_c(BG,\s H^kK),F}|\alpha|^s
=\sum_{r,k\in\bb Z}\quad\sum_{H^r_c(BG)\otimes
V_k,F}|\alpha|^s \\
&=\sum_{k\in\bb Z}\ \sum_{r\in\bb 
Z}\quad\sum_{H^r_c(BG)\otimes
V_k,F}|\alpha|^s=\sum_{k\in\bb Z}q^{-ds}\big(
\sum_{V_k,F}|\alpha|^s\big)\prod_{i,j}\frac{1}
{1-|\alpha_{ij}|^{-s}} \\
&=\big(\sum_{k\in\bb Z}\quad\sum_{V_k,F}|\alpha|^s\big)
\big(q^{-ds}\prod_{i,j}\frac{1}{1-|\alpha_{ij}|^{-s}}\big),
\end{split}
\end{equation*}
where the first factor is convergent by assumption, and the
product in the second factor is taken over finitely many
indices. This shows that $f_!K_0$ is a convergent complex.
\end{proof}

Let $E_{\lambda}$ be a finite subextension of 
$\overline{\bb Q}_{\ell}/\bb Q_{\ell}$ with ring of 
integers $\Lambda$ and residue field $\Lambda_0,$ and let 
$(\s S,\mathcal L)$ be a pair on $\s X$ defined 
by simple lcc $\Lambda_0$-sheaves on strata. A complex 
$K_0\in W(\s X_0,\overline{\bb Q}_{\ell})$ is 
said to be \textit{$(\s S,\mathcal L)$-stratifiable} 
(or \textit{trivialized by $(\s S,\mathcal L)$}), if 
$K$ is defined over $E_{\lambda},$ with an integral model 
over $\Lambda$ trivialized by $(\s S,\mathcal L).$ 

\begin{lemma}\label{L4.11}
Let $X_0/\bb F_q$ be a geometrically connected variety,
$E_{\lambda}$ a finite subextension of $\overline{\bb 
Q}_{\ell}/\bb Q_{\ell}$ with ring of integers $\Lambda,$
and let $\mathcal L$ be a finite set of simple lcc
$\Lambda_0$-sheaves on $X.$ Then (\ref{T4.3}i) holds for the
structural map $f:X_0\to\emph{Spec }\bb F_q$ and 
all lisse $\iota$-mixed convergent complexes $K_0$ on $X_0$ 
that are trivialized by $(\{X\},\mathcal L).$
\end{lemma}

\begin{proof}
Let $N=\dim X_0.$ From the spectral sequence
$$
E_2^{rk}=H^r_c(X,\s H^kK)\Longrightarrow H^{r+k}_c(X,K)
$$
we see that
$$
\sum_{n\in\bb Z}\quad\sum_{H^n_c(X,K),F}|\alpha|^s\le
\sum_{n\in\bb Z}
\quad\sum_{r+k=n}\quad\sum_{H^r_c(X,\s H^kK),F}
|\alpha|^s=\sum_{\substack{0\le r\le2N \\ k\in\bb
Z}}\quad\sum_{H^r_c(X,\s H^kK),F}|\alpha|^s,
$$
therefore it suffices to show that, for each $0\le
r\le2N,$ the series $\sum_{k\in\bb Z}\sum_{H^r_c(X,\s 
H^kK),F}|\alpha|^s$ converges.

Let $d$ be the number in (\ref{L3.11}) for $\mathcal L.$
Each cohomology sheaf $\mathscr H^nK_0$ is $\iota$-mixed and
lisse on $X_0,$ so by (\ref{T2.4}i) we have the decomposition
$$
\s H^nK_0=\bigoplus_{b\in\bb R/\bb Z}(\s H^nK_0)(b)
$$
according to the weights mod $\bb Z,$ defined over
$E_{\lambda}$ (\ref{R2.5}ii). For each coset $b,$ we choose a
representative $b_0\in b,$ and take a
$b_1\in\overline{\bb Q}_{\ell}^*$ such that
$w_q(b_1)=-b_0.$ Then the sheaf $(\s 
H^nK_0)(b)^{(b_1)}$ deduced by twist is lisse with integer
punctual weights. Let $W$ be the filtration by punctual
weights (\ref{T2.4}ii) of $(\s H^nK_0)(b)^{(b_1)}.$ For
every $v\ge1$ and $x\in X_0(\bb F_{q^v}),$ and every real
$s>0,$ we have
\begin{equation*}
\begin{split}
\sum_{n\in\bb Z}\quad\sum_{(\s 
H^nK_0)_{\overline{x}},F_x}|\alpha|^{s/v}
&=\sum_{\substack{n\in\bb Z \\ b\in\bb R/\bb
Z}}\quad\sum_{(\s H^nK_0)(b)_{\overline{x}},F_x}
|\alpha|^{s/v} \\
&=\sum_{\substack{n\in\bb Z \\
b\in\bb R/\bb Z}}\quad\sum_{(\s 
H^nK_0)(b)^{(b_1)}_{\overline{x}},F_x}|\alpha
/b_1^v|^{s/v} \\
&=\sum_{\substack{n\in\bb Z \\ b\in\bb R/
\bb Z}}q^{b_0s/2}\sum_{(\s H^nK_0)(b)^{(b_1)}
_{\overline{x}},F_x}|\alpha|^{s/v} \\
&=\sum_{\substack{n\in\bb Z \\
b\in\bb R/\bb Z}}q^{b_0s/2}\sum_{i\in\bb 
Z}\qquad\sum_{\text{Gr}_i^W((\s 
H^nK_0)(b)^{(b_1)})_{\overline{x}},F_x}|\alpha|^{s/v} \\
&=\sum_{\substack{n\in\bb Z \\ b\in\bb R/
\bb Z}}q^{b_0s/2}\sum_{i\in\bb Z}q^{is/2}\cdot\text
{rank(Gr}_i^W((\s H^nK_0)(b)^{(b_1)})).
\end{split}
\end{equation*}
Since $K_0$ is a convergent complex, this series is
convergent.

For each $n\in\bb Z,$ every direct summand $(\s 
H^nK_0)(b)$ of $\s H^nK_0$ is trivialized
by $(\{X\},\mathcal L).$ The filtration $W$ of
each $(\s H^nK_0)(b)^{(b_1)}$ gives a filtration of
$(\s H^nK_0)(b)$ (also denoted $W$) by twisting back,
and it is clear that this latter filtration is defined over
$E_{\lambda}.$ We have $\text{Gr}^W_i((\s 
H^nK_0)(b)^{(b_1)})=(\text{Gr}^W_i(\s 
H^nK_0)(b))^{(b_1)},$ and each $\text{Gr}^W_i((\s 
H^nK_0)(b))$ is trivialized by $(\{X\},\mathcal L).$ By 
(\ref{L3.11}),
\begin{equation*}
\begin{split}
h^r_c(X,\text{Gr}^W_i((\s H^nK)(b)^{(b_1)}))
&=h^r_c(X,\text{Gr}^W_i((\s
H^nK)(b)))\qquad\qquad\text{(\cite{LO2}, 9.1.i)} \\
&\le d\cdot\text{rank}(\text{Gr}^W_i((\s H^nK)(b))) \\
&=d\cdot\text{rank}(\text{Gr}^W_i((\s H^nK)(b)^{(b_1)})).
\end{split}
\end{equation*}
Therefore
\begin{equation*}
\begin{split}
\sum_{n\in\bb Z}\quad\sum_{H^r_c(X,\s 
H^nK),F}|\alpha|^s &=\sum_{\substack{n\in\bb Z \\
b\in\bb R/\bb Z}}\quad\sum_{H^r_c(X,(\s 
H^nK)(b)),F}|\alpha|^s \\
&=\sum_{\substack{n\in\bb Z \\ b\in\bb
R/\bb Z}}\quad\sum_{H^r_c(X,(\s H^nK)(b)^{(b_1)}),
F}|b_1^{-1}\alpha|^s \\
&\le\sum_{\substack{n\in\bb Z \\ b\in\bb
R/\bb Z}}q^{b_0s/2}\sum_{i\in\bb Z}\quad\sum
_{H^r_c(X,\text{Gr}_i^W((\s
H^nK)(b)^{(b_1)})),F}|\alpha|^s \\
&\le\sum_{\substack{n\in\bb Z \\ b\in\bb
R/\bb Z}}q^{b_0s/2}\sum_{i\in\bb 
Z}q^{(i+r)s/2}\cdot h^r_c(X,\text{Gr}_i^W((\s 
H^nK)(b)^{(b_1)})) \\
&\le q^{rs/2}d\sum_{\substack{n\in\bb Z \\
b\in\bb R/\bb Z}}q^{b_0s/2}\sum_{i
\in\bb Z}q^{is/2}\cdot\text{rank(Gr}_i^W((\s 
H^nK)(b)^{(b_1)})),
\end{split}
\end{equation*}
and it converges.
\end{proof}

Now we prove (\ref{T4.3}i) in general.

\begin{proof}
We may assume all stacks involved are reduced. From
(\ref{T2.11}) and (\ref{T3.10}) we know that $f_!K_0\in
W^{-,\text{stra}}_m(\s Y_0,\overline{\bb
Q}_{\ell}).$

Let $y\in\s Y_0(\bb F_{q^v}),$ 
and we want to show that
$((f_!K_0)_{\overline{y}},F_y)$ is a convergent complex. 
Since the property of being convergent depends only on the 
cohomology sheaves, by base change (\cite{LO2}, 12.5.3) we 
reduce to the case when 
$\s Y_0=\text{Spec }\bb F_{q^v}.$ Replacing $q$ by
$q^v,$ we may assume $v=1.$ By (\ref{R4.2}iii) we only need
to show that $(R\Gamma_c(\s X,K),F)$ is convergent.

If $j:\s U_0\hookrightarrow\s X_0$ is an open
substack with complement $i:\s Z_0\hookrightarrow\s 
X_0,$ then we have an exact triangle
$$
\xymatrix@C=.5cm{
j_!j^*K_0 \ar[r] & K_0 \ar[r] & i_*i^*K_0 \ar[r] &}
$$
in $W^-(\s X_0,\overline{\bb Q}_{\ell}),$ which
induces an exact triangle
$$
\xymatrix@C=.5cm{
R\Gamma_c(\s U_0,j^*K_0) \ar[r] & R\Gamma_c(\s X_0,
K_0) \ar[r] & R\Gamma_c(\s Z_0,i^*K_0) \ar[r] &}
$$
in $W^-(\text{Spec }\bb F_q,\overline{\bb
Q}_{\ell}).$ So by (\ref{C4.6}) and noetherian induction, it
suffices to prove (\ref{T4.3}i) for a nonempty open substack.
By (\cite{Beh2}, 5.1.14) we may assume that the inertia stack
$\s I_0$ is flat over $\s X_0.$ Then we can form
the rigidification $\pi:\s X_0\to X_0$ with respect to
$\s I_0$ (\cite{Ols2}, $\S1.5),$ where $X_0$ is an
$\bb F_q$-algebraic space of quasi-compact diagonal.
$X_0$ contains an open dense subscheme (\cite{Knu},
II, 6.7). Replacing $\s X_0$ by the inverse image of
this scheme, we can assume $X_0$ is a scheme.

If (\ref{T4.3}i) holds for two composable morphisms $f$ and
$g,$ then it holds for their composition $g\circ f.$ Since
$R\Gamma_c(\mathscr X_0,-)=R\Gamma_c(X_0,-)\circ\pi_!,$ we
reduce to prove (\ref{T4.3}i) for these two morphisms. For
every $x\in X_0(\mathbb F_{q^v}),$ the fiber of $\pi$ over
$x$ is a gerbe over $\text{Spec }k(x).$ Extending the base 
$k(x)$ (\ref{R4.2}iii) one can assume it is a neutral 
gerbe (in fact all gerbes over a finite field are 
neutral; see (\cite{Beh2}, 6.4.2)). This means the 
following diagram is 2-Cartesian:
$$
\xymatrix@C=.5cm{ 
B\text{Aut}_x \ar[r] \ar[d] & \mathscr X_0 \ar[d]^-{\pi} \\ 
\text{Spec }\mathbb F_{q^v} \ar[r]^-x & X_0.} 
$$
So we reduce to two cases: $\mathscr X_0=BG_0$ for an 
$\mathbb F_q$-algebraic group $G_0,$ or $\mathscr X_0=X_0$ is 
an $\mathbb F_q$-scheme. The first case is proved in 
(\ref{L4.10}).

For the second case, given a convergent complex $K_0\in
W_m^{-,\text{stra}}(X_0,\overline{\mathbb Q}_{\ell}),$
defined over some $E_{\lambda}$ with ring of integers
$\Lambda,$ and trivialized by a pair $(\mathscr
S,\mathcal L)$ ($\mathcal L$ being defined over $\Lambda_0$)
on $X,$ we can refine this pair so that every stratum is
connected, and replace $X_0$ by models of the strata defined 
over some finite extension of $\mathbb F_q$ (\ref{R4.2}iii). 
This case is proved in (\ref{L4.11}).
\end{proof}

\section{Trace formula for stacks}

We prove two special cases of (\ref{T4.3}ii) in the following
two lemmas.

\begin{proposition}\label{L5.3}
Let $G_0$ be a finite \'etale group scheme over $\bb 
F_q,$ and $\s F_0$ a sheaf on $BG_0.$ Then
$$
c_1(BG_0,\s F_0)=c_1(\emph{Spec }\bb 
F_q,R\Gamma_c(BG_0,\s F_0)).
$$
\end{proposition}

\begin{proof}
This is a special case of (\cite{Ols1}, 8.6), on 
correspondences given by group homomorphisms, due to Olsson.
\end{proof}

\begin{proposition}\label{L5.4}
Let $G_0$ be a connected $\bb F_q$-algebraic group, and
let $\s F_0$ be a sheaf on $BG_0.$ Then
$$
c_1(BG_0,\s F_0)=c_1(\emph{Spec }\bb
F_q,R\Gamma_c(BG_0,\s F_0)).
$$
\end{proposition}

\begin{proof}
Let $f:BG_0\to\text{Spec }\bb F_q$ be the structural map
and $d=\dim G_0.$ By (\ref{alg-gp}ii), the sheaf 
$\s F_0$ on $BG_0$ takes the form $f^*V,$ for some 
sheaf $V$ on $\text{Spec }\bb F_q.$ By 
(\ref{alg-gp}iii), we have 
$$
c_1(BG_0,\s F_0)=\frac{1}{\#G_0(\bb F_q)}\text
{Tr}(F_x,\s F_{\overline{x}})=\frac{\text{Tr}(F,V)}
{\#G_0(\bb F_q)}.
$$
By the projection formula we have $H^n_c(BG,\s F)\simeq
H^n_c(BG)\otimes V,$ so 
$$
\text{Tr}(F,H^n_c(BG,\s 
F))=\text{Tr}(F,H^n_c(BG))\cdot\text{Tr}(F,V).
$$
Then
\begin{equation*}
\begin{split}
c_1(\text{Spec }\bb F_q,R\Gamma_c(BG_0,\s 
F_0)) &=\sum_n(-1)^n\text{Tr}(F,H^n_c(BG,\s F)) \\
&=\text{Tr}(F,V)\sum_n(-1)^n\text{Tr}(F,H^n_c(BG)),
\end{split}
\end{equation*}
so we can assume $\s F_0=\overline{\bb Q}_{\ell}.$
Using the notations in (\ref{4.6.1}) we have
\begin{equation*}
\begin{split}
\sum_n(-1)^n\text{Tr}(F,H^n_c(BG)) &=\sum_{r\ge0}\text{Tr}
(F,H^{-2r-2d}_c(BG))=\sum_{r\ge0}\quad\sum_{\substack{\sum 
m_{ij}(i+1)=2r \\ 
m_{ij}\ge0}}q^{-d}\prod_{i,j}\alpha_{ij}^{-m_{ij}} \\ 
&=q^{-d}\sum_{m_{ij}\ge0}\quad\prod_{i,j}\alpha_{ij}^{-m_{ij}}
=q^{-d}\prod_{i,j}(1+\alpha_{ij}^{-1}+\alpha_{ij}^{-2}+\cdots)
\\
&=q^{-d}\prod_{i,j}\frac{1}{1-\alpha_{ij}^{-1}}.
\end{split}
\end{equation*}
It remains to show
$$
\#G_0(\bb F_q)=q^d\prod_{i,j}(1-\alpha_{ij}^{-1}).
$$

In (\ref{4.6.1}), we saw that if each $N^i$ has an ordered 
basis $v_{i1},\cdots,v_{in_i}$ with respect to which $F$ is upper
triangular, then since $H^*(G)=\bigwedge N,\ H^i(G)$ has basis
$$
v_{i_1j_1}\wedge v_{i_2j_2}\wedge\cdots\wedge v_{i_mj_m},
$$
such that $\sum_{r=1}^mi_r=i,\ i_r\le i_{r+1},$ and if
$i_r=i_{r+1},$ then $j_r<j_{r+1}.$ The eigenvalues of $F$ on
$H^i(G)$ are $\alpha_{i_1j_1}\cdots\alpha_{i_mj_m}$ for such
indices. By Poincar\'e duality, the eigenvalues of $F$ on
$H^{2d-i}_c(G)$ are $q^d(\alpha_{i_1j_1}\cdots\alpha_{i_mj_m})^{-1}.$ 
Note that all the $i_r$ are odd, so
$$
2d-i\equiv i=\sum_{r=1}^mi_r\equiv m\mod2.
$$
Applying the classical trace formula to $G_0,$ we have
$$
\#G_0(\bb F_q)=\sum(-1)^mq^d\alpha_{i_1j_1}^{-1}\cdots
\alpha_{i_mj_m}^{-1}=q^d\prod_{i,j}(1-\alpha_{ij}^{-1}).
$$
This finishes the proof.
\end{proof}

\begin{anitem}\label{5.5}
Note that, in (\ref{L5.3}) and (\ref{L5.4}) we did not make 
explicit use of the fact that $\s F_0$ is $\iota$-mixed.
\end{anitem}

Now we prove (\ref{T4.3}ii) in general.

\begin{proof}
Since $c_v(\s X_0,K_0)=c_1(\s X_0\otimes\bb
F_{q^v},K_0\otimes\bb F_{q^v}),$ we can assume $v=1.$
We shall reduce to proving (\ref{T4.3}ii) for all fibers of
$f$ over $\bb F_q$-points of $\s Y_0,$ following
the approach of Behrend (\cite{Beh2}, 6.4.9).

Let $y\in\s Y_0(\bb F_q)$ and $(\s X_0)_y$
be the fiber over $y.$ Then $(\s X_0)_y(\bb F_q)$ is
the groupoid of pairs $(x,\alpha),$ where $x\in\s
X_0(\bb F_q)$ and $\alpha:f(x)\to y$ is an isomorphism in
$\s Y_0(\bb F_q).$ Suppose $(\s X_0)_y(\bb 
F_q)\ne\emptyset,$ and fix an $x\in(\s X_0)_y(\bb
F_q).$ Then $\text{Isom}(f(x),y)(\bb F_q)$ is a trivial left
$\text{Aut}_y(\bb F_q)$-torsor. There is also a natural
right action of $\text{Aut}_x(\bb F_q)$ on
$\text{Isom}(f(x),y)(\bb F_q),$ namely
$\varphi\in\text{Aut}_x(\bb F_q)$ takes $\alpha$ to
$\alpha\circ f(\varphi).$ Under this action, for
$\alpha$ and $\alpha'$ to be in the same orbit, there
should be a $\varphi\in\text{Aut}_x(\bb F_q)$ such
that the diagram
$$
\xymatrix@C=.6cm{
f(x) \ar[rr]^-{f(\varphi)} \ar[rd]_-{\alpha'} && f(x)
\ar[ld]^-{\alpha} \\ & y &}
$$
commutes, and this is the definition for $(x,\alpha)$ to be
isomorphic to $(x,\alpha')$ in $(\s X_0)_y(\bb
F_q).$ So the set of orbits $\text{Isom}(f(x),y)(\bb
F_q)/\text{Aut}_x(\bb F_q)$ is identified with the
inverse image of the class of $x$ along the map $[(\s
X_0)_y(\bb F_q)]\to[\s X_0(\bb F_q)].$ The
stabilizer group of $\alpha\in\text{Isom}(f(x),y)(\bb
F_q)$ is $\text{Aut}_{(x,\alpha)}(\bb F_q),$ the
automorphism group of $(x,\alpha)$ in $(\s
X_0)_y(\bb F_q).$

In general, if a finite group $G$ acts on a finite set $S,$
then we have
$$
\sum_{[x]\in S/G}\frac{\#G}{\#\text{Stab}_G(x)}=\sum
_{[x]\in S/G}\#\text{Orb}_G(x)=\#S.
$$
Now for $S=\text{Isom}(f(x),y)(\bb F_q)$ and
$G=\text{Aut}_x(\bb F_q),$ we have
$$
\sum_{\substack{(x,\alpha)\in[(\s X_0)_y(\bb F_q)]
\\ (x,\alpha)\mapsto x}}\frac{\#\text{Aut}_x(\bb
F_q)}{\#\text{Aut}_{(x,\alpha)}(\bb F_q)}=\#\text{Isom}
(f(x),y)(\bb F_q)=\#\text{Aut}_y(\bb F_q);
$$
the last equality follows from the fact that $S$ is a trivial
$\text{Aut}_y(\bb F_q)$-torsor.

If we assume (\ref{T4.3}ii) holds for the fibers 
$f_y:(\s X_0)_y\to\text{Spec }\bb F_q$ of $f,$ 
for all $y\in\s Y_0(\bb F_q),$ then 
\begin{equation*}
\begin{split}
c_1(\s Y_0,f_!K_0)
&=\sum_{y\in[\s Y_0(\bb F_q)]}\frac{\text{Tr}(F_y,
(f_!K)_{\overline{y}})}{\#\text{Aut}_y(\bb F_q)} \\
&=\sum_{y\in[\s Y_0(\bb F_q)]}\frac{\text{Tr}(F_y,
(f_{y!}K)_{\overline{y}})}{\#\text{Aut}_y(\bb F_q)} 
\qquad\qquad(\cite{LO2}, 12.5.3) \\ 
&=\sum_{y\in[\s Y_0(\bb F_q)]}\frac{1}{\#\text{Aut}
_y(\bb F_q)}\sum_{(x,\alpha)\in[(\s X_0)_y(\bb
F_q)]}\frac{\text{Tr}(F_x,K_{\overline{x}})}{\#\text
{Aut}_{(x,\alpha)}(\bb F_q)} \\
&=\sum_{y\in[\s Y_0(\bb F_q)]}\frac{1}{\#\text
{Aut}_y(\bb F_q)}\sum_{\substack{x\in[\s X_0(\bb
F_q)] \\ x\mapsto y}}\Big(\sum_{\substack{(x,\alpha)\in[(\s
X_0)_y(\bb F_q)] \\ (x,\alpha)\mapsto x}}\frac{\text
{Tr}(F_x,K_{\overline{x}})}{\#\text{Aut}_{(x,\alpha)}(\bb
F_q)}\Big) \\
&=\sum_{y\in[\s Y_0(\bb F_q)]}\frac{1}{\#\text{Aut}_y
(\bb F_q)}\sum_{\substack{x\in[\s X_0(\bb F_q)] \\
x\mapsto y}}\frac{1}{\#\text{Aut}_x(\bb F_q)} \\
& \qquad\qquad\qquad\qquad\qquad\Big(\sum
_{\substack{(x,\alpha)\in[(\s X_0)_y(\bb F_q)] \\
(x,\alpha)\mapsto x}}\frac{\#\text{Aut}_x(\bb F_q)}{\#
\text{Aut}_{(x,\alpha)}(\bb F_q)}\Big)\text{Tr}(F_x,K
_{\overline{x}}) \\
&=\sum_{y\in[\s Y_0(\bb F_q)]}\frac{1}{\#\text{Aut}
_y(\bb F_q)}\sum_{\substack{x\in[\s X_0(\bb
F_q)] \\ x\mapsto y}}\frac{\text{Tr}(F_x,K_{\overline{x}})}
{\#\text{Aut}_x(\bb F_q)}\#\text{Aut}_y(\bb F_q) \\
&=\sum_{x\in[\s X_0(\bb F_q)]}\frac{\text{Tr}(F_x,
K_{\overline{x}})}{\#\text{Aut}_x(\bb F_q)}=:c_1(\s
X_0,K_0).
\end{split}
\end{equation*}
So we reduce to the case when $\mathscr Y_0=\text{Spec
}\mathbb F_q.$

If $K_0'\to K_0\to K_0''\to K_0'[1]$ is an exact triangle of
convergent complexes in $W_m^{-,\text{stra}}(\s 
X_0,\overline{\bb Q}_{\ell}),$ then by (\ref{C4.6}) and
(\ref{T4.3}i) we have
$$
c_1(\s X_0,K_0)=c_1(\s X_0,K_0')+c_1(\s X_0,K_0'')
$$
and
$$
c_1(\s Y_0,f_!K_0)=c_1(\s Y_0,f_!K_0')+c_1
(\s Y_0,f_!K_0'').
$$
If $j:\s U_0\to\s X_0$ is an open substack with
complement $i:\s Z_0\to\s X_0,$ then
$$
c_1(\s X_0,j_!j^*K_0)=c_1(\s U_0,j^*K_0)\ \text{
and }\ c_1(\s X_0,i_*i^*K_0)=c_1(\s Z_0,i^*K_0).
$$
By noetherian induction we can shrink $\s X_0$ to a
nonempty open substack. So as before we may assume the
inertia stack $\s I_0$ is flat over $\s X_0,$
with rigidification $\pi:\s X_0\to X_0,$ where $X_0$ is
a scheme. If (\ref{T4.3}ii) holds for two composable
morphisms $f$ and $g,$ then it holds for $g\circ f.$ So we
reduce to two cases as before: $\s X_0=X_0$ is a
scheme, or $\s X_0=BG_0,$ where $G_0$ is either a
connected algebraic group, or a finite \'etale algebraic
group over $\bb F_q.$ We may assume $X_0$ is separated,
by further shrinking (for instance to an affine open
subscheme).

For a complex of sheaves $K_0$ and an integer $n,$ we have 
an exact triangle
$$
\xymatrix@C=.5cm{
\tau_{<n}K_0 \ar[r] & \tau_{<n+1}K_0 \ar[r] & \s 
H^n(K_0)[-n] \ar[r] &,}
$$
so
\begin{equation*}
\begin{split}
c_1(\tau_{<n+1}K_0) &=c_1(\tau_{<n}K_0)+c_1(\s 
H^n(K_0)[-n]) \\
&=c_1(\tau_{<n}K_0)+(-1)^nc_1(\s H^n(K_0)).
\end{split}
\end{equation*}
Since $K_0$ is bounded above, $\tau_{<N}K_0\simeq
K_0$ for $N\gg0.$ Since $K_0$ is convergent, $c_1(\tau_{<
n}K_0)\to0$ absolutely as $n\to-\infty,$ so the series
$\sum_{n\in\bb Z}(-1)^nc_1(\s H^n(K_0))$ converges
absolutely to $c_1(K_0).$

Applying $R\Gamma_c$ we get an exact triangle
$$
\xymatrix@C=.5cm{
R\Gamma_c(\s X_0,\tau_{<n}K_0) \ar[r] &
R\Gamma_c(\s X_0,\tau_{<n+1}K_0) \ar[r] &
R\Gamma_c(\s X_0,\s H^nK_0)[-n] \ar[r] &}
$$
in $W^-(\text{Spec }\bb F_q,\overline{\bb
Q}_{\ell}).$ We claim that, for $\s X_0=X_0$ a scheme,
or $BG_0,$ we have
$$
\lim_{n\to-\infty}c_1(\text{Spec }\bb
F_q,R\Gamma_c(\s X_0,\tau_{<n}K_0))=0
$$
absolutely. Recall that $c_1(R\Gamma_c(\tau_{<n}K_0))=
\sum_{i\in \bb Z}(-1)^i\iota\text{Tr}(F,H^i_c(\s 
X,\tau_{<n}K)),$ so we need to show
$$
\sum_{i\in\bb Z}\quad\sum_{H^i_c(\s X,\tau_{<
n}K),F}|\alpha|\to0\qquad\text{as}\qquad n\to-\infty.
$$

From the spectral sequence
$$
H^r_c(\s X,\s H^k\tau_{<n}K)
\Longrightarrow H^{r+k}_c(\s X,\tau_{<n}K)
$$
we see that
\begin{equation*}
\begin{split}
\sum_{i\in\bb Z}\quad\sum_{H^i_c(\s X,\tau_{<
n}K),F}|\alpha| &\le\sum_{i\in\bb
Z}\quad\sum_{r+k=i}\quad\sum_{H^r_c(\s X,\s 
H^k\tau_{<n}K),F}|\alpha| \\
&=\sum_{i\in\bb Z}\quad\sum_{\substack{r+k=i \\ k<
n}}\quad\sum_{H^r_c(\s X,\s H^kK),F}|\alpha|.
\end{split}
\end{equation*}
Let $d=\dim\s X_0$ (cf. \ref{dim}). In the cases where
$\s X_0$ is a scheme or $BG_0,$ we have $H^r_c(\s 
X,\s F)=0$ for every sheaf $\s F$ unless $r\le2d$
(cf. (\ref{4.6.1}) and (\ref{L4.8})). Therefore
$$
\sum_{i\in\bb Z}\quad\sum_{\substack{r+k=i \\ k<n}}\quad
\sum_{H^r_c(\s X,\s H^kK),F}|\alpha|\le
\sum_{i<n+2d}\quad\sum_{r+k=i}\quad\sum_{H^r_c(\s 
X,\s H^kK),F}|\alpha|,
$$
and it suffices to show that the series
$$
\sum_{i\in\bb Z}\quad\sum_{r+k=i}\qquad\sum_{H^r_c
(\s X,\s H^kK),F}|\alpha|
$$
converges. This is proved for $BG_0$ in (\ref{L4.10}), and
for schemes $X_0$ in (\ref{L4.11}) (we may shrink $X_0$ so
that the assumption in (\ref{L4.11}) is satisfied).

Note that in the two cases $\s X_0=X_0$ or $BG_0,$ 
(\ref{T4.3}ii) holds when $K_0$ is a sheaf
concentrated in degree 0. For separated schemes $X_0,$ this
is a classical result of Grothendieck and Verdier \cite{Gro,
Ver}; for $BG_0,$ this is done in (\ref{L5.3}) and
(\ref{L5.4}). Therefore, for a general convergent complex
$K_0,$ we have
\begin{equation*}
\begin{split}
c_1(R\Gamma_c(\tau_{<n+1}K_0)) &=c_1(R\Gamma_c(\tau_{<
n}K_0))+c_1(R\Gamma_c(\s H^nK_0)[-n]) \\
&=c_1(R\Gamma_c(\tau_{<n}K_0))+(-1)^nc_1(\s H^nK_0),
\end{split}
\end{equation*}
and so
$$
c_1(R\Gamma_c(K_0))=\sum_{n\in\bb Z}(-1)^nc_1(\s 
H^nK_0)+\lim_{n\to-\infty}c_1(R\Gamma_c(\tau_{<n}K_0))=
c_1(K_0).
$$
\end{proof}

\begin{corollary}\label{C5.5}
Let $f:\s X_0\to\s Y_0$ be a morphism of 
$\bb F_q$-algebraic stacks, and let $K_0\in 
W_m^{-,\emph{stra}}(\s X_0,\overline{\bb 
Q}_{\ell})$ be a convergent complex of sheaves. Then
$$
L(\s X_0,K_0,t)=L(\s Y_0,f_!K_0,t).
$$
\end{corollary}

\section{Infinite products}

For a convergent complex $K_0$ on $\s X_0,$ the series
$\sum_{v\ge1}c_v(K_0)t^v/v$ (and hence the $L$-series 
$L(\s X_0,K_0,t)$) usually has a finite radius of 
convergence. For instance, we have the following lemma. 

\begin{lemma}\label{radius}
Let $X_0/\bb F_q$ be a variety of dimension $d.$ Then 
the radius of convergence of $\sum_{v\ge1}c_v(X_0)t^v/v$ 
is $1/q^d.$ 
\end{lemma}

\begin{proof}
Let $f_{X_0}(t)=\sum_{v\ge1}c_v(X_0)t^v/v.$ Let $Y_0$ be an 
irreducible component of $X_0$ with complement $U_0.$ Then 
$c_v(X_0)=c_v(Y_0)+c_v(U_0),$ and since all the $c_v$-terms 
are non-negative, we see that the radius of convergence of 
$f_{X_0}(t)$ is the minimum of that of $f_{Y_0}(t)$ and 
that of $f_{U_0}(t).$ Since 
$\max\{\dim(Y_0),\dim(U_0)\}=d,$ and $U_0$ has fewer 
irreducible component than $X_0,$ by induction we can 
assume $X_0$ is irreducible. 

Then there exists an open dense subscheme $U_0\subset X_0$ 
that is smooth over $\text{Spec }\bb F_q.$ Let 
$Z_0=X_0-U_0,$ then $\dim(Z_0)<\dim(X_0)=d.$ From the 
cohomology sequence 
$$
\xymatrix@C=.5cm{ 
H^{2d-1}_c(Z) \ar[r] & H^{2d}_c(U) \ar[r] & H^{2d}_c(X) 
\ar[r] & H^{2d}_c(Z)}
$$
we see that $H^{2d}_c(X)=H^{2d}_c(U)=\overline{\bb 
Q}_{\ell}(-d).$ The Frobenius eigenvalues 
$\{\alpha_{ij}\}_j$ on $H^i_c(X)$ have $\iota$-weights 
$\le i,$ for $0\le i<2d$ (\cite{Del2}, 3.3.4). By the fixed 
point formula, 
$$
\frac{c_v(X_0)}{c_{v+1}(X_0)}=\frac{q^{vd}+\sum_{0\le 
i<2d}(-1)^i\sum_j\alpha_{ij}^v}{q^{(v+1)d}+\sum_{0\le 
i<2d}(-1)^i\sum_j\alpha_{ij}^{v+1}}=\frac{\frac{1}{q^d}+ 
\frac{1}{q^d}\sum_{0\le i<2d}(-1)^i\sum_j(\frac{\alpha 
_{ij}}{q^d})^v} {1+\sum_{0\le i<2d}(-1)^i\sum_j(\frac 
{\alpha_{ij}}{q^d})^{v+1}},
$$
which converges to $1/q^d$ as $v\to\infty,$ therefore the 
radius of convergence of $f_{X_0}(t)$ is 
$$
\lim_{v\to\infty}\frac{c_v(X_0)/v}{c_{v+1}(X_0)/(v+1)}= 
\frac{1}{q^d}.
$$
\end{proof}

In order to prove the meromorphic
continuation (\ref{T8.1}), we want to express the $L$-series
as a possibly infinite product. For schemes, if we consider
only bounded complexes, the $L$-series can be expressed as
a finite alternating product of polynomials
$P_n(X_0,K_0,t),$ so it is rational \cite{Gro}. In the stack
case, even for the sheaf $\overline{\bb Q}_{\ell},$
there might be infinitely many nonzero compact cohomology
groups, and we need to consider the issue of convergence of
the coefficients in an infinite products.

\begin{definition}\label{D6.1}
Let $f_n(t)=\sum_{k\ge0}a_{nk}t^k\in\bb C[[t]]$ be a
sequence of power series over $\bb C.$ The sequence
is said to be \emph{convergent term by term}, if for each
$k,$ the sequence $(a_{nk})_n$ converges, and the series
$$
\lim_{n\to\infty}f_n(t):=\sum_{k\ge0}t^k\lim_{n\to\infty}a_{nk}
$$
is called the limit of the sequence $\big(f_n(t)\big)_n.$
\end{definition}

\begin{anitem}\label{}
Strictly speaking, a series (resp. infinite product) is
defined to be a sequence $(a_n)_n,$ usually written as an
``infinite sum" (resp. ``infinite product") so that $(a_n)_n$
is the sequence of finite partial sums (resp. finite partial
products) of it. So the definition above applies to series
and infinite products.
\end{anitem}

Recall that $\log(1+g)=\sum_{m\ge1}(-1)^{m+1}g^m/m$ for 
$g\in t\mathbb C[[t]].$ 

\begin{lemma}\label{L6.2}
(i) Let $f_n(t)=1+\sum_{k\ge1}a_{nk}t^k\in\bb C[[t]]$ be a
sequence of power series. Then $\big(f_n(t)\big)_n$ is
convergent term by term if and only if $\big(\log
f_n(t)\big)_n$ is convergent term by term, and
$$
\lim_{n\to\infty}\log f_n(t)=\log\lim_{n\to\infty}f_n(t).
$$

(ii) Let $f$ and $g$ be two power series with constant term 1.
Then
$$
\log(fg)=\log(f)+\log(g).
$$

(iii) Let $f_n(t)\in1+t\bb C[[t]]$ be a sequence as in
(i). Then the infinite product $\prod_{n\ge1}f_n(t)$
converges term by term if and only if the series
$\sum_{n\ge1}\log f_n(t)$ converges term by term, and
$$
\sum_{n\ge1}\log f_n(t)=\log\prod_{n\ge1}f_n(t).
$$
\end{lemma}

\begin{proof}
(i) We have 
\begin{equation*}
\begin{split}
\log f_n(t) &=\sum_{m\ge1}(-1)^{m+1}
\big(\sum_{k\ge1}a_{nk}t^k\big)^m/m \\
&=t\cdot a_{n1}+t^2(a_{n2}-\frac{a_{n1}^2}{2})+t^3(a_{n3}
-a_{n1}a_{n2}+\frac{a_{n1}^3}{3}) \\
& \qquad+t^4(a_{n4}-a_{n1}a_{n3}
-\frac{a_{n2}^2}{2}+a_{n1}^2a_{n2})+\cdots \\
&=:\sum_{k\ge1}A_{nk}t^k.
\end{split}
\end{equation*}
In particular, for each $k,\ A_{nk}-a_{nk}=h(a_{n1},\cdots,
a_{n,k-1})$ is a polynomial in $a_{n1},\cdots,a_{n,k-1}$ with
rational coefficients. So if $(a_{nk})_n$ converges for each
$k,$ then $(A_{nk})_n$ also converges, and by induction the
converse also holds. If $\lim_{n\to\infty}a_{nk}=a_k,$ then
$\lim_{n\to\infty}A_{nk}=a_k+h(a_1,\cdots,a_{k-1}),$
and
$$
\log\lim_{n\to\infty}f_n(t)=\log(1+\sum_{k\ge1}a_kt^k)=
\sum_{k\ge1}(a_k+h(a_1,\cdots,a_{k-1}))t^k=\lim_{n\to\infty}
\log f_n(t).
$$

(ii) $\log$ and $\exp$ are inverse to each other on power
series, so it suffices to prove that for $f$ and $g\in
t\bb C[[t]],$ we have
$$
\exp(f+g)=\exp(f)\exp(g).
$$
This follows from the binomial formula:
\begin{equation*}
\begin{split}
\exp(f+g) &=\sum_{n\ge0}(f+g)^n/n!=\sum_{n\ge0}\frac{1}{n!}
\sum_{k=0}^n\binom{n}{k}f^kg^{n-k}=\sum_{n\ge0}\quad
\sum_{k=0}^n\frac{f^k}{k!}\cdot\frac{g^{n-k}}{(n-k)!} \\
&=\sum_{i,j\ge0}
\frac{f^i}{i!}\cdot\frac{g^j}{j!}=\big(\sum_{i\ge0}f^i/i!
\big)\big(\sum_{j\ge0}g^j/j!\big)=\exp(f)\exp(g).
\end{split}
\end{equation*}

(iii) Let $F_N(t)=\prod_{n=1}^Nf_n(t).$ Applying (i) to the
sequence $(F_N(t))_N,$ we see that the infinite product
$\prod_{n\ge1}f_n(t)$ converges term by term if and only if
(by definition) $\big(F_N(t)\big)_N$ converges term by term,
if and only if the sequence $\big(\log F_N(t)\big)_N$
converges term by term, if and only if (by definition) the
series $\sum_{n\ge1}\log f_n(t)$ converges term by term, since
by (ii)
$$
\log\prod_{n=1}^Nf_n(t)=\sum_{n=1}^N\log f_n(t)
$$
Also
$$
\log\prod_{n\ge1}f_n(t)=\log\lim_{N\to\infty}F_N(t)=
\lim_{N\to\infty}\log F_N(t)=\lim_{N\to\infty}\sum_{n=1}
^N\log f_n(t)=:\sum_{n\ge1}\log f_n(t).
$$
\end{proof}

\begin{blank}\label{}
For a complex of sheaves $K_0$ on $\s X_0$ and
$n\in\bb Z,$ define
$$
P_n(\s X_0,K_0,t):=\det(1-Ft,H^n_c(\s X,K)).
$$
We regard $P_n(\s X_0,K_0,t)^{\pm1}$ as a complex power
series with constant term 1 via $\iota.$
\end{blank}

\begin{proposition}\label{P6.5}
For every convergent complex of sheaves $K_0\in
W_m^{-,\emph{stra}}(\s X_0,\overline{\bb
Q}_{\ell}),$ the infinite product
$$
\prod_{n\in\bb Z}P_n(\s X_0,K_0,t)^{(-1)^{n+1}}
$$
is convergent term by term to the $L$-series $L(\s 
X_0,K_0,t).$
\end{proposition}

\begin{proof}
The complex $R\Gamma_c(\s X,K)$ is bounded above, so
$P_n(\s X_0,K_0,t)=1$ for $n\gg0,$ and the infinite
product is a one-direction limit, namely $n\to-\infty.$

Let $\alpha_{n1},\cdots,\alpha_{nm_n}$ be the eigenvalues
(counted with multiplicity) of $F$ on $H^n_c(\s 
X,K),$ regarded as complex numbers via $\iota,$ so that
$$
P_n(t)=P_n(\s 
X_0,K_0,t)=(1-\alpha_{n1}t)\cdots(1-\alpha_{nm_n}t).
$$
By (\ref{L6.2}iii) it suffices to show that the series
$$
\sum_{n\in\bb Z}(-1)^{n+1}\log P_n(t)
$$
converges term by term to $\sum_{v\ge1}c_v(K_0)t^v/v.$

We have
\begin{equation*}
\begin{split}
\sum_{n\in\bb Z}(-1)^{n+1}\log P_n(t) &=\sum_{n\in\bb
Z}(-1)^{n+1}\log\prod_i(1-\alpha_{ni}t)=\sum_{n\in\bb
Z}(-1)^n\sum_i\sum_{v\ge1}\frac{\alpha_{ni}^vt^v}{v} \\
&=\sum_{v\ge1}\frac{t^v}{v}\sum_{n\in\bb
Z}(-1)^n\sum_i\alpha_{ni}^v=\sum_{v\ge1}\frac{t^v}{v}c_v
(R\Gamma_c(K_0)),
\end{split}
\end{equation*}
which converges term by term (\ref{T4.3}i), and is equal
(\ref{T4.3}ii) to $\sum_{v\ge1}c_v(K_0)t^v/v.$
\end{proof}

\begin{subremark}\label{R6.6}
In particular we have
$$
Z(\s X_0,t)=\prod_{n\in\bb Z}P_n(\s X_0,t)^{(-1)^{n+1}},
$$
where $P_n(\s X_0,t)=P_n(\s X_0,\overline 
{\bb Q}_{\ell},t).$ This generalizes the 
classical result for schemes (\cite{Gro}, 5.1). When
we want to emphasize the dependence on the prime $\ell,$ we
will write $P_{n,\ell}(\s X_0,t).$ 

If $G_0$ is a connected $\bb F_q$-algebraic group,
(\ref{4.6.1}) shows that the zeta function of $BG_0$ is given
by
\begin{equation*}
\begin{split}
Z(BG_0,t) &=\prod_{r\ge0}\quad\prod_{\substack{m_{ij}\ge0 \\
\sum_{i,j}m_{ij}(i+1)=2r}}\Big(1-q^{-d}\prod
_{i,j}\alpha_{ij}^{-m_{ij}}\cdot t\Big)^{-1} \\
&=\prod_{m_{ij}\ge0}\Big(1-q^{-d}\prod_{i,j}\alpha_{ij}^
{-m_{ij}}\cdot t\Big)^{-1}.
\end{split}
\end{equation*}
\end{subremark}

\section{Examples of zeta functions}

In this section we compute the zeta functions of some stacks,
and in each example we do it in two ways: counting rational
points and computing cohomology groups. Also we investigate
some analytic properties.

\begin{example}\label{example1}
$B\bb G_m.$ By (\ref{alg-gp}iii) we have $c_v(B\bb
G_m)=1/c_v(\bb G_m),$ so the zeta function is
$$
Z(B\bb G_m,t)=\exp\Big(\sum_{v\ge1}c_v(B\bb G_m)\frac
{t^v}{v}\Big)=\exp\Big(\sum_{v\ge1}\frac{1}{q^v-1}\frac{t^v}
{v}\Big).
$$
Using Borel's theorem (\ref{Borel}) one can show (or see
(\cite{LMB}, 19.3.2)) that the cohomology ring $H^*(B\bb
G_m)$ is a polynomial ring $\overline{\bb Q}_{\ell}[T],$
generated by a variable $T$ of degree 2, and the global
Frobenius action is given by $FT^n=q^nT^n.$ So by
Poincar\'e duality, we have
\begin{equation*}
\begin{split}
\text{Tr}(F,H^{-2n-2}_c(B\bb G_m))
&=\text{Tr}(F,H^{-2n-2}_c(B\bb G_m,\overline{\bb
Q}_{\ell}(-1)))/q \\
&=\text{Tr}(F^{-1},H^{2n}(B\bb G_m))/q=q^{-n-1}.
\end{split}
\end{equation*}
This gives
$$
\prod_{n\in\bb Z}P_n(B\bb G_m,t)^{(-1)^{n+1}}=
\prod_{n\ge1}(1-q^{-n}t)^{-1}.
$$
It is easy to verify (\ref{R6.6}) directly:
\begin{equation*}
\begin{split}
\exp\Big(\sum_{v\ge1}\frac{1}{q^v-1}\frac{t^v}{v}\Big) &=
\exp\Big(\sum_{v\ge1}\frac{1/q^v}{1-1/q^v}\frac{t^v}{v}
\Big)=\exp\Big(\sum_{v\ge1}\frac{t^v}{v}\sum_{n\ge1}\frac{1}
{q^{nv}}\Big) \\
&=\prod_{n\ge1}\exp\Big(\sum_{v\ge1}\frac{(t/q^n)^v}{v}\Big)
=\prod_{n\ge1}(1-t/q^n)^{-1}.
\end{split}
\end{equation*}

There is also a functional equation
$$
Z(B\bb G_m,qt)=\frac{1}{1-t}Z(B\bb G_m,t),
$$
which implies that $Z(B\bb G_m,t)$ has a meromorphic
continuation to the whole complex plane, with simple poles at
$t=q^n,$ for $n\ge1.$

$H^{-2n-2}_c(B\bb G_m)$ is pure of weight $-2n-2.$ A
natural question is whether Deligne's theorem of weights
(\cite{Del2}, 3.3.4) still holds for algebraic stacks. Olsson
told me that it does not hold in general, as the following
example shows.
\end{example}

\begin{example}\label{example2}
$BE,$ where $E$ is an elliptic curve over $\bb F_q.$
Again by (\ref{alg-gp}iii) we have 
$$
c_v(BE)=\frac{1}{\#E(\mathbb F_{q^v})}.
$$
Let $\alpha$ and $\beta$ be the roots of the reciprocal
characteristic polynomial of the Frobenius on $H^1(E):$
\begin{equation}\label{E2}
x^2-(1+q-c_1(E))x+q=0.
\end{equation}
Then for every $v\ge1,$ we have
$c_v(E)=1-\alpha^v-\beta^v+q^v=(1-\alpha^v)(1-\beta^v).$ So
\begin{equation*}
\begin{split}
c_v(BE) &=\frac{1}{(1-\alpha^v)(1-\beta^v)}=\frac{\alpha^{-v}}
{1-\alpha^{-v}}\cdot\frac{\beta^{-v}}{1-\beta^{-v}} \\
&=\big(\sum_{n\ge1}\alpha^{-nv}\big)\big(\sum_{m\ge1}\beta^
{-nv}\big)=\sum_{n,m\ge1}\Big(\frac{1}{\alpha^n\beta^m}\Big)^v,
\end{split}
\end{equation*}
and the zeta function is
$$
Z(BE,t)=\exp\Big(\sum_{v\ge1}c_v(BE)\frac{t^v}{v}\Big)=
\exp\Big(\sum_{\substack{n,m\ge1 \\ v\ge1}}
\big(\frac{t}{\alpha^n\beta^m}\big)^v/v\Big)=
\prod_{n,m\ge1}\big(1-\frac{t}{\alpha^n\beta^m}\big)^{-1}.
$$

To compute its cohomology, one can apply Borel's theorem
(\ref{Borel}) to $E,$ and we have $N=N^1=H^1(E),$ so $N[-1]$
is a 2-dimensional vector space sitting in degree 2, on which
$F$ has eigenvalues $\alpha$ and $\beta.$ Then $H^*(BE)$ is a
polynomial ring $\overline{\bb Q}_{\ell}[a,b]$ in two
variables, both sitting in degree 2, and the basis $a,b$ can
be chosen so that the Frobenius action $F$ on $H^2(BE)$ is
upper triangular (or even diagonal)
$$
\begin{bmatrix}\alpha & \gamma \\
& \beta\end{bmatrix}.
$$
Then $F$ acting on
$$
H^{2n}(BE)=\text{Sym}^nN[-1]=\overline{\bb
Q}_{\ell}\langle a^n,a^{n-1}b,\cdots,b^n\rangle
$$
can be represented by
$$
\begin{bmatrix}\alpha^n & * & * & * \\
& \alpha^{n-1}\beta & * & * \\
&& \ddots & * \\
&&& \beta^n\end{bmatrix},
$$
with eigenvalues $\alpha^n,\alpha^{n-1}\beta,\cdots,\beta^n.$
So the eigenvalues of $F$ on $H^{-2-2n}_c(BE)$ are
$$
q^{-1}\alpha^{-n},q^{-1}\alpha^{1-n}\beta^{-1},\cdots,
q^{-1}\beta^{-n},
$$
and $\prod_{n\in\bb Z}P_n(BE,t)^{(-1)^{n+1}}$ is
$$
\frac{1}{(1-q^{-1}t)[(1-q^{-1}\alpha^{-1}t)(1-q^{-1}
\beta^{-1}t)][(1-q^{-1}\alpha^{-2}t)(1-q^{-1}\alpha^{-1}
\beta^{-1}t)(1-q^{-1}\beta^{-2}t)]\cdots},
$$
which is the same as $Z(BE,t)$ above (since $\alpha\beta=q$).

Let $Z_1(t):=Z(BE,qt).$ Its radius of convergence is 1, 
since by (\ref{radius}) 
$$
\lim_{v\to\infty}\frac{c_v(BE)}{c_{v+1}(BE)}=\lim_{v\to
\infty}\frac{c_{v+1}(E)}{c_v(E)}=q.
$$
There is also a functional equation
$$
Z_1(\alpha t)=\frac{1}{1-\alpha t}Z_1(t)Z_2(t),
$$
where
$$
Z_2(t)=\frac{1}{(1-\alpha\beta^{-1}t)(1-\alpha\beta^{-2}t)
(1-\alpha\beta^{-3}t)\cdots}.
$$
$Z_2(t)$ is holomorphic in the open unit disk and satisfies
the functional equation
$$
Z_2(\beta t)=\frac{1}{1-\alpha t}Z_2(t).
$$
Therefore $Z_2(t),$ and hence $Z(BE,t),$ has a meromorphic
continuation to the whole complex $t$-plane with the obvious
poles.

\begin{subremark}\label{R7.1}
$H^{-2-2n}_c(BE)$ is pure of weight $-2-n,$ which is not
$\le-2-2n$ unless $n=0.$ So the upper bound of weights for
schemes fails for $BE.$ This also leads to the failure of 
the decomposition theorem for $BE;$ see (\cite{Decom}, 
$\S1$) for the example of a pure complex on $BE$ which is 
not geometrically semi-simple. 

Also note that, the equation
(\ref{E2}) is independent of $\ell,$ so the polynomials 
$P_{n,\ell}(BE,t)$ are independent of $\ell.$
\end{subremark}
\end{example}

\begin{example}\label{example3}
$BG_0,$ where $G_0$ is a finite \'etale $\bb F_q$-group 
scheme, corresponding to a finite group $G$ and a Frobenius 
automorphism $\sigma$ on it. Then we have $BG_0(\bb 
F_{q^v})\simeq 
G/\rho^{(v)},$ where $\rho^{(v)}$ is the right action of $G$ 
on the set $G$ given by $h:g\mapsto\sigma^v(h^{-1})gh.$ So 
$$
c_v(BG_0)=\sum_{[g]\in G/\rho^{(v)}}\frac{1}{\#\text
{Stab}_{\rho^{(v)}}(g)}=\frac{\#G}{\#G}=1,
$$
and the zeta function is
$$
Z(BG_0,t)=\frac{1}{1-t}.
$$
Its cohomology groups are given in (\ref{L4.8}):
$H^0_c(BG)=\overline{\mathbb Q}_{\ell},$ and other $H^i_c=0.$
This verifies (\ref{R6.6}).

Note that $Z(BG_0,t)$ is the same as the zeta function
of its coarse moduli space $\text{Spec }\mathbb F_q.$ As a
consequence, for every $\mathbb F_q$-algebraic stack $\mathscr
X_0,$ with finite inertia $\mathscr I_0\to\mathscr X_0$ and
coarse moduli space $\pi:\mathscr X_0\to X_0$ (\cite{Crd},
1.1), we have $Z(\mathscr X_0,t)=Z(X_0,t),$ and hence it is a
rational function. This is because for every $x\in X_0(\mathbb
F_{q^v}),$ the fiber $\pi^{-1}(x)$ is a neutral gerbe over
$\text{Spec }k(x),$ and from the above we see that
$c_v(\pi^{-1}(x))=1,$ and hence $c_v(\mathscr X_0)=c_v(X_0).$
The fact that $Z(X_0,t)$ is a rational function follows from
(\cite{Knu}, II, 6.7) and noetherian induction. More
generally, we have the following.

\begin{subproposition}\label{C7.2}
Let $\mathscr X_0$ be an $\mathbb F_q$-algebraic stack. 
Suppose that $\mathscr X_0$ either has finite inertia, or 
is Deligne-Mumford (not necessarily separated). Then for 
every $K_0\in W^b(\mathscr X_0,\overline{\mathbb 
Q}_{\ell}),$ the $L$-series $L(\mathscr X_0,K_0,t)$ is a 
rational function.
\end{subproposition}

\begin{proof}
It suffices to show that (\ref{T4.3}) holds for the structural 
map $\s X_0\to\text{Spec }\bb F_q$ and $K_0\in 
W^b(\s X_0,\overline{\bb Q}_{\ell})$ in these two cases. 
We will not make explicit use of the fact (\ref{Laff}) that 
$K_0$ is $\iota$-mixed. 

\textbf{Case when $\mathscr X_0$ has finite inertia.} Let 
$\pi:\mathscr X_0\to X_0$ be its coarse moduli space. 
For any sheaf $\mathscr F_0$ on $\mathscr X_0,$ by 
(\ref{L4.8}) we have isomorphisms $H^r_c(X,R^0\pi_!\mathscr 
F)\simeq H^r_c(\mathscr{X,F}),$ so $R\Gamma_c(\mathscr
X_0,\mathscr F_0)$ is a bounded complex, hence a convergent
complex. To prove the trace formula for $\mathscr X_0\to\text
{Spec }\mathbb F_q$ and the sheaf $\mathscr F_0,$ it suffices
to prove it for $\mathscr X_0\to X_0$ and $X_0\to\text{Spec
}\mathbb F_q.$ The first case, when passing to fibers, is
reduced to $BG_0,$ and when passing to fibers again, it is
reduced to the two subcases: when $G_0$ is finite, or when
$G_0$ is connected. In both of these two cases as well as the
case of an algebraic space $X_0\to\text{Spec }\mathbb F_q,$ the
trace formula can be proved without using $\iota$-mixedness 
(\ref{5.5}). Therefore, (\ref{T4.3}) holds 
for $\mathscr X_0\to\text{Spec }\mathbb F_q$ and any sheaf,
hence any bounded complex, on $\mathscr X_0.$

The trace formula is equivalent to the equality of power
series
$$
L(\mathscr X_0,K_0,t)=\prod_{i\in\mathbb Z}P_i(\mathscr
X_0,K_0,t)^{(-1)^{i+1}},
$$
and the right-hand side is a finite product, because 
$R\Gamma_c(\mathscr X_0,K_0)$ is bounded. 
Therefore, $L(\mathscr X_0,K_0,t)$ is rational.

\textbf{Case when $\mathscr X_0$ is Deligne-Mumford.} For 
both (i) and (ii) of (\ref{T4.3}), we may replace $\mathscr 
X_0$ by a non-empty open substack, hence by (\cite{LMB}, 
6.1.1) we may assume $\mathscr X_0$ is the quotient stack 
$[X_0'/G],$ where $X_0'$ is an affine $\mathbb F_q$-scheme 
of finite type and $G$ is a finite group acting on $X_0'.$ 
This stack has finite diagonal, and hence finite inertia, 
so by the previous case we are done. Also, we know that 
$R\Gamma_c(\mathscr X_0,K_0)$ is bounded, therefore 
$L(\mathscr X_0,K_0,t)$ is rational. 
\end{proof}

If one wants to use Poincar\'e duality to get a functional
equation for the zeta function, (\cite{Ols1}, 5.17) and
(\cite{LO2}, 9.1.2) suggest that we should assume $\mathscr
X_0$ to be proper smooth and of finite diagonal. Under these
assumptions, one gets the expected functional equation for
the zeta function, as well as the independence of $\ell$ for
the coarse moduli space, which is proper but possibly
singular. Examples of such stacks include the moduli stack 
of pointed stable curves $\overline{\mathscr
M}_{g,n}$ over $\mathbb F_q.$

\begin{subproposition}\label{T7.3}
Let $\mathscr X_0$ be a proper smooth $\mathbb F_q$-algebraic
stack of equidimension $d,$ with finite diagonal, and let
$\pi:\mathscr X_0\to X_0$ be its coarse moduli space. Then
$Z(X_0,t)$ satisfies the usual functional equation
$$
Z(X_0,\frac{1}{q^dt})=\pm q^{d\chi/2}t^{\chi}Z(X_0,t),
$$
where $\chi:=\sum_{i=0}^{2d}(-1)^i\deg P_{i,\ell}(X_0,t).$
Moreover, $H^i(X)$ is pure of weight $i,$ for every
$0\le i\le2d,$ and the reciprocal roots of each
$P_{i,\ell}(X_0,t)$ are algebraic integers independent of
$\ell.$
\end{subproposition}

\begin{proof}
First we show that the adjunction map $\overline{\mathbb
Q}_{\ell}\to\pi_*\pi^*\overline{\mathbb
Q}_{\ell}=\pi_*\overline{\mathbb Q}_{\ell}$ is an isomorphism.
Since $\pi$ is quasi-finite and proper (\cite{Crd}, 1.1),
we have $\pi_*=\pi_!$ (\cite{Ols1}, 5.1) and
$R^r\pi_!\overline{\mathbb Q}_{\ell}=0$ for $r\ne0$
(\cite{Ols1}, 5.8). The natural map $\overline{\mathbb
Q}_{\ell}\to R^0\pi_*\overline{\mathbb Q}_{\ell}$ is an
isomorphism, since the geometric fibers of $\pi$ are
connected.

Therefore $R\Gamma(\mathscr X_0,\overline{\mathbb
Q}_{\ell})=R\Gamma(X_0,\pi_*\overline{\mathbb Q}_{\ell})
=R\Gamma(X_0,\overline{\mathbb Q}_{\ell}),$ and hence
(\cite{Ols1}, 5.17) $H^i(\mathscr X)\simeq H^i_c(\mathscr
X)\simeq H^i(X)\simeq H^i_c(X)$ for all $i.$ Let $P_i(t)=
P_i(\mathscr X_0,t)=P_i(X_0,t).$ Since $X_0$ is an algebraic
space of dimension $d,\ P_i(t)=1$ if $i\notin[0,2d].$ Since
$\mathscr X_0$ is proper and smooth, Poincar\'e duality gives
a perfect pairing
$$
\xymatrix@C=.5cm{
H^i(\mathscr X)\times H^{2d-i}(\mathscr X) \ar[r] &
\overline{\mathbb Q}_{\ell}(-d).}
$$
Following the standard proof for proper smooth varieties (as
in (\cite{Mil1}, 27.12)) we get the expected functional
equation for $Z(\mathscr X_0,t)=Z(X_0,t).$

$H^i(X)$ is mixed of weights $\le i$ (\cite{Del2}, 3.3.4),
so by Poincar\'e duality, it is pure of weight $i.$ Following
the proof in (\cite{Del1}, p.276), this purity implies that
the polynomials $P_{i,\ell}(X_0,t)$ 
have integer coefficients independent of $\ell.$
\end{proof}

\begin{subremark}\label{R7.4}
Weizhe Zheng suggested (\ref{C7.2}) to me. He also suggested
that we give a functional equation relating $L(\mathscr
X_0,DK_0,t)$ and $L(\mathscr X_0,K_0,t),$ for $K_0\in
W^b(\mathscr X_0,\overline{\mathbb Q}_{\ell}),$ where
$\mathscr X_0$ is a proper $\mathbb F_q$-algebraic stack with
finite diagonal, of equidimension $d,$ but not necessarily
smooth. Here is the functional equation:
$$
L(\mathscr X_0,K_0,t^{-1})=t^{\chi_c}\cdot Q\cdot L(\mathscr
X_0,DK_0,t),
$$
where $\chi_c=\sum_{i=0}^{2d}(-1)^ih^i_c(\mathscr X,K)$ and
$Q=(t^{\chi_c}L(\mathscr X_0,K_0,t))|_{t=\infty}.$ Note that
the rational function $L(\mathscr X_0,K_0,t)$ has degree
$-\chi_c,$ hence $Q$ is well-defined. The proof is similar to
the above.
\end{subremark}
\end{example}

\begin{example}\label{example4}
$BGL_N.$ We have $\#GL_N(\mathbb
F_{q^v})=(q^{vN}-1)(q^{vN}-q^v)\cdots(q^{vN}-q^{v(N-1)}),$ so
one can use $c_v(BGL_N)=1/c_v(GL_N)$ to compute $Z(BGL_N,t).$
One can also compute the cohomology groups of $BGL_N$ using
Borel's theorem (\ref{Borel}). We refer to (\cite{Beh1},
2.3.2) for the result. Let us consider the case $N=2$ only.
The general case is similar.

We have
$$
c_v(BGL_2)=\frac{1}{q^{4v}}\bigg(1+\frac{1}{q^v}+\frac{2}
{q^{2v}}+\frac{2}{q^{3v}}+\frac{3}{q^{4v}}+\frac{3}{q^{5v}}
+\cdots\bigg),
$$
and therefore
\begin{equation*}
\begin{split}
Z(BGL_2,t) &=\exp\Big(\sum_v\frac{(t/q^4)^v}{v}\Big)\cdot
\exp\Big(\sum_v\frac{(t/q^5)^v}{v}\Big)\cdot\exp\Big(\sum_v
\frac{2(t/q^6)^v}{v}\Big)\cdots \\
&=\frac{1}{1-t/q^4}\cdot\frac{1}{1-t/q^5}\cdot\Big(\frac{1}
{1-t/q^6}\Big)^2\cdot\Big(\frac{1}{1-t/q^7}\Big)^2\cdot
\Big(\frac{1}{1-t/q^8}\Big)^3\cdots.
\end{split}
\end{equation*}
So $Z(BGL_2,qt)=Z(BGL_2,t)\cdot Z_1(t),$ where
$$
Z_1(t)=\frac{1}{(1-t/q^3)(1-t/q^5)(1-t/q^7)(1-t/q^9)\cdots}.
$$
$Z_1(t)$ satisfies the functional equation
$$
Z_1(q^2t)=\frac{1}{1-t/q}\cdot Z_1(t),
$$
So $Z_1(t),$ and hence $Z(BGL_2,t),$ has a meromorphic
continuation with the obvious poles.

The non-zero compactly supported cohomology groups of $BGL_2$
are given as follows:
$$
H^{-8-2n}_c(BGL_2)=\overline{\mathbb Q}_{\ell}(n+4)^{\oplus
\big(\left\lfloor\frac{n}{2}\right\rfloor+1\big)},\ n\ge0.
$$
This gives
$$
\prod_{n\in\mathbb Z}P_n(BGL_2,t)^{(-1)^{n+1}}
=\frac{1}{(1-t/q^4)(1-t/q^5)(1-t/q^6)^2(1-t/q^7)^2\cdots},
$$
and (\ref{R6.6}) is verified. Note that the eigenvalues are
$1/q^{n+4},$ which are independent of $\ell.$
\end{example}

\section{Analytic continuation}

We state and prove a generalized version of (\ref{T1.3}).

\begin{theorem}\label{T8.1}
Let $\mathscr X_0$ be an $\mathbb F_q$-algebraic stack, and
let $K_0\in W_m^{-,\emph{stra}}(\mathscr 
X_0,\overline{\mathbb Q}_{\ell})$ be a convergent complex. 
Then $L(\mathscr X_0,K_0,t)$ has a meromorphic continuation 
to the whole complex $t$-plane, and its poles can only be 
zeros of the polynomials $P_{2n}(\mathscr X_0,K_0,t)$ for 
some integers $n.$
\end{theorem}

We need a preliminary lemma. For an open subset $U\subset
\mathbb C,$ let $\mathscr O(U)$ be the set of analytic
functions on $U.$ There exists a sequence $\{K_n\}_{n\ge1}$ of
compact subsets of $U$ such that $U=\bigcup_nK_n$ and
$K_n\subset(K_{n+1})^{\circ}.$ For $f$ and $g$ in $\mathscr
O(U),$ define
$$
\rho_n(f,g)=\sup\{|f(z)-g(z)|;z\in K_n\}\quad\text{and}
\quad\rho(f,g)=\sum_{n=1}^{\infty}\Big(\frac{1}{2}\Big)^n
\frac{\rho_n(f,g)}{1+\rho_n(f,g)}.
$$
Then $\rho$ is a metric on $\mathscr O(U)$ and the topology
is independent of the subsets $\{K_n\}_n$ chosen
(cf. (\cite{Con}, VII, $\S1$)).

The following lemma is (\cite{Con}, p.167, 5.9).

\begin{lemma}\label{L8.2}
Let $U\subset\mathbb C$ be connected and open and let
$(f_n)_n$ be a sequence in $\mathscr O(U)$ such that no $f_n$
is identically zero. If $\sum_n(f_n(z)-1)$ converges
absolutely and uniformly on compact subsets of $U,$ then
$\prod_{n\ge1}f_n(z)$ converges in $\mathscr O(U)$ to an
analytic function $f(z).$ If $z_0$ is a zero of $f,$ then
$z_0$ is a zero of only a finite number of the functions
$f_n,$ and the multiplicity of the zero of $f$ at $z_0$ is
the sum of the multiplicities of the zeros of the functions
$f_n$ at $z_0.$
\end{lemma}

Now we prove (\ref{T8.1}).

\begin{proof}
Factorize $P_n(\mathscr X_0,K_0,t)$ as
$\prod_{j=1}^{m_n}(1-\alpha_{nj}t)$ in $\mathbb C.$ Since
$R\Gamma_c(\mathscr X_0,K_0)$ is a convergent complex
(\ref{T4.3}i), the series $\sum_{n,j}|\alpha_{nj}|$
converges.

By (\ref{P6.5}) we have
$$
L(\mathscr X_0,K_0,t)=\prod_{n\in\mathbb Z}\big(\prod_{j=1}
^{m_n}(1-\alpha_{nj}t)\big)^{(-1)^{n+1}}
$$
as formal power series. To apply (\ref{L8.2}), take $U$ to be
the region $\mathbb C-\{\alpha_{nj}^{-1};n\text{ even}\}.$
Take the lexicographical order on the set of all factors
$$
1-\alpha_{nj}t,\text{ for }n\text{ odd};\
\frac{1}{1-\alpha_{nj}t},\text{ for }n\text{ even}.
$$
Each factor is an analytic function on $U.$ The sum
$``\sum_n(f_n(z)-1)"$ here is equal to
$$
\sum_{n\text{ odd},j}(-\alpha_{nj}t)+\sum_{n\text{ even},j}
\frac{\alpha_{nj}t}{1-\alpha_{nj}t}.
$$
Let
$$
g_n(t)=\begin{cases}\sum_{j=1}^{m_n}|\alpha_{nj}t|,\
n\text{ odd,} \\ \sum_{j=1}^{m_n}\frac{|\alpha_{nj}t|}
{|1-\alpha_{nj}t|},\ n\text{ even.}\end{cases}
$$
We need to show that $\sum_ng_n(t)$ is pointwise convergent,
uniformly on compact subsets of $U.$ Precisely, we want
to show that for any compact subset $B\subset U,$ and for any
$\varepsilon>0,$ there exists a constant $N_B\in\mathbb Z$
such that
$$
\sum_{n\le N}g_n(t)<\varepsilon
$$
for all $N\le N_B$ and $t\in B.$ Since $g_n(t)$ are
non-negative, it suffices to do this for $N=N_B.$ There
exists a constant $M_B$ such that $|t|<M_B$ for all $t\in B.$
Since $\sum_{n,j}|\alpha_{nj}|$ converges,
$|\alpha_{nj}|\to0$ as $n\to-\infty,$ and there exists a
constant $L_B\in\mathbb Z$ such that $|\alpha_{nj}|<1/(2M_B)$
for all $n<L_B.$ So
$$
g_n(t)\le2\sum_{j=1}^{m_n}|\alpha_{nj}t|
$$
for all $n<L_B$ and $t\in B.$ There exists a constant
$N_B<L_B$ such that
$$
\sum_{n\le N_B}\quad\sum_j|\alpha_{nj}|<\varepsilon/(2M_B)
$$
and hence
$$
\sum_{n\le N_B}g_n(t)\le2\sum_{n\le N_B}\
\sum_j|\alpha_{nj}t|\le
2M_B\sum_{n\le N_B}\ \sum_j|\alpha_{nj}|<\varepsilon.
$$
By (\ref{L8.2}), $L(\mathscr X_0,K_0,t)$ extends to an
analytic function on $U.$ By the second part of (\ref{L8.2}),
the $\alpha_{nj}^{-1}$'s, for $n$ even, are at worst poles
rather than essential singularities, therefore the $L$-series
is meromorphic on $\mathbb C.$
\end{proof}

Now $L(\mathscr X_0,K_0,t)$ can be called an 
``$L$-function". 

\section{Weight theorem for algebraic stacks}

\begin{blank}\label{dim}
We prove (\ref{T1.4}) in this section. For the reader's
convenience, we briefly review the definition of the
\textit{dimension} of a locally noetherian $S$-algebraic
stack $\mathcal X$ from (\cite{LMB}, chapter 11).

If $X$ is a locally noetherian $S$-algebraic space and $x$ is
a point of $X,$ the dimension $\dim_x(X)$ of $X$ at $x$ is
defined to be $\dim_{x'}(X'),$ for any pair $(X',x')$ 
where $X'$ is an $S$-scheme \'etale over $X$ and $x'\in
X'$ maps to $x.$ This is independent of the choice of the
pair. If $f:X\to Y$ is a morphism of $S$-algebraic spaces,
locally of finite type, and $x$ is a point of $X$ with image
$y$ in $Y,$ then the relative dimension $\dim_x(f)$ of $f$ at
$x$ is defined to be $\dim_x(X_y).$

Let $P:X\to\mathcal X$ be a presentation of an $S$-algebraic
stack $\mathcal X,$ and let $x$ be a point of $X.$ Then the
relative dimension $\dim_x(P)$ of $P$ at $x$ is defined to be
the relative dimension at $(x,x)$ of the smooth morphism of
$S$-algebraic spaces $\text{pr}_1:X\times_{\mathcal X}X\to
X.$

If $\mathcal X$ is a locally noetherian $S$-algebraic stack
and if $\xi$ is a point of $\mathcal X,$ the dimension of 
$\mathcal X$ at $\xi$ is defined to be $\dim_{\xi}(\mathcal
X)=\dim_x(X)-\dim_x(P),$ where $P:X\to\mathcal X$ is an
arbitrary presentation of $\mathcal X$ and $x$ is an
arbitrary point of $X$ lying over $\xi.$ This definition is
independent of all the choices made. At last one defines 
$\dim\mathcal X=\sup_{\xi}\dim_{\xi}\mathcal X.$ For quotient
stacks we have $\dim[X/G]=\dim X-\dim G.$
\end{blank}

Now we prove (\ref{T1.4}).

\begin{proof}
If $j:\mathscr U_0\to\mathscr X_0$ is an open substack with
complement $i:\mathscr Z_0\to\mathscr X_0,$ then we have an
exact sequence
$$
\xymatrix@C=.5cm{
\cdots \ar[r] & H^n_c(\mathscr U,j^*\mathscr F) \ar[r] &
H^n_c(\mathscr X,\mathscr F) \ar[r] & H^n_c(\mathscr
Z,i^*\mathscr F) \ar[r] & \cdots}.
$$
If both $H^n_c(\mathscr U,j^*\mathscr F)$ and $H^n_c(\mathscr
Z,i^*\mathscr F)$ are zero (resp. have all 
$\iota$-weights $\le m$ for some number $m$), then so is
$H^n_c(\mathscr X,\mathscr F).$ Since the dimensions of
$\mathscr U_0$ and $\mathscr Z_0$ are no more than that of
$\mathscr X_0,$ and the set of punctual $\iota$-weights of
$i^*\mathscr F_0$ and of $j^*\mathscr F_0$ is the same as 
that of $\mathscr F_0,$ we may shrink $\mathscr X_0$ to a 
non-empty open substack. We can also make any finite base 
change on $\mathbb F_q.$ To simplify notation, we may use 
twist (\ref{Weil-cplx}) and projection formula to assume 
$w=0.$ As before,
we reduce to the case when $\mathscr X_0$ is geometrically
connected, and the inertia $f:\mathscr I_0\to\mathscr X_0$ is
flat, with rigidification $\pi:\mathscr X_0\to X_0,$ where
$X_0$ is a scheme. The squares in the following diagram are
2-Cartesian:
$$
\xymatrix@C=.7cm{
&& \mathscr I_0 \ar[d]_f & \text{Aut}_y \ar[l] \ar[d] \\
B\text{Aut}_{\overline{x}} \ar[r]
\ar[d]_{\pi_{\overline{x}}} & B\text{Aut}_x \ar[r]
\ar[d]_{\pi_x} & \mathscr X_0 \ar[d]^{\pi} & \text{Spec }
\mathbb F_{q^v} \ar[l]_-y \\
\text{Spec }\mathbb F \ar[r]_-{\overline{x}} &
\text{Spec }\mathbb F_{q^v} \ar[r]_-x & X_0 &.}
$$
We have $(R^k\pi_!\mathscr F_0)_{\overline{x}}=H^k_c(B
\text{Aut}_{\overline{x}},\mathscr F).$ Since $f$ is
representable and flat, and $\mathscr X_0$ is connected,
all automorphism groups $\text{Aut}_x$ have the same
dimension, say $d.$

Assume (\ref{T1.4}) holds for all $BG_0,$ where $G_0$ are
$\mathbb F_q$-algebraic groups. Then $R^k\pi_!\mathscr
F_0=0$ for $k>-2d,$ and for $k\le-2d,$ the punctual
$\iota$-weights of $R^k\pi_!\mathscr F_0$ are
$\le\frac{k}{2}-d,$ hence by (\cite{Del2}, 3.3.4), the
punctual $\iota$-weights of $H^r_c(X,R^k\pi_!\mathscr
F)$ are $\le\frac{k}{2}-d+r.$ Consider the Leray spectral
sequence
$$
E_2^{rk}=H^r_c(X,R^k\pi_!\mathscr F)
\Longrightarrow H^{r+k}_c(\mathscr X,\mathscr F).
$$
If we maximize $\frac{k}{2}-d+r$ under the constraints
$$
r+k=n,\ 0\le r\le2\dim X_0,\text{ and }k\le-2d,
$$
we find that $H^n_c(\mathscr X,\mathscr F)=0$ for $n>
2\dim X_0-2d=2\dim\mathscr X_0,$ and for $n\le2\dim\mathscr
X_0,$ the punctual $\iota$-weights of $H^n_c(\mathscr
X,\mathscr F)$ are $\le\dim X_0+\frac{n}{2}-d=\dim\mathscr
X_0+\frac{n}{2}.$

So we reduce to the case $\mathscr X_0=BG_0.$ The Leray
spectral sequence for $h:BG_0\to B\pi_0(G_0)$ degenerates
(by (\ref{L4.8})) to isomorphisms
$$
H^0_c(B\pi_0(G),R^nh_!\mathscr F)\simeq H^n_c(BG,\mathscr
F).
$$
The fibers of $h$ are isomorphic to $BG_0^0,$ so by base
change and (\ref{L4.8}) we reduce to the case when $G_0$ is
connected. Let $d=\dim G_0$ and $f:BG_0\to\text{Spec
}\mathbb F_q$ be the structural map. In this case, $\mathscr
F_0\cong f^*V$ for some $\overline{\mathbb
Q}_{\ell}$-representation $V$ of $W(\mathbb F_q),$
and hence $\mathscr F_0$ and $V$ have the same punctual
$\iota$-weights. Using the natural isomorphism
$H^n_c(BG)\otimes V\simeq H^n_c(BG,\mathscr F),$
we reduce to the case when $\mathscr F_0=\overline{\mathbb
Q}_{\ell}.$ In (\ref{4.6.1}) we see that, if
$\alpha_{i1},\cdots,
\alpha_{in_i}$ are the eigenvalues of $F$ on $N^i,\ i\ge1$
odd, then the eigenvalues of $F$ on $H^{-2k-2d}_c(BG)$ are
$$
q^{-d}\prod_{i,j}\alpha_{ij}^{-m_{ij}},\ \text{where }
\sum_{i,j}m_{ij}(i+1)=2k.
$$
Since $i\ge1,$ we have $\sum im_{ij}\ge k;$ together with
$|\alpha_{ij}|\ge q^{i/2}$ one deduces
$$
|q^{-d}\prod_{i,j}\alpha_{ij}^{-m_{ij}}|\le
q^{\frac{-k-2d}{2}},
$$
so the punctual $\iota$-weights of $H^{-2k-2d}_c(BG)$ are
$\le-k-2d$ for $k\ge0,$ and the other compactly supported
cohomology groups are zero.

It is clear from the proof and (\cite{Del2}, 3.3.10) that the
weights of $H^n_c(\mathscr X,\mathscr F)$ differ from the
weights of $\mathscr F_0$ by integers.

Recall that $H^{2k}(BG)$ is pure of weight $2k,$ for a linear
algebraic group $G_0$ over $\mathbb F_q$ (\cite{Del3}, 9.1.4).
Therefore, $H^{-2k-2d}_c(BG)$ is pure of weight $-2k-2d,$ and
following the same proof as above, we are done.
\end{proof}

\begin{remark}\label{R9.1}
When $\mathscr X_0=X_0$ is a scheme, and $n\le2\dim X_0,$ we
have $\dim X_0+\frac{n}{2}+w\ge n+w,$ so our bound for weights
is worse than the bound in (\cite{Del2}, 3.3.4). For an
$\mathbb F_q$-abelian variety $A,$ our bound for the weights of
$H^n_c(BA)$ is sharp: the weights are exactly
$\dim(BA)+\frac{n}{2},$ whenever the cohomology group is
non-zero.
\end{remark}

We hope (\ref{T1.4}) has some useful and interesting
applications, for instance for generalizing the 
decomposition theorem of Beilinson-Bernstein-Deligne-Gabber (cf. 
\cite{Decom}) to stacks with affine 
stabilizers, and for studying the Hasse-Weil zeta functions of
Artin stacks over number fields. For instance, it 
implies that the Hasse-Weil zeta function is analytic in 
some right half complex $s$-plane. 

Using (\ref{T1.4}) we can show certain stacks have $\mathbb
F_q$-points.

\begin{example}\label{example9.2}
Let $\s X_0$ be a form of $B\bb G_m,$ i.e.,
$\s X\cong B\bb G_{m,\bb F}$ over $\bb F.$ 
Then all the 
automorphism group schemes in $\s X_0$ are affine, and
$h^{-2-2n}_c(\s X)=h^{-2-2n}_c(B\bb G_m)=1,$ for all
$n\ge0.$ Let $\alpha_{-2-2n}$ be the eigenvalue of $F$ on
$H^{-2-2n}_c(\s X).$ Then by (\ref{T1.4}) we have
$|\alpha_{-2-2n}|\le q^{-1-n}.$ Smoothness is fppf local on
the base, so $\s X_0$ is smooth and connected, hence
$H^{-2}_c(\s X)=\overline{\bb Q}_{\ell}(1)$ and
$\alpha_{-2}=q^{-1}.$ So
\begin{equation*}
\begin{split}
\#\s X_0(\bb F_q) &=\sum_{n\ge0}\text{Tr}(F,
H^{-2-2n}_c(\s X))=q^{-1}+\alpha_{-4}+\alpha_{-6}+
\cdots \\
&\ge q^{-1}-q^{-2}-q^{-3}+\cdots=q^{-1}-\frac{q^{-1}}{q-1}>0
\end{split}
\end{equation*}
when $q\ne2.$ In fact, since there exists an integer 
$r\ge1$ such that $\mathscr X_0\otimes\mathbb F_{q^r}\cong 
B\bb G_{m,\bb F_{q^r}},$ we see that all 
cohomology groups $H^{-2-2n}_c(\s X)$ are pure, i.e. 
$|\alpha_{-2-2n}|=q^{-1-n}.$

In fact, one can classify the forms of $B\bb 
G_{m,\bb F_q}$ as follows. If $\s X_0$ is a form, 
then it is also a gerbe over $\text{Spec }\bb F_q,$ 
hence a neutral gerbe $BG_0$ for some algebraic group 
$G_0$ by (\cite{Beh2}, 6.4.2). By comparing the 
automorphism groups, we see that $G_0$ is a form of 
$\bb G_{m,\bb F_q}.$ There is only 
one nontrivial form of $\bb G_{m,\bb F_q},$ because 
$$
H^1(\bb F_q,\text{Aut}(\bb G_m))=H^1(\bb 
F_q,\bb Z/2\bb Z)=\bb Z/2\bb Z, 
$$
and this form is the kernel $R^1_{\bb F_{q^2}/\bb 
F_q}\bb G_{m,\bb F_{q^2}}$ of the norm map 
$$
\xymatrix@C=1cm{
R_{\bb F_{q^2}/\bb F_q}\bb G_{m,\bb F_{q^2}} 
\ar[r]^-{\text{Nm}} & \bb G_{m,\bb F_q},}
$$
where $R_{\bb F_{q^2}/\bb F_q}$ is the operation of 
Weil's restriction of scalars. Therefore, the only 
non-trivial form of $B\bb G_{m,\bb F_q}$ is 
$B(R^1_{\bb F_{q^2}/\bb F_q}\bb G_{m,\bb 
F_{q^2}}).$ In particular, they all have $\bb 
F_q$-points, even when $q=2.$
\end{example}

\begin{example}\label{example9.3}
Consider the projective line $\bb P^1$ with the following 
action of $\bb G_m:$ it acts by multiplication on the open 
part $\bb A^1\subset\bb P^1,$ and leaves the point 
$\infty$ fixed. So we get a quotient stack
$[\bb P^1/\bb G_m]$ over $\bb F_q.$ Let $\s 
X_0$ be a form of $[\bb P^1/\bb G_m].$ 
We want to find an $\bb F_q$-point on $\s 
X_0,$ or even better, an $\bb F_q$-point on $\s X_0$
which, when considered as a point in $\s X(\bb
F)\cong[\bb P^1/\bb G_m](\bb F),$ lies in the open
dense orbit $[\bb G_m/\bb G_m](\bb F).$

\begin{anitem}
Consider the following general situation. Let $G_0$ be a
connected $\bb F_q$-algebraic group, and let $X_0$ be a
proper smooth variety with a $G_0$-action over $\bb
F_q.$ Let
$$
\xymatrix@C=.5cm{[X_0/G_0] \ar[r]^-f & BG_0 \ar[r]^-g &
\text{Spec }\bb F_q}
$$
be the natural maps, and let $\s X_0$ be a form of
$[X_0/G_0].$ Then $f$ is representable and proper. For every
$k,\ R^kf_*\overline{\bb Q}_{\ell}$ is a lisse sheaf, and
takes the form $g^*V_k$ for some sheaf $V_k$ on $\text{Spec
}\bb F_q.$ Consider the Leray spectral sequence
$$
E^{rk}_2=R^rg_!R^kf_*\overline{\bb
Q}_{\ell}\Longrightarrow R^{r+k}(gf)_!\overline{\bb
Q}_{\ell}.
$$
Since $R^rg_!R^kf_*\overline{\bb Q}_{\ell}=R^rg_!(g^*
V_k)=(R^rg_!\overline{\bb Q}_{\ell})\otimes V_k,$ we
have
$$
h^n_c(\s X)=h^n_c([X/G])\le\sum_{r+k=n}h^r_c(BG)\cdot
\dim V_k=\sum_{r+k=n}h^r_c(BG)\cdot h^k(X).
$$
\end{anitem}

Now we return to $[\bb P^1/\bb G_m].$ Since $h^0(\bb
P^1)=h^2(\bb P^1)=1$ and $h_c^{-2i}(B\bb G_m)=1$ for
$i\ge1,$ we see that $h_c^n(\s X)=0$ for $n$ odd and
$$
h^{2n}_c(\s X)\le h^0(\bb P^1)h^{2n}_c(B\bb
G_m)+h^2(\bb P^1)h^{2n-2}_c(B\bb G_m)=\begin{cases}
0,\ n\ge1, \\ 1,\ n=0, \\ 2,\ n<0.\end{cases}
$$
Since $\s X_0$ is connected and smooth of dimension 
0, we have $H^0_c(\s X)=\overline{\bb Q}_{\ell}.$ By
(\ref{T1.4}), the $\iota$-weights of
$H^{2n}_c(\s X)$ are $\le2n.$ The trace formula gives
\begin{equation*}
\begin{split}
\#\s X_0(\bb F_q) &=\sum_{n\le0}\text{Tr}(F,H^{2n}_c(\s 
X))=1+\sum_{n<0}\text{Tr}(F,H^{2n}_c(\s X)) \\ 
&\ge1-2\sum_{n<0}q^n=1-\frac{2}{q-1}>0
\end{split}
\end{equation*}
when $q\ge4.$

In order for the rational point to be in the open dense
orbit, we need an upper bound for the number of $\bb 
F_q$-points on the closed orbits. When passing to $\bb 
F,$ there are 2 closed orbits, both having stabilizer 
$\bb G_{m,\bb F}.$ So in 
$[\s X_0(\bb F_q)]$ there are at most 2 points
whose automorphism groups are forms of the algebraic group
$\bb G_{m,\bb F_q}.$ From the cohomology sequence
$$
\xymatrix@C=.7cm{
1 \ar[r] & (R^1_{\bb F_{q^2}/\bb F_q}\bb
G_{m,\bb F_{q^2}})(\bb F_q) \ar[r] & \bb 
F_{q^2}^* \ar[r]^-{\text{Nm}} & \bb F_q^*}
$$
we see that $\#(R^1_{\bb F_{q^2}/\bb F_q}\bb
G_{m,\bb F_{q^2}})(\bb F_q)=q+1.$ Since $\frac{1}{q+1}
\le\frac{1}{q-1},$ the space that the closed orbits can take 
is at most $\frac{2}{q-1},$ and equality holds only when the 
two closed orbits are both defined over $\bb F_q$ with 
stabilizer $\bb G_m.$ In order for there to exist an 
$\bb F_q$-point in the open dense orbit, we need
$$
1-\frac{2}{q-1}>\frac{2}{q-1},
$$
and this is so when $q\ge7.$
\end{example}

\section{About independence of $\ell$}\label{sec-indep}

The coefficients of the expansion of the infinite product
$$
Z(\mathscr X_0,t)=\prod_{i\in\mathbb Z}P_{i,\ell}(\mathscr
X_0,t)^{(-1)^{i+1}}
$$
are rational numbers and are independent of $\ell,$ because
the $c_v(\mathscr X_0)$'s are rational numbers 
independent of $\ell.$ A famous conjecture is that this is 
also true for 
each $P_{i,\ell}(\mathscr X_0,t).$ First we show that the
roots of $P_{i,\ell}(\mathscr X_0,t)$ are Weil $q$-numbers.
Note that $P_{i,\ell}(\mathscr X_0,t)\in\mathbb Q_{\ell}[t].$

\begin{definition}\label{D10.1}
An algebraic number is called a \emph{Weil $q$-number} if
all of its conjugates have the same weight relative to $q,$
and this weight is a rational integer. It is called a
\emph{Weil $q$-integer} if in addition it is an algebraic
integer. A number in $\overline{\mathbb Q}_{\ell}$ is
called a \emph{Weil $q$-number} if it is a Weil $q$-number via
$\iota.$
\end{definition}

For $\alpha\in\overline{\mathbb Q}_{\ell},$ being a Weil
$q$-number or not is independent of $\iota;$ in fact the
images in $\mathbb C$ under various $\iota$'s are conjugate.

For an $\mathbb F_q$-variety $X_0,$ not necessarily smooth or
proper, (\cite{Del2}, 3.3.4) implies all Frobenius
eigenvalues of $H^i_c(X)$ are Weil $q$-integers. The following
lemma generalizes this.

\begin{lemma}\label{L10.2}
For every $\mathbb F_q$-algebraic stack $\mathscr X_0,$ and a
prime number $\ell\ne p,$ the roots of each
$P_{i,\ell}(\mathscr X_0,t)$ are Weil $q$-numbers. In
particular, the coefficients of $P_{i,\ell}(\mathscr X_0,t)$
are algebraic numbers in $\mathbb Q_{\ell}$ (i.e. algebraic
over $\mathbb Q$).
\end{lemma}

\begin{proof}
For an open immersion $j:\mathscr U_0\to\mathscr X_0$ with
complement $i:\mathscr Z_0\to\mathscr X_0,$ we have an exact
sequence
$$
\xymatrix@C=.6cm{
\cdots \ar[r] & H^i_c(\mathscr U) \ar[r] & H^i_c(\mathscr
X) \ar[r] & H^i_c(\mathscr Z) \ar[r] & \cdots,}
$$
thus we may shrink to a non-empty open substack. In particular,
(\ref{L10.2}) holds for algebraic spaces, by (\cite{Knu}, II
6.7) and (\cite{Del2}, 3.3.4).

We may assume $\mathscr X_0$ is smooth and
connected. By Poincar\'e duality, it suffices to show that
the Frobenius eigenvalues of $H^i(\mathscr X)$ are Weil
$q$-numbers, for all $i.$ Take a presentation $X_0\to\mathscr
X_0$ and consider the associated strictly simplicial smooth covering
$X^{\bullet}_0\to\mathscr X_0$ by algebraic spaces. Then
there is a spectral sequence (\cite{LO2}, 10.0.9)
$$
E_1^{rk}=H^k(X^r)\Longrightarrow H^{r+k}(\mathscr X),
$$
and the assertion for $\mathscr X_0$ follows from the
assertion for algebraic spaces.
\end{proof}

\begin{problem}\label{Q10.3}
Is each
$$
P_{i,\ell}(\mathscr X_0,t)=\det(1-Ft,H^i_c(\mathscr
X,\mathbb Q_{\ell}))
$$
a polynomial with coefficients in $\mathbb Q,$ and the
coefficients are independent of $\ell?$
\end{problem}

\begin{subremark}\label{R10.3.1}
(i) Note that, unlike the case for varieties, we cannot expect
the coefficients to be integers (for instance, for $B\mathbb
G_m,$ the coefficients are $1/q^i$).

(ii) (\ref{Q10.3}) is known to be true for smooth proper
varieties (\cite{Del2}, 3.3.9), and (coarse moduli spaces of)
proper smooth algebraic stacks of finite diagonal
(\ref{T7.3}). It remains open for general varieties. Even the
Betti numbers are not known to be independent of $\ell$ for a
general variety. See \cite{Ill}.
\end{subremark}

Let us give positive answer to (\ref{Q10.3}) in some special
cases of algebraic stacks. In $\S7$ we see that it holds for
$BE$ and $BGL_N.$ We can generalize these two cases as
follows.

\begin{lemma}\label{L10.4}
(\ref{Q10.3}) has a positive answer for 

(i) $BA,$ where $A$ is an $\mathbb F_q$-abelian variety;

(ii) $BG_0,$ where $G_0$ is a linear algebraic group over 
$\mathbb F_q.$
\end{lemma}

\begin{proof}
(i) Let $g=\dim A.$ Then $N=H^1(A)$ is a $2g$-dimensional
vector space, with eigenvalues $\alpha_1,\cdots,\alpha
_{2g}$ for the Frobenius action $F,$ and $N$ is pure of
weight 1. Let $a_1,\cdots,a_{2g}$ be a basis for $N$
so that $F$ is upper-triangular
$$
\begin{bmatrix}\alpha_1 & * & * \\
& \ddots & * \\
&& \alpha_{2g}\end{bmatrix}.
$$
Then $H^*(BA)=\text{Sym}^*N[-1]=\overline{\mathbb
Q}_{\ell}[a_1,\cdots,a_{2g}],$ where each $a_i$ sits
in degree 2. In degree $2n,\ H^{2n}(BA)=\overline{
\mathbb Q}_{\ell}\langle a_{i_1}\cdots a_{i_n}|1\le i_1,
\cdots,i_n\le2g\rangle,$ and the eigenvalues are
$\alpha_{i_1}\cdots\alpha_{i_n}.$ By Poincar\'e duality
$$
H_c^{-2n-2g}(BA)=H^{2n}(BA)^{\vee}\otimes\overline{\mathbb
Q}_{\ell}(g)
$$
we see that the eigenvalues of $F$ on $H^{-2g-2n}_c(BA)$ are
$$
q^{-g}\cdot\alpha_{i_1}^{-1}\cdots\alpha_{i_n}^{-1}.
$$
Each factor
$$
P_{-2g-2n}(q^gt)=\prod_{1\le
i_1,\cdots,i_n\le2g}\big(1-(\alpha_{i_1}
\cdots\alpha_{i_n})^{-1}t\big)
$$
stays unchanged if we permute the $\alpha_i$'s arbitrarily,
so the coefficients are symmetric polynomials in the
$\alpha_i^{-1}$'s with integer coefficients, hence are
polynomials in the elementary symmetric functions, which are
coefficients of $\prod_{i=1}^{2g}(t-\alpha_i^{-1}).$ The
polynomial
$$
\prod_{i=1}^{2g}(1-\alpha_it)=\det\big(1-Ft,H^1(A,
\mathbb Q_{\ell})\big)
$$
also has roots $\alpha_i^{-1},$ and this is a polynomial with
integer coefficients, independent of $\ell,$ since $A$ is
smooth and proper. Let $m=\pm q^g$ be leading coefficient of
it. Then
$$
\prod_{i=1}^{2g}(t-\alpha_i^{-1})=\frac{1}{m}\prod_{i=1}
^{2g}(1-\alpha_it).
$$
This verifies (\ref{Q10.3}) for $BA.$

(ii) Let $d=\dim G_0.$ For every $k\ge0,\ H^{2k}(BG)$ is pure
of weight $2k$ (\cite{Del3}, 9.1.4), hence by Poincar\'e
duality, $H^{-2d-2k}_c(BG)$ is pure of weight $-2d-2k.$
The entire function
$$
\frac{1}{Z(BG_0,t)}=\prod_{k\ge0}P_{-2d-2k}(BG_0,t)\in
\mathbb Q[[t]]
$$
is independent of $\ell,$ and invariant under the action of
$\text{Gal}(\mathbb Q)$ on the coefficients of the Taylor
expansion. Therefore the roots of $P_{-2d-2k}(BG_0,t)$ can be
described as
$$
\text{``zeros of }\frac{1}{Z(BG_0,t)}\text{ that have
weight }2d+2k\text{ relative to }q",
$$
which is a description independent of $\ell,$ and these roots
(which are algebraic numbers) are permuted under
$\text{Gal}(\mathbb Q).$ Hence $P_{-2d-2k}(BG_0,t)$ has
rational coefficients.
\end{proof}

The following proposition generalizes both (\ref{T7.3}) and
(\ref{L10.4}ii).

\begin{proposition}\label{P10.6}
Let $X_0$ be the coarse moduli space of a proper smooth
$\mathbb F_q$-algebraic stack of finite diagonal, and let
$G_0$ be a linear $\mathbb F_q$-algebraic group that acts on
$X_0,$ and let $\mathscr X_0$ be a form of the quotient stack
$[X_0/G_0].$ Then (\ref{Q10.3}) is verified for $\mathscr X_0.$
\end{proposition}

\begin{proof}
It suffices to show that $H^n_c(\mathscr X)$ is pure of 
weight $n,$ for every $n.$ To show this, we can make a 
finite extension of the base field $\mathbb F_q,$ so we 
may assume $\mathscr X_0=[X_0/G_0].$ Let
$$
\xymatrix@C=.5cm{
\mathscr X_0 \ar[r]^-f & BG_0 \ar[r]^-h & B\pi_0(G_0)}
$$
be the natural maps. 

Let $d=\dim G_0.$ Consider the spectral sequence
$$
H^{-2d-2r}_c(BG,R^kf_!\overline{\mathbb Q}_{\ell})
\Longrightarrow H^{-2d-2r+k}_c(\mathscr X).
$$
The $E_2$-terms can be computed from the degenerate Leray
spectral sequence for $h:$
$$
H^{-2d-2r}_c(BG,R^kf_!\overline{\mathbb Q}_{\ell})\simeq
H^0_c(B\pi_0(G),R^{-2d-2r}h_!R^kf_!\overline{\mathbb
Q}_{\ell}).
$$
The restriction of $R^{-2d-2r}h_!R^kf_!\overline{\mathbb
Q}_{\ell}$ along the natural projection $\text{Spec }\mathbb 
F_q\to B\pi_0(G_0)$ is isomorphic to the Galois module 
$H^{-2d-2r}_c(BG^0,R^kf_!\overline{\mathbb Q}_{\ell}),$ 
and since $G^0_0$ is connected, $(R^kf_!\overline{\mathbb 
Q}_{\ell})|_{BG_0^0}$ is the inverse image of some sheaf 
$V_k$ via the structural map $BG^0_0\to\text{Spec }\mathbb 
F_q.$ By base change, we see that the sheaf $V_k,$ regarded 
as a $\text{Gal}(\mathbb F_q)$-module, is $H^k(X).$ 
By projection formula we have
$$
H^{-2d-2r}_c(BG^0,R^kf_!\overline{\mathbb
Q}_{\ell})\simeq H^{-2d-2r}_c(BG^0)\otimes H^k(X)
$$
as representations of $\text{Gal}(\mathbb F_q),$ and by
(\ref{T7.3}), the right hand side is pure of weight
$-2d-2r+k.$ By (\ref{L4.8}),
$H^{-2d-2r}_c(BG,R^kf_!\overline{\mathbb
Q}_{\ell})$ is also pure of weight $-2d-2r+k,$ therefore
$H^n_c(\mathscr X)$ is pure of weight $n,$ for every $n.$
\end{proof}

\begin{blank}\label{10.7}
Finally, let us consider the following much weaker version of
independence of $\ell.$ For $\mathscr X_0$ and $i\in\mathbb
Z,$ let $\Psi(\mathscr X_0,i)$ be the following property: the
Frobenius eigenvalues of $H^i_c(\mathscr X,\overline{\mathbb
Q}_{\ell}),$ counted with multiplicity, for all $\ell\ne p,$
are contained in a finite set of algebraic numbers with
multiplicities assigned, and this set together with the
assignment of multiplicity, depends only on $\mathscr X_0$ and
$i.$ In particular it is independent of $\ell.$ In other
words, there is a finite decomposition of the set of all prime
numbers $\ell\ne p$ into disjoint union of some
subsets, such that the Frobenius eigenvalues of
$H^i_c(\mathscr X,\overline{\mathbb Q}_{\ell})$ depends only
on the subset that $\ell$ belongs to. If this property holds,
we also denote such a finite set of algebraic numbers (which
is not unique) by $\Psi(\mathscr X_0,i),$ if there is no
confusion.

\begin{subproposition}\label{P10.7}
The property $\Psi(\mathscr X_0,i)$ holds for every $\mathscr
X_0$ and $i.$
\end{subproposition}

\begin{proof}
If $\mathscr U_0$ is an open substack of $\mathscr X_0$
with complement $\mathscr Z_0,$ and properties
$\Psi(\mathscr U_0,i)$ and $\Psi(\mathscr Z_0,i)$ hold,
then $\Psi(\mathscr X_0,i)$ also holds, and the finite set
$\Psi(\mathscr X_0,i)$ a subset of $\Psi(\mathscr
U_0,i)\cup\Psi(\mathscr Z_0,i).$

Firstly we prove this for schemes $X_0.$ By shrinking $X_0$
we can assume it is a connected smooth variety. By Poincar\'e
duality it suffices to prove the similar statement
$\Psi^*(X_0,i)$ for ordinary cohomology, i.e. with $H^i_c$
replaced by $H^i,$ for all $i.$ This follows from \cite{deJ}
and (\cite{Del2}, 3.3.9). Therefore it also holds for all
algebraic spaces.

For a general algebraic stack $\mathscr X_0,$ by shrinking
it we can assume it is connected smooth. By Poincar\'e
duality, it suffices to prove $\Psi^*(\mathscr X_0,i)$ for
all $i.$ This can be done by taking a hypercover by
simplicial algebraic spaces, and considering the associated
spectral sequence.
\end{proof}
\end{blank}

\end{document}